# Representations of Sklyanin-type algebras, the elliptic modular triple and the one-body quantum Ruijsenaars-van Diejen model

A. Mostovskii A. Zotov

*Steklov Mathematical Institute of the Russian Academy of Sciences,*
*8 Gubkina Street, 119991 Moscow, Russia*

E-mails: mostovskii.am21@physics.msu.ru, zotov@mi-ras.ru

**Abstract**

In a central quotient of the four-dimensional Sklyanin algebra, we derive closed quadratic relations generated by quadratic combinations of the original algebra. Then we construct a one-parameter deformation of the induced representation by elliptic difference operators, and obtain representations of the original Sklyanin algebra and its $\mathcal{B}C_1$ generalization by specialization. A linear combination of the deformed generators gives the eight-parameter one-body Ruijsenaars-van Diejen Hamiltonian. Next, we propose an elliptic modular triple by extending and unifying Spiridonov's two elliptic modular doubles. In particular, this construction yields three pairwise commuting Ruijsenaars-van Diejen Hamiltonians with parameters transformed as in the modular triple, extending the pair of Rains and Ruijsenaars. In the three-dimensional Sklyanin algebra modulo its central cubic, the same construction yields quadratic elements satisfying the Sklyanin relations with certain transformed parameters.

## Contents

# 1 Introduction

We begin by recalling a number of known constructions and then turn to the aims and results of this paper.

## 1.1 An overview

**Sklyanin-type algebras.** The Sklyanin algebra provides an algebraic framework for elliptic quantum integrable systems studied by the quantum inverse scattering method [37, 38]. The original and best-known version is defined as follows. The Sklyanin algebra $\mathcal{A}$ is the quotient of the free associative algebra on four generators $S_0, S_1, S_2, S_3$ by the two-sided ideal generated by the following quadratic relations. For every cyclic permutation $(i, j, k)$ of $(1, 2, 3)$, we have

$$\mathcal{A}: \quad \begin{aligned} &[S_i, S_j]_- = \imath [S_0, S_k]_+ \,, \\ &[S_0, S_k]_- = \imath \frac{J_j - J_i}{J_k} [S_i, S_j]_+ \,, \end{aligned} \tag{1.1}$$

where $[X, Y]_\pm = XY \pm YX$ and $J_1, J_2, J_3 \neq 0$ are $\mathbb{C}$-valued parameters. These relations can be deduced from $R$-matrix quantization using the quantum exchange relation $R_{12}(z_1 - z_2)L_1(z_1)L_2(z_2) = L_2(z_2)L_1(z_1)R_{12}(z_1 - z_2)$ with Baxter's eight-vertex $R$-matrix. The Sklyanin algebra is an $\mathbb{N}$-graded $\mathbb{C}$-algebra. Thus,

$$\mathcal{A} = \bigoplus_{n \geq 0} \mathcal{A}_n \,, \tag{1.2}$$

where $\mathcal{A}_0 = \mathbb{C}$, $\mathcal{A}_1 = \mathrm{span}\{S_0, S_i, S_j, S_k\}$ and the component $\mathcal{A}_2$ is spanned by the quadratic monomials $S_\alpha S_\beta$, $\alpha, \beta = 0, ..., 3$.

The next quadratic algebra we study is the generalization of (1.1) proposed in [19], which includes terms linear in $S_i$:

$$\mathcal{BC}_1: \quad \begin{aligned} &[\mathbb{S}_i, \mathbb{S}_j]_- = \imath[\mathbb{S}_0, \mathbb{S}_k]_+ \,, \\ &[\mathbb{S}_0, \mathbb{S}_k]_- = -\imath\frac{K_j - K_i}{K_k}[\mathbb{S}_i, \mathbb{S}_j]_+ + \frac{2\imath}{K_k}\left(\widetilde{\nu}_i\rho_i\mathbb{S}_j - \widetilde{\nu}_j\rho_j\mathbb{S}_i\right). \end{aligned} \tag{1.3}$$

It has three additional constants $\widetilde{\nu}_1, \widetilde{\nu}_2, \widetilde{\nu}_3 \in \mathbb{C}$. We call it the $\mathcal{BC}_1$ Sklyanin algebra since it is related to the Calogero-Inozemtsev and Ruijsenaars-van Diejen models of type $\mathrm{BC}_1$ (see [19, 24]). The algebra (1.3) also arises from $R$-matrix quantization (of the Zhukovsky-Volterra gyrostat), this time through the reflection equation [39] $R_{12}(z_1-z_2)L_1(z_1)R_{12}(z_1+z_2)L_2(z_2) = L_2(z_2)R_{12}(z_1+z_2)L_1(z_1)R_{12}(z_1-z_2)$. We recall this construction in Appendix B. When $\widetilde{\nu}_j = 0$, (1.3) reduces to the original algebra (1.1) up to a redefinition of the parameters $K_i$.[1]

Finally, we use a more general quadratic algebra with additional linear terms:

$$\mathcal{A}(\mathcal{A}_2): \quad \begin{aligned} &[\mathcal{Z}_i, \mathcal{Z}_j]_- = \imath\,[\mathcal{Z}_0, \mathcal{Z}_k]_+ + \imath u_0\mu_k\mathcal{Z}_k - \imath u_k\lambda_k\mathcal{Z}_0 \,, \\ &[\mathcal{Z}_0, \mathcal{Z}_k]_- = -\imath\,\frac{K_j - K_i}{K_k}\,[\mathcal{Z}_i, \mathcal{Z}_j]_+ + 2\imath\frac{u_i}{K_k}\,\rho_i\mu_j\mathcal{Z}_j - 2\imath\frac{u_j}{K_k}\,\rho_j\mu_i\mathcal{Z}_i \,. \end{aligned} \tag{1.4}$$

Here $\mathcal{Z}_0, \mathcal{Z}_1, \mathcal{Z}_2, \mathcal{Z}_3$ are the generators, and $u_\alpha, \mu_k, \lambda_k, \rho_k$ are additional parameters subject to constraints. We denote this algebra by $\mathcal{A}(\mathcal{A}_2)$ since it arises as a Sklyanin-type algebra generated by certain quadratic elements of the original algebra $\mathcal{A}$ (1.1). This algebra plays an auxiliary role. We consider special values of its parameters and study its degenerations to $\mathcal{A}$ (1.1) or $\mathcal{BC}_1$ (1.3). The algebra (1.4) can be described as a specialization of a more general quadratic algebra recently proposed by A. Mikhailov and T. Skrypnyk in [23].

**Representations and their relation to the Ruijsenaars-van Diejen Hamiltonian.** In this paper we mainly discuss representations of Sklyanin algebras in terms of difference operators. In his original work [38], Sklyanin proposed the following representation on the space of meromorphic functions:

$$\hat{S}_\alpha(z, \hbar, \eta|\tau) = \imath^{\delta_{\alpha,2}}\frac{\theta_{\alpha+1}(\hbar)\theta_{\alpha+1}(-2\eta)}{\theta_1'(0)}\left(\varphi_\alpha(2z, -2\eta + \omega_\alpha)e^{\hbar\partial_z} + \varphi_\alpha(-2z, -2\eta + \omega_\alpha)e^{-\hbar\partial_z}\right) \tag{1.5}$$

for $\alpha \in \{0, 1, 2, 3\}$, where $\theta_a(x)$ are the theta functions (A.1) and $\varphi_\alpha(x, y+\omega_\alpha)$ are the functions defined in (A.25). It is well known [14] that the generator with index 0 is the Hamiltonian of the one-body quantum elliptic Ruijsenaars-Schneider model [35]:

$$\hat{S}_0(z, \hbar, \eta|\tau) = \hat{H}^{\mathrm{RS}} \,. \tag{1.6}$$

For fixed $\hbar$ and $\tau$, this expression depends on a single parameter $\eta$. It plays the role of the spin parameter in the XYZ model described by the Sklyanin algebra and of the coupling constant in the Ruijsenaars-Schneider model. The one-body Hamiltonian (1.6) corresponds to the root system of type $A_1$. Its generalization to the $\mathrm{BC}_1$ root system is known as the one-body Ruijsenaars-van Diejen model

[1]When $\widetilde{\nu}_j = 0$, (1.1) and (1.3) have the same form except for the sign on the right-hand side of the second line. For this reason, the parameters $J_i$ and $K_i$ are related by a change of variables. See Section 2 for details.

[9], and its Hamiltonian as the elliptic $q$-difference Heun operator [36]. The integrability of the general $\mathrm{BC}_n$ model was proved in part in [9] and in full in [17]. In this paper we use the form of the one-body Ruijsenaars-van Diejen Hamiltonian proposed in [4]:

$$\begin{gathered}\hat{H}^{8\mathrm{vD}}(z|\tau;\hbar,\nu,\bar{\nu},\eta,\xi)=\bar{v}(z,-2\eta)v(z+\hbar,-2\xi)e^{2\hbar\partial_z}+\bar{v}(-z,-2\eta)v(-z+\hbar,-2\xi)e^{-2\hbar\partial_z}+\\+\bar{v}(-z,-2\eta)v(z-\hbar,-2\eta)+\bar{v}(z,-2\eta)v(-z-\hbar,-2\eta)\,.\end{gathered}\tag{1.7}$$

The function $v(z,\eta)$ is defined in (A.32) as a linear combination of the functions $\varphi_a(2z,u+\omega_a)$:

$$v(z,\eta)\equiv v(z,\eta|\tau;\nu_0,\nu_1,\nu_2,\nu_3)=\sum_{a=0}^{3}\nu_a\varphi_a(2z,\eta+\omega_a)\,,\tag{1.8}$$

and $\bar{v}(z,\eta)=v(z,\eta|\tau;\bar{\nu})$ is the same function with the four parameters $\bar{\nu}_a$ in place of $\nu_a$. The number 8 in (1.7) denotes the number of independent constants[2]. When $\xi=\eta$, (1.7) factors as

$$\hat{H}^{7\mathrm{vD}}=\Delta^{\mathrm{R}}_{\bar{\nu}}\Delta^{\mathrm{R}}_{\nu}\,,\qquad\Delta^{\mathrm{R}}_{\nu}=v(z,-2\eta)e^{\hbar\partial_z}+v(-z,-2\eta)e^{-\hbar\partial_z}\,,\tag{1.9}$$

where $\Delta^{\mathrm{R}}_{\nu}=\Delta^{\mathrm{R}}(z|\tau;\hbar,\nu,\eta)$ is the Rains $\Delta$-operator [31, 34]. A similar relation was observed at the classical level in [24]. In that paper, the classical analogue of $\Delta^{\mathrm{R}}_{\nu}$ was identified with the classical counterpart of the generator $\mathbb{S}_0$ of the $\mathcal{BC}_1$ algebra related to the Zhukovsky-Volterra gyrostat. Another important result relating the Sklyanin algebra and the Ruijsenaars-van Diejen Hamiltonian $\hat{H}^{8\mathrm{vD}}$ was formulated by Rains and Ruijsenaars in [32]. It states that $\hat{H}^{8\mathrm{vD}}$ can be expressed as a linear combination of ordered products of two Sklyanin generators:

$$\hat{H}^{8\mathrm{vD}}=\sum_{\alpha,\beta=0}^{3}C_{\alpha\beta}\,\hat{S}_\alpha(z,\hbar,\frac{\xi+\eta}{2}\mid\tau)\hat{S}_\beta(z,\hbar,\frac{\xi+\eta}{2}\mid\tau)\,,\tag{1.10}$$

for some constants $C_{\alpha\beta}$. Their statement is an existence result.

Some applications of the Sklyanin algebras and Ruijsenaars-van Diejen model in theoretical physics can be found in [33, 21]. Its relation to Painlevé equations was discussed in [27].

**Additional quantum Hamiltonian.** Another important result of [32] is the construction of a quantum Hamiltonian commuting with the Ruijsenaars-van Diejen operator $\hat{H}^{8\mathrm{vD}}$ (1.7). Namely, define

$$\hat{H}^{8\mathrm{vD}}_1:=\hat{H}^{8\mathrm{vD}}\left(z|2\hbar;\frac{\tau}{2},\widetilde{\nu}^1,\bar{\widetilde{\nu}}^1,\eta_1,\xi_1\right)\,.\tag{1.11}$$

Then

$$[\hat{H}^{8\mathrm{vD}},\hat{H}^{8\mathrm{vD}}_1]=0\,,\qquad\eta_1+\xi_1-\eta-\xi\in\frac{1}{4}\mathbb{Z}\,.\tag{1.12}$$

The two sets of coupling constants $\widetilde{\nu}^1$ and $\bar{\widetilde{\nu}}^1$ can be chosen independently of $\nu$ and $\bar{\nu}$. In [32], commutativity was proved using the quasi-periodicity of the coefficient functions.

[2]The rescaling $(\nu,\bar{\nu})\mapsto(a\nu,a^{-1}\bar{\nu})$ leaves the operator invariant. An independent rescaling multiplies the entire operator by a constant. Thus eight parameters remain when the overall normalization is disregarded.

**Elliptic modular double.** In [42], V. Spiridonov proposed two elliptic analogues of Faddeev's modular double for $U_q(sl_2)$ algebra [11]. Consider the representation $\hat{S}_\alpha(z,\hbar,\eta|\tau)$ in (1.5) of the Sklyanin algebra (1.1). Set

$$\hat{\tilde{S}}^1_\alpha = \hat{S}_\alpha(z,\frac{\tau}{2},\eta|2\hbar)\,, \qquad \hat{\tilde{S}}^2_\alpha = \hat{S}_\alpha(\frac{z}{2\hbar},\frac{1}{4\hbar},\frac{\eta}{2\hbar}|\frac{\tau}{2\hbar}) \tag{1.13}$$

and introduce a pair of matrices

$$\varepsilon_1 = \begin{pmatrix} 1 & -1 & -1 & 1 \\ -1 & 1 & 1 & -1 \\ -1 & 1 & 1 & -1 \\ 1 & -1 & -1 & 1 \end{pmatrix}, \qquad \varepsilon_2 = \begin{pmatrix} 1 & 1 & -1 & -1 \\ 1 & 1 & -1 & -1 \\ -1 & -1 & 1 & 1 \\ -1 & -1 & 1 & 1 \end{pmatrix}. \tag{1.14}$$

Both sets of operators $\hat{\tilde{S}}^1_\alpha$ and $\hat{\tilde{S}}^2_\alpha$ satisfy the Sklyanin relations (1.1) since they are obtained from $\hat{S}_\alpha$ with different parameters. Spiridonov [42] showed that the following cross-relations hold:

$$\hat{S}_\alpha\hat{\tilde{S}}^1_\beta = \varepsilon_1^{\alpha\beta}\,\hat{\tilde{S}}^1_\beta\hat{S}_\alpha\,, \qquad \hat{S}_\alpha\hat{\tilde{S}}^2_\beta = \varepsilon_2^{\alpha\beta}\,\hat{\tilde{S}}^2_\beta\hat{S}_\alpha\,. \tag{1.15}$$

Then one can define an elliptic modular double as a pair of Sklyanin algebras (1.1) with the cross-relations $S_\alpha\tilde{S}^1_\beta = \varepsilon_1^{\alpha\beta}\,\tilde{S}^1_\beta S_\alpha$ (the first version of the elliptic modular double) or $S_\alpha\tilde{S}^2_\beta = \varepsilon_2^{\alpha\beta}\,\tilde{S}^2_\beta S_\alpha$ (the second version of the elliptic modular double). Formulae (1.13) give representations of these two doubles by difference operators. The first version does not have a well-defined trigonometric limit $\mathrm{Im}(\tau)\to+\infty$, but the second one in the trigonometric limit reduces to known representations of Faddeev's modular double [11, 15].

**Trigonometric background.** A standard trigonometric degeneration of the Sklyanin algebra, followed by a suitable change of generators and normalization of a central element, yields

$$\mathcal{U}_q(sl_2):\quad \begin{array}{c} AB = qBA, \qquad BD = qDB, \qquad DC = qCD, \qquad CA = qAC, \\ [B,C] = \dfrac{A^2-D^2}{q-q^{-1}}\,, \qquad AD = DA = 1\,, \end{array} \tag{1.16}$$

with central element

$$\Omega^{\mathcal{U}_q} = q^{-1}A^2 + qD^2 + (q-q^{-1})^2BC\,. \tag{1.17}$$

Other trigonometric degenerations were obtained in [10]. Algebra (1.16) can also be derived from the $RLL$ relation with the standard trigonometric $XXZ$ $R$-matrix [12]. Recall that the $\mathcal{B}C_1$ algebra is defined by the reflection equation with elliptic $R$-matrix (B.1) (see Appendix B) and, as was shown in [46], for the $XXZ$ $R$-matrix, the reflection equation gives the Askey-Wilson-type algebra [47] with three generators $H_1$, $H_2$, $H_3$ and relation between them

$$\mathcal{AW}:\qquad q^{-1}H_\alpha H_\beta - qH_\beta H_\alpha = g_\gamma H_\gamma + h_\gamma\,, \tag{1.18}$$

where $\{\alpha,\beta,\gamma\}$ stands for any cyclic permutation of $\{1,2,3\}$ and $g_\gamma$, $h_\gamma$ are some $\mathbb{C}$-valued constants. Algebra (1.18) has been realized by quadratic elements of $\mathcal{U}_q(sl_2)$

$$\begin{aligned} H_1 &= \big(\chi_1B+\chi_2C+\chi_3A\big)A, \\ H_2 &= D\big(\kappa_1B+\kappa_2C+\kappa_3D\big), \\ H_3 &= (q^2-q^{-2})\left(\frac{\chi_1\kappa_2A^2+\chi_2\kappa_1D^2}{(q-q^{-1})^2} - q\chi_2\kappa_2C^2 - \chi_2\kappa_3DC - \chi_3\kappa_2CA\right). \end{aligned} \tag{1.19}$$

The structure constants $g_\gamma$ and $h_\gamma$ can be expressed in terms of coefficients $\chi_i$ and $\kappa_i$ and the value of the central element (1.17), see [45], [46] for details. This construction provides a trigonometric example of the quadratic realizations studied below for elliptic Sklyanin-type algebras. Furthermore, the Askey-Wilson algebra can be embedded in the $BC_1$ DAHA [18]. This fact partially explains the name of algebra (1.3).

## 1.2 Purpose of the paper

**Representations of $\mathcal{A}(\mathcal{A}_2)$.** We begin by computing the commutation relations between quadratic elements $X_\alpha(U)$ of $\mathcal{A}$ of the form

$$X_\alpha(U) = US_\alpha \in \mathcal{A}_2\,, \quad Y_\alpha(U) = S_\alpha U \in \mathcal{A}_2\,, \qquad U = \sum_{\beta=0}^{3} u_\beta S_\beta \in \mathcal{A}_1\,, \qquad u_\beta \in \mathbb{C}\,. \tag{1.20}$$

We prove that the normalized elements $Y_\alpha(U)$, denoted here by $Z_\alpha$, satisfy the relations (1.4) of $\mathcal{A}(\mathcal{A}_2)$ with parameters depending on the values $\delta_1, \delta_2 \in \mathbb{C}$ of the quadratic Casimir operators $\Omega_1 = S_0^2 + S_1^2 + S_2^2 + S_3^2$ and $\Omega_2 = J_1 S_1^2 + J_2 S_2^2 + J_3 S_3^2$ of the original algebra $\mathcal{A}$. In this way we obtain a realization of the algebra (1.4) by the classes of the quadratic elements $Z_\alpha$ in the central quotient $\mathcal{A}/(\Omega_1 - \delta_1, \Omega_2 - \delta_2)$. This realization allows us to construct representations of $\mathcal{BC}_1$ and $\mathcal{A}$ since these algebras are obtained from $\mathcal{A}(\mathcal{A}_2)$ in particular cases. Indeed, we have $Z_\alpha = Z_\alpha(S)$. Then the difference operators $\hat{Z}_\alpha = Z_\alpha(\hat{S})$ satisfy the same relations.

However, the Casimir operators $\Omega_1^{\mathcal{Z}}$ and $\Omega_2^{\mathcal{Z}}$ of the algebra $\mathcal{A}(\mathcal{A}_2)$ both vanish in the realization described above: $\Omega_1^{\mathcal{Z}}\Big|_{\mathcal{Z}_\alpha = Z_\alpha} = \Omega_2^{\mathcal{Z}}\Big|_{\mathcal{Z}_\alpha = Z_\alpha} = 0$. Therefore, their values in the representation $\hat{Z}_\alpha = Z_\alpha(\hat{S})$ are also zero: $\delta_1^{\mathcal{Z}} = \delta_2^{\mathcal{Z}} = 0$. For this reason, we propose a one-parameter generalization of the representation $\hat{Z}_\alpha = Z_\alpha(\hat{S})$. We denote the additional parameter by $\xi \in \mathbb{C}$. As a result, we obtain the following representation of the algebra $\mathcal{A}(\mathcal{A}_2)$ defined by (1.4):

$$\begin{gathered}\hat{Z}_\alpha^\xi = \imath^{\delta_{\alpha,2}} \frac{\theta_1(\hbar)\theta_{\alpha+1}(2\hbar)\theta_{\alpha+1}(-2\eta)}{\theta_1'(0)\theta_1(2\hbar)} \Big( \varphi_\alpha(2z, -2\eta + \omega_\alpha) v(z+\hbar, -2\xi) e^{2\hbar\partial_z} + \\ + \varphi_\alpha(-2z, -2\eta + \omega_\alpha) v(z-\hbar, -2\eta) + \big(z \to -z\big)\Big)\,, \qquad \alpha \in \{0,1,2,3\}\,.\end{gathered} \tag{1.21}$$

The values $\delta_1^{\mathcal{Z}}, \delta_2^{\mathcal{Z}}$ of the Casimir elements $\Omega_1^{\mathcal{Z}}$ and $\Omega_2^{\mathcal{Z}}$ in this representation are given by

$$\delta_1^{\mathcal{Z}} = \delta_1(2\hbar, \eta) \frac{\theta_1^2(\hbar)}{\theta_1^2(2\hbar)} \sum_{\alpha=0}^{3} \nu_\alpha^2 \Big( \wp(2\xi + \omega_\alpha) - \wp(2\eta + \omega_\alpha) \Big)\,, \quad \nu_\alpha = \imath^{\delta_{\alpha,2}} \frac{\theta_{\alpha+1}(\hbar)\theta_{\alpha+1}(-2\eta)}{\theta_1'(0)}\, u_\alpha\,, \tag{1.22}$$

$$\delta_2^{\mathcal{Z}} = -64\, \delta_2(2\hbar, \eta) \frac{\theta_1^{10}(\hbar)\theta_1^2(4\hbar)}{\theta_1^{12}(2\hbar)} \sum_{\alpha=0}^{3} \nu_\alpha^2 \Big( \wp(2\xi + \omega_\alpha) - \wp(2\eta + \omega_\alpha) \Big)\,. \tag{1.23}$$

Since this $\xi$-deformed representation is not induced by the map $\hat{Z}_\alpha = Z_\alpha(\hat{S})$, we prove (1.21)–(1.23) directly. Explicit dependence on $\xi$ of the Casimirs values in (1.22)–(1.23) means that different values of $\xi$ provide non-equivalent representations. When $\xi = \eta$, we recover the induced representation $\hat{Z}_\alpha = Z_\alpha(\hat{S})$.

**Representations of $\mathcal{B}C_1$.** Given a representation of $\mathcal{A}(\mathcal{A}_2)$, we may consider particular cases corresponding to $\mathcal{B}C_1$ or $\mathcal{A}$. For example, for the $\mathcal{B}C_1$ algebra (1.3) we obtain the following representation:

$$\hat{\mathbb{S}}_0^\xi = \frac{\theta_1(\hbar)\theta_1(2\hbar)}{\theta_1'(0)} \Big(\phi(2z,2\hbar)v(z+\hbar,-2\xi)e^{2\hbar\partial_z} + \phi(-2z,2\hbar)v(-z+\hbar,-2\xi)e^{-2\hbar\partial_z}\Big)\,, \tag{1.24}$$

$$\hat{\mathbb{S}}_k^\xi = \imath^{\delta_{k,2}}\frac{\theta_1(\hbar)\theta_{k+1}^2(2\hbar)}{\theta_1'(0)\theta_1(2\hbar)} \Big(\varphi_k(2z,2\hbar+\omega_k)v(z+\hbar,-2\xi)e^{2\hbar\partial_z}+ \tag{1.25}$$

$$+\varphi_k(-2z,2\hbar+\omega_k)v(z-\hbar,2\hbar) + \big(z\to -z\big)\Big)\,, \qquad k=1,2,3\,.$$

The Casimir elements $\Omega_1^{\mathcal{B}C_1}$ and $\Omega_2^{\mathcal{B}C_1}$ of the $\mathcal{B}C_1$ algebra (1.3) are

$$\Omega_1^{\mathcal{B}C_1} = \mathbb{S}_0^2 + \sum_{k=1}^{3}\mathbb{S}_k^2\,, \qquad \Omega_2^{\mathcal{B}C_1} = \sum_{k=1}^{3}\Big(K_k\big(K_k-K_i-K_j\big)\mathbb{S}_k^2 + 2\widetilde{\nu}_k\rho_k K_k\mathbb{S}_k\Big)\,. \tag{1.26}$$

In this representation, their values $\delta_1^{\mathcal{B}C_1}$ and $\delta_2^{\mathcal{B}C_1}$ are given by

$$\delta_1^{\mathcal{B}C_1} = 4u_0^2\theta_1^4(\hbar)\,, \tag{1.27}$$

$$\delta_2^{\mathcal{B}C_1} = 256\,\frac{\theta_1^{10}(\hbar)\theta_1^2(4\hbar)}{\theta_1'(0)^2\theta_1^{10}(2\hbar)}\sum_{\alpha=0}^{3}(-1)^{\delta_{\alpha,2}}u_\alpha^2\theta_{\alpha+1}^2(\hbar)\theta_{\alpha+1}^2(2\hbar)\Big(\wp(2\xi+\omega_\alpha)-\wp(2\hbar+\omega_\alpha)\Big)\,. \tag{1.28}$$

The constants $\nu_\alpha$ in the definition of the function $v$ in (1.24)–(1.25) are given by (1.22) with $\eta=-\hbar$: $\nu_\alpha = i^{\delta_{\alpha,2}}\theta_{\alpha+1}(\hbar)\theta_{\alpha+1}(2\hbar)u_\alpha/\theta_1'(0)$. At the end of the paper, we consider the classical limit and recover our results from [24] for the classical quadratic $\mathrm{BC}_1$ Sklyanin algebra related to the Zhukovsky-Volterra gyrostat.

The property, which consists in the fact that some quadratic combination of generators of algebra form a closed algebra of the same type, seems very intriguing, it would be interesting to understand its geometric meaning and to find another algebras which possesses it. Trigonometric (and rational) degenerations of these results will be studied in a separate publication.

**Elliptic modular triple.** The representations (1.13) of the elliptic modular doubles are formulated in terms of a single parameter $\eta$. For $\hat{S}_\alpha = \hat{S}_\alpha(z,\hbar,\eta|\tau)$, we consider generalizations of $\hat{\tilde{S}}{}^1_\alpha$ and $\hat{\tilde{S}}{}^2_\alpha$ obtained by introducing the parameters $\eta_1,\eta_2\in\mathbb{C}$:

$$\hat{\tilde{S}}{}^1_\alpha = \hat{S}_\alpha(z,\frac{\tau}{2},\eta_1|2\hbar)\,, \qquad \hat{\tilde{S}}{}^2_\alpha = \hat{S}_\alpha(\frac{z}{2\hbar},\frac{1}{4\hbar},\frac{\eta_2}{2\hbar}|\frac{\tau}{2\hbar})\,. \tag{1.29}$$

First, we show that when $\eta_1 = \eta + \frac{k}{4}$, $k\in\mathbb{Z}$, the relations of both doubles remain the same except for the first cross-relation in (1.15), which acquires an additional factor $(-1)^k$:

$$\hat{S}_\alpha\hat{\tilde{S}}{}^1_\beta = (-1)^k\varepsilon_1^{\alpha\beta}\,\hat{\tilde{S}}{}^1_\beta\hat{S}_\alpha\,. \tag{1.30}$$

Second, we prove that when $\eta_2 = \frac{N\hbar}{2}$, $N\in\mathbb{Z}$, there are also closed cross-commutation relations between $\hat{\tilde{S}}{}^1_\alpha$ and $\hat{\tilde{S}}{}^2_\beta$:

$$\hat{\tilde{S}}{}^1_\alpha\hat{\tilde{S}}{}^2_\beta = (-1)^N\varepsilon^{\alpha\beta}\,\hat{\tilde{S}}{}^2_\beta\hat{\tilde{S}}{}^1_\alpha \tag{1.31}$$

with

$$\varepsilon = \begin{pmatrix} 1 & -1 & -1 & 1 \\ 1 & -1 & -1 & 1 \\ -1 & 1 & 1 & -1 \\ -1 & 1 & 1 & -1 \end{pmatrix} . \tag{1.32}$$

This allows us to combine the two modular doubles into an elliptic modular triple. It consists of three copies of the Sklyanin algebra (1.1) with cross-relations

$$S_\alpha \overset{1}{\widetilde{S}}_\beta = (-1)^k \varepsilon_1^{\alpha\beta} \overset{1}{\widetilde{S}}_\beta S_\alpha \,, \qquad S_\alpha \overset{2}{\widetilde{S}}_\beta = \varepsilon_2^{\alpha\beta} \overset{2}{\widetilde{S}}_\beta S_\alpha \,, \qquad \overset{1}{\widetilde{S}}_\alpha \overset{2}{\widetilde{S}}_\beta = (-1)^N \varepsilon^{\alpha\beta} \overset{2}{\widetilde{S}}_\beta \overset{1}{\widetilde{S}}_\alpha \tag{1.33}$$

for integers $k, N \in \mathbb{Z}$ and the three matrices in (1.14) and (1.32):

$$\varepsilon_1 = \begin{pmatrix} 1 & -1 & -1 & 1 \\ -1 & 1 & 1 & -1 \\ -1 & 1 & 1 & -1 \\ 1 & -1 & -1 & 1 \end{pmatrix} , \; \varepsilon_2 = \begin{pmatrix} 1 & 1 & -1 & -1 \\ 1 & 1 & -1 & -1 \\ -1 & -1 & 1 & 1 \\ -1 & -1 & 1 & 1 \end{pmatrix} , \; \varepsilon = \begin{pmatrix} 1 & -1 & -1 & 1 \\ 1 & -1 & -1 & 1 \\ -1 & 1 & 1 & -1 \\ -1 & 1 & 1 & -1 \end{pmatrix} . \tag{1.34}$$

A representation of these relations is given by $\hat{S}_\alpha = \hat{S}_\alpha(z, \hbar, \eta|\tau)$ and (1.29) with $\eta_1 = \eta + \frac{k}{4}$, $\eta_2 = \frac{N\hbar}{2}$.

The term modular triple was already used in paper [26], where it means extension of the modular double for $q$-Virasoro algebra to $\mathrm{SL}(3, \mathbb{Z})$ symmetry. Transformations of parameters in (1.29) can also be related with $\mathrm{SL}(3, \mathbb{Z})$.

$\mathrm{BC}_1$ **Ruijsenaars-van Diejen model.** For the Hamiltonian (1.7) we first show that it is written as a linear combination of representations of generators $\hat{Z}^\xi_\alpha$ (1.21) of the algebra $\mathcal{A}(\mathcal{A}_2)$ (1.4)

$$\hat{H}^{\text{8vD}}(z|\tau; \hbar, \nu, \bar{\nu}, \eta, \xi) = \sum_{\alpha=0}^{3} \bar{u}_\alpha \hat{Z}^\xi_\alpha \tag{1.35}$$

with the following identification of parameters:

$$\nu_\alpha = \imath^{\delta_{\alpha,2}} \frac{\theta_{\alpha+1}(\hbar)\theta_{\alpha+1}(-2\eta)}{\theta_1'(0)} u_\alpha \,, \qquad \bar{\nu}_\alpha = \imath^{\delta_{\alpha,2}} \frac{\theta_1(\hbar)\theta_{\alpha+1}(2\hbar)\theta_{\alpha+1}(-2\eta)}{\theta_1'(0)\theta_1(2\hbar)} \bar{u}_\alpha \,, \tag{1.36}$$

where $\bar{u}_\alpha$ are constants from (1.35), $u_\alpha$ are parameters of (1.4) and the rest of parameters entering (1.4) are

$$\begin{gathered} \sqrt{K_i} = \Delta_i = J_j + J_k - J_i \,, \quad \rho_i = \sqrt{J_j J_k \Delta_j \Delta_k} \,, \quad \text{for every cyclic triple } i,j,k \in \{1,2,3\} \,, \\ \mu_k = \frac{\delta_2(\hbar,\eta)}{J_k} - \delta_1(\hbar,\eta) \,, \qquad \lambda_k = \frac{\rho_k}{2J_1J_2J_3}\delta_2(\hbar,\eta) \,, \end{gathered} \tag{1.37}$$

where $\delta_1(\hbar, \eta)$, $\delta_2(\hbar, \eta)$ are the Casimir values of the original Sklyanin representation

$$\delta_1(\hbar, \eta) = 4\theta_1^2(2\eta + \hbar) \,, \qquad \delta_2 = 4\theta_1(2\eta)\theta_1(2\eta + 2\hbar) \,. \tag{1.38}$$

The case $\xi = \eta$ is considered separately. Then we deal with 7-constant Hamiltonian, and it is represented in the form (1.10). The statement (1.9) is reformulated as $\hat{H}^{\text{7vD}}(z|\tau; 2\hbar, \nu, \bar{\nu}, \xi - \hbar) = \hat{H}^{\text{8vD}}(z|\tau; 2\hbar, \nu, \bar{\nu}, \xi - \hbar, \xi - \hbar) = \hat{\mathbb{S}}^\xi_0(z|\bar{u})\hat{\mathbb{S}}^\xi_0(z|u)$ (see (5.8)).

Finally, using the modular triple, we find one more additional commuting quantum Hamiltonian for (1.11)-(1.12). Namely, consider

$$\hat{H}_1^{\text{8vD}} := \hat{H}^{\text{8vD}}\left(z|2\hbar;\frac{\tau}{2},\widetilde{\nu}^1,\bar{\widetilde{\nu}}^1,\eta_1,\xi_1\right), \qquad \hat{H}_2^{\text{8vD}} := \hat{H}^{\text{8vD}}\left(\frac{z}{2\hbar}|\frac{\tau}{2\hbar};\frac{1}{4\hbar},\widetilde{\nu}^2,\bar{\widetilde{\nu}}^{\,2},\frac{\eta_2}{2\hbar},\frac{\xi_2}{2\hbar}\right). \tag{1.39}$$

We show that

$$\begin{aligned}
&\left[\hat{H}^{\text{8vD}},\hat{H}_1^{\text{8vD}}\right]_- = 0\,, && \text{for } \eta_1+\xi_1 = \eta+\xi+\frac{k}{4}\,, k\in\mathbb{Z}\,,\\
&\left[\hat{H}^{\text{8vD}},\hat{H}_2^{\text{8vD}}\right]_- = 0\,, && \text{for any } \eta,\xi,\eta_2,\xi_2\,,\\
&\left[\hat{H}_1^{\text{8vD}},\hat{H}_2^{\text{8vD}}\right]_- = 0\,, && \text{for any } \eta_1,\xi_1 \text{ and } \eta_2+\xi_2 = \frac{N\hbar}{2}\,, N\in\mathbb{Z}\,.
\end{aligned} \tag{1.40}$$

This is a generalization of the previous result (1.11)-(1.12) by Rains and Ruijsenaars [32].

The degenerations of this elliptic one-body model include the important cases considered in [43, 2]. We plan to study degenerations of this elliptic construction in future work.

**Three-dimensional Sklyanin algebra** The three-dimensional Sklyanin algebra is generated by $x, y, z$ with defining relations

$$Skl_3: \qquad \begin{aligned} &ayz + bzy + cx^2 = 0\,,\\ &azx + bxz + cy^2 = 0\,,\\ &axy + byx + cz^2 = 0\,, \end{aligned} \tag{1.41}$$

where $a, b, c \in \mathbb{C}$. This algebra first appeared in [1] and was rediscovered by Feigin-Odesskii [28] during the construction of the more general elliptic algebras, which unify three- and four-dimensional Sklyanin algebra into a common family. These algebras also have a lot of important properties (see surveys [30], [40] and references therein). For generic parameters, the three-dimensional Sklyanin algebra is Koszul and has Hilbert series $(1-t)^{-3}$ [13]. From a physical point of view, finite-dimensional representations of $Skl_3$ arise in the description of vacua of marginal deformations of four-dimensional $\mathcal{N}=4$ super Yang–Mills theory [44].

This algebra is interesting to us because it has the property discussed above for $\mathcal{U}_q(sl_2)$ and four-dimensional Sklyanin algebra. Namely, in the last section, we show that in the quotient $Skl_3/G$ for the central element $G$ (6.6), the quadratic elements $L_x(u)$, $L_y(u)$, $L_z(u)$ and $R_x(u)$, $R_y(u)$, $R_z(u)$ of the form

$$\begin{aligned}
&L_x(u) = (u_x x + u_y y + u_z z)x, \quad L_y(u) = (u_x x + u_y y + u_z z)y, \quad L_z(u) = (u_x x + u_y y + u_z z)z,\\
&R_x(u) = x(u_x x + u_y y + u_z z), \quad R_y(u) = y(u_x x + u_y y + u_z z), \quad R_z(u) = z(u_x x + u_y y + u_z z)
\end{aligned} \tag{1.42}$$

satisfy the defining relations of $Skl_3$ (1.41) with parameters $\big(b(c^3-a^3), a(b^3-c^3), c(a^3-b^3)\big)$ instead of $\big(a, b, c\big)$. It would be interesting to determine which other quadratic algebras admit a similar construction.

**Organization of the paper.** In Section 2 we derive the defining relations and the two Casimir elements of $\mathcal{A}(\mathcal{A}_2)$ and identify this algebra as a specialization of the Mikhailov-Skrypnyk algebra.

Section 3 constructs its representations by elliptic difference operators, introduces the deformation parameter $\xi$, and treats the specializations to $\mathcal{A}$ and $\mathcal{B}C_1$. In Section 4 we derive the cross-relations defining the elliptic modular triple and prove that quadratic elements from distinct copies of the Sklyanin algebra commute. Section 5 applies these results to the one-body Ruijsenaars-van Diejen model and recovers the Poisson algebra of the Zhukovsky-Volterra gyrostat in the classical limit. Section 6 treats the analogous quadratic construction for the three-dimensional Sklyanin algebra. A list of notations is given after Section 6. Appendix A collects the elliptic identities used in the proofs. Appendix B recalls how the quantum exchange relation and the reflection equation give the defining relations of $\mathcal{A}$ and $\mathcal{B}C_1$, respectively.

# 2 Commutation relations in $\mathcal{A}_2$

## 2.1 Casimir elements

**Sklyanin algebra $\mathcal{A}$.** The four-dimensional Sklyanin algebra (over $\mathbb{C}$) is generated by $S_0, S_1, S_2, S_3$ with defining relations

$$[S_i, S_j]_- = \imath[S_0, S_k]_+ \,, \tag{2.1}$$

$$[S_0, S_k]_- = \imath\frac{J_j - J_i}{J_k}\,[S_i, S_j]_+ \,, \tag{2.2}$$

for every cyclic permutation $(i, j, k)$ of $(1, 2, 3)$. The two central quadratic elements are

$$\begin{aligned} \Omega_1 &= S_0^2 + S_1^2 + S_2^2 + S_3^2 \,, \\ \Omega_2 &= J_1 S_1^2 + J_2 S_2^2 + J_3 S_3^2 \,. \end{aligned} \tag{2.3}$$

That is, $[\Omega_1, S_\alpha]_- = [\Omega_2, S_\alpha]_- = 0$ for any $S_\alpha$. The centers of the three- and four-dimensional Sklyanin algebras are studied in [6], [41]. The commutation relations (2.1)-(2.2) can also be written as

$$[S_i, S_j]_- = \imath\,[S_0, S_k]_+ \,, \tag{2.4}$$

$$[S_0, S_k]_- = -\imath\,\frac{\mathcal{K}_j - \mathcal{K}_i}{\mathcal{K}_k}\,[S_i, S_j]_+ \,. \tag{2.5}$$

The parameters $\mathcal{K}_i$ describe the same algebra (2.1)-(2.2) with the opposite sign convention in the second relation. They are related to $J_i$ by

$$\frac{\mathcal{K}_j - \mathcal{K}_i}{\mathcal{K}_k} = -\frac{J_j - J_i}{J_k} \,. \tag{2.6}$$

Since (2.6) is invariant under a common rescaling of either triple, the normalization can be chosen so that

$$J_i = \mathcal{K}_i(\mathcal{K}_i - \mathcal{K}_j - \mathcal{K}_k) \,. \tag{2.7}$$

To state the inverse transformation, define

$$\Delta_i = J_j + J_k - J_i \neq 0 \,. \tag{2.8}$$

Then

$$\mathcal{K}_i = C J_i(J_i - J_j - J_k) \,, \qquad C^2 = -\frac{1}{\Delta_1\Delta_2\Delta_3} \,. \tag{2.9}$$

The elliptic parameterization is specified below in (3.9)–(3.10). Under this change, the Casimir element $\Omega_2$ in (2.3) becomes

$$\Omega_2 = \mathcal{K}_1(\mathcal{K}_1 - \mathcal{K}_2 - \mathcal{K}_3)S_1^2 + \mathcal{K}_2(\mathcal{K}_2 - \mathcal{K}_3 - \mathcal{K}_1)S_2^2 + \mathcal{K}_3(\mathcal{K}_3 - \mathcal{K}_1 - \mathcal{K}_2)S_3^2 \,. \tag{2.10}$$

$\mathcal{B}C_1$ **Sklyanin algebra.** In [19], a three-parameter deformation of the four-dimensional Sklyanin algebra was introduced. The generators $\mathbb{S}_\alpha$ of this algebra satisfy the following relations

$$[\widetilde{\nu}_i, \widetilde{\nu}_j]_- = 0, \quad [\widetilde{\nu}_i, \mathbb{S}_\alpha]_- = 0\,, \tag{2.11}$$

$$\imath[\mathbb{S}_0, \mathbb{S}_k]_+ = [\mathbb{S}_i, \mathbb{S}_j]_-\,, \tag{2.12}$$

$$[\mathbb{S}_0, \mathbb{S}_k]_- = -\imath\frac{K_j - K_i}{K_k}\,[\mathbb{S}_i, \mathbb{S}_j]_+ + \frac{2\imath}{K_k}\,(\widetilde{\nu}_i\rho_i\mathbb{S}_j - \widetilde{\nu}_j\rho_j\mathbb{S}_i)\,. \tag{2.13}$$

As in (1.3), the three additional constants are $\widetilde{\nu}_i$. If $\widetilde{\nu}_i = 0$ for all $i$, the algebra (2.11)–(2.13) coincides with the standard Sklyanin algebra (2.4)–(2.5) upon replacing $K_i$ by $\mathcal{K}_i$.[3] The parameters $\rho_i$ are normalization factors. They can be absorbed into the definition of $\widetilde{\nu}_i$, but we keep them explicit for convenience.

The quadratic Casimir elements were obtained in [19] by calculating the quantum determinant, as for the original Sklyanin algebra $\mathcal{A}$:

$$\begin{aligned}\Omega_1^{\mathcal{B}C_1} &= \mathbb{S}_0^2 + \sum_{k=1}^{3}\mathbb{S}_k^2\,,\\ \Omega_2^{\mathcal{B}C_1} &= \sum_{k=1}^{3}\Big(K_k\big(K_k - K_i - K_j\big)\mathbb{S}_k^2 + 2\widetilde{\nu}_k\rho_k K_k\mathbb{S}_k\Big)\,.\end{aligned} \tag{2.14}$$

## 2.2 Quadratic relations for elements in $\mathcal{A}_2$

The component $\mathcal{A}_2$ in the decomposition (1.2) is spanned by the quadratic monomials $S_\alpha S_\beta$ for $\alpha, \beta = 0, ..., 3$. Consider the full $4\times 4$ matrix of quadratic monomials

$$Q_{\alpha\beta} = S_\alpha S_\beta\,, \quad \alpha, \beta = 0, ..., 3\,, \tag{2.15}$$

and let $U$ be an arbitrary element of $\mathcal{A}_1$:

$$U = \sum_{\beta=0}^{3} u_\beta S_\beta\,, \qquad u_\beta \in \mathbb{C}\,. \tag{2.16}$$

For this $U$, define

$$X_\alpha(U) := US_\alpha \in \mathcal{A}_2, \quad \alpha = 0, 1, 2, 3\,. \tag{2.17}$$

Then

$$(X_0(U), X_1(U), X_2(U), X_3(U)) = (u_0, u_1, u_2, u_3)Q\,. \tag{2.18}$$

In the same way, we define

$$(Y_0(U), Y_1(U), Y_2(U), Y_3(U))^{\mathbf{T}} = Q(u_0, u_1, u_2, u_3)^{\mathbf{T}}\,, \tag{2.19}$$

for

$$Y_\alpha(U) := S_\alpha U, \quad \alpha = 0, 1, 2, 3\,. \tag{2.20}$$

The next step is to describe the commutation relations between the elements in (2.17) and (2.20). Fix a cyclic permutation $(i, j, k)$ of $(1, 2, 3)$ and use the notation

$$J_{ij} := \frac{J_j - J_i}{J_k}, \quad \Delta_k := J_i + J_j - J_k \neq 0\,. \tag{2.21}$$

[3] In the elliptic parameterization, (3.4)–(3.11) show that $\mathcal{K}_i \propto \mathbf{K}_i(\hbar)$ while $K_i \propto \mathbf{K}_i(2\hbar)$.

In the free algebra $\mathbb{C}\langle S_0, S_1, S_2, S_3\rangle$, introduce the six polynomials defining the relations:

$$\begin{aligned}
\mathcal{F}_{ij} &:= S_iS_j - S_jS_i - \imath(S_0S_k + S_kS_0)\,,\\
\mathcal{F}_{0k} &:= S_0S_k - S_kS_0 - \imath J_{ij}(S_iS_j + S_jS_i)\,,\\
\mathcal{F}_{jk} &:= S_jS_k - S_kS_j - \imath(S_0S_i + S_iS_0)\,,\\
\mathcal{F}_{0i} &:= S_0S_i - S_iS_0 - \imath J_{jk}(S_jS_k + S_kS_j)\,,\\
\mathcal{F}_{ki} &:= S_kS_i - S_iS_k - \imath(S_0S_j + S_jS_0)\,,\\
\mathcal{F}_{0j} &:= S_0S_j - S_jS_0 - \imath J_{ki}(S_kS_i + S_iS_k)\,.
\end{aligned} \tag{2.22}$$

Each $\mathcal{F}$ vanishes in $\mathcal{A}$ by (2.1) and (2.2). For any $U \in \mathcal{A}_1$, define two cubic expressions

$$G_1(U) := \Delta_k(S_iUS_j - S_jUS_i) - 2\imath J_k(S_0US_k + S_kUS_0) \in \mathcal{A}_3\,, \tag{2.23}$$

$$G_2(U) := \Delta_k(S_0US_k - S_kUS_0) - 2\imath(J_j - J_i)(S_iUS_j + S_jUS_i) \in \mathcal{A}_3\,. \tag{2.24}$$

To compute the maps $G_1$ and $G_2$, it suffices to evaluate them on the generators $S_\alpha$ for $\alpha = 0, ..., 3$ because (2.23) and (2.24) are linear in $U$.

**Lemma 2.1** *The values of (2.23) and (2.24) on the generators are given in the following table:*

| $U$ | $G_1(U)$ | $G_2(U)$ |
|---|---|---|
| $S_0$ | $2\imath(\Omega_2 - J_k\Omega_1)S_k$ | $0$ |
| $S_i$ | $0$ | $2\imath(\Omega_2 - J_j\Omega_1)S_j$ |
| $S_j$ | $0$ | $-2\imath(\Omega_2 - J_i\Omega_1)S_i$ |
| $S_k$ | $-2\imath\Omega_2S_0$ | $0$ |

**Proof.** To prove the statement, it suffices to verify the eight relations below.

For $U = S_0$:

$$\begin{aligned}
G_1(S_0) - 2\imath(\Omega_2 - J_k\Omega_1)S_k = {}& -J_iS_j\mathcal{F}_{0i} + J_jS_i\mathcal{F}_{0j} + J_k\mathcal{F}_{ij}S_0 + \imath J_k\mathcal{F}_{0k}S_0 - \\
& -\imath(J_j - J_k)S_j\mathcal{F}_{jk} + \imath(J_i - J_k)S_i\mathcal{F}_{ki}\,,
\end{aligned} \tag{2.25}$$

$$\begin{aligned}
G_2(S_0) = {}& \imath\Delta_k(S_0\mathcal{F}_{ij} - \mathcal{F}_{ij}S_0) + (J_j - J_k)(S_j\mathcal{F}_{jk} + \mathcal{F}_{jk}S_j) + \\
& +\imath J_i(S_j\mathcal{F}_{0i} - \mathcal{F}_{0i}S_j) + (J_k - J_i)(S_i\mathcal{F}_{ki} + \mathcal{F}_{ki}S_i) - \\
& -\imath J_j(S_i\mathcal{F}_{0j} - \mathcal{F}_{0j}S_i)\,.
\end{aligned} \tag{2.26}$$

For $U = S_i$:

$$\begin{aligned}
G_1(S_i) = {}& \Delta_k(S_i\mathcal{F}_{ij} + \mathcal{F}_{ij}S_i) + (J_k - J_j)(S_k\mathcal{F}_{jk} + \mathcal{F}_{jk}S_k) + \\
& +\imath J_i(S_k\mathcal{F}_{0i} - \mathcal{F}_{0i}S_k) + \imath J_j(S_0\mathcal{F}_{ki} - \mathcal{F}_{ki}S_0) - \\
& -J_j(S_0\mathcal{F}_{0j} + \mathcal{F}_{0j}S_0)\,,
\end{aligned} \tag{2.27}$$

$$\begin{aligned}
G_2(S_i) - 2\imath(\Omega_2 - J_j\Omega_1)S_j = {}& \imath(J_j - J_i)\mathcal{F}_{ij}S_i + J_k\mathcal{F}_{0k}S_i + \imath(J_k - J_j)S_k\mathcal{F}_{jk} + \\
& +J_iS_k\mathcal{F}_{0i} - J_jS_0\mathcal{F}_{ki} + \imath J_jS_0\mathcal{F}_{0j}\,.
\end{aligned} \tag{2.28}$$

For $U = S_j$:

$$G_1(S_j) = \Delta_k(S_j\mathcal{F}_{ij} + \mathcal{F}_{ij}S_j) - \imath J_i(S_0\mathcal{F}_{jk} - \mathcal{F}_{jk}S_0) + \\ + J_i(S_0\mathcal{F}_{0i} + \mathcal{F}_{0i}S_0) + (J_k - J_i)(S_k\mathcal{F}_{ki} + \mathcal{F}_{ki}S_k) + \\ + \imath J_j(S_k\mathcal{F}_{0j} - \mathcal{F}_{0j}S_k)\,, \tag{2.29}$$

$$G_2(S_j) + 2\imath(\Omega_2 - J_i\Omega_1)S_i = -\imath(J_j - J_i)\mathcal{F}_{ij}S_j + J_k\mathcal{F}_{0k}S_j + J_iS_0\mathcal{F}_{jk} - \imath J_iS_0\mathcal{F}_{0i} + \\ + \imath(J_k - J_i)S_k\mathcal{F}_{ki} + J_jS_k\mathcal{F}_{0j}\,. \tag{2.30}$$

For $U = S_k$:

$$G_1(S_k) + 2\imath\Omega_2S_0 = J_k\mathcal{F}_{ij}S_k - \imath J_k\mathcal{F}_{0k}S_k - J_iS_i\mathcal{F}_{jk} - \imath J_iS_i\mathcal{F}_{0i} - \\ - J_jS_j\mathcal{F}_{ki} - \imath J_jS_j\mathcal{F}_{0j}\,, \tag{2.31}$$

$$G_2(S_k) = -\imath\Delta_k(S_k\mathcal{F}_{ij} - \mathcal{F}_{ij}S_k) - \imath J_i(S_i\mathcal{F}_{jk} - \mathcal{F}_{jk}S_i) - \\ - J_i(S_i\mathcal{F}_{0i} + \mathcal{F}_{0i}S_i) - \imath J_j(S_j\mathcal{F}_{ki} - \mathcal{F}_{ki}S_j) - \\ - J_j(S_j\mathcal{F}_{0j} + \mathcal{F}_{0j}S_j)\,, \tag{2.32}$$

where the polynomials $\mathcal{F}$ are defined in (2.22) and vanish in $\mathcal{A}$. Therefore, the right-hand sides of the above identities vanish in $\mathcal{A}$. This completes the proof. ∎

Next, we compute the commutation relations.

**Theorem 2.1** *Let $U = \sum_{\alpha=0}^{3} u_\alpha S_\alpha \in \mathcal{A}_1$, $V = \sum_{\alpha=0}^{3} v_\alpha S_\alpha \in \mathcal{A}_1$ and*

$$X_\alpha(U) = US_\alpha, \quad X_\alpha(V) = VS_\alpha, \qquad \alpha = 0,1,2,3\,, \\ Y_\alpha(V) = S_\alpha V, \quad Y_\alpha(U) = S_\alpha U, \qquad \alpha = 0,1,2,3\,. \tag{2.33}$$

*Then, for every cyclic permutation $(i,j,k)$ of $(1,2,3)$ with $\Delta_k := J_i + J_j - J_k \neq 0$, the following relations hold:*

$$X_i(V)X_j(U) - X_j(V)X_i(U) = 2\imath\,\frac{J_k}{\Delta_k}\big(X_0(V)X_k(U) + X_k(V)X_0(U)\big) + \\ + 2\imath\frac{u_0}{\Delta_k}(\Omega_2 - J_k\Omega_1)X_k(V) - 2\imath\frac{u_k}{\Delta_k}\,\Omega_2X_0(V)\,, \tag{2.34}$$

$$X_0(V)X_k(U) - X_k(V)X_0(U) = 2\imath\,\frac{J_j - J_i}{\Delta_k}\big(X_i(V)X_j(U) + X_j(V)X_i(U)\big) + \\ + 2\imath\frac{u_i}{\Delta_k}(\Omega_2 - J_j\Omega_1)X_j(V) - 2\imath\frac{u_j}{\Delta_k}(\Omega_2 - J_i\Omega_1)X_i(V)\,, \tag{2.35}$$

$$Y_i(V)Y_j(U) - Y_j(V)Y_i(U) = 2\imath\,\frac{J_k}{\Delta_k}\big(Y_0(V)Y_k(U) + Y_k(V)Y_0(U)\big) + \\ + 2\imath\frac{v_0}{\Delta_k}(\Omega_2 - J_k\Omega_1)Y_k(U) - 2\imath\frac{v_k}{\Delta_k}\,\Omega_2Y_0(U)\,, \tag{2.36}$$

$$Y_0(V)Y_k(U) - Y_k(V)Y_0(U) = 2\imath\,\frac{J_j - J_i}{\Delta_k}\big(Y_i(V)Y_j(U) + Y_j(V)Y_i(U)\big) + \\ + 2\imath\frac{v_i}{\Delta_k}(\Omega_2 - J_j\Omega_1)Y_j(U) - 2\imath\frac{v_j}{\Delta_k}(\Omega_2 - J_i\Omega_1)Y_i(U)\,, \tag{2.37}$$

$$X_i(U)Y_j(V) - X_j(U)Y_i(V) = \imath\big(X_0(U)Y_k(V) + X_k(U)Y_0(V)\big)\,, \tag{2.38}$$

$$X_0(U)Y_k(V) - X_k(U)Y_0(V) = \imath\frac{J_j - J_i}{J_k}\big(X_i(U)Y_j(V) + X_j(U)Y_i(V)\big)\,. \tag{2.39}$$

*In particular, when $U = V$, we have*

$$[X_i(U), X_j(U)]_- = 2\imath\,\frac{J_k}{\Delta_k}\,[X_0(U), X_k(U)]_+ + \tag{2.40}$$
$$+2\imath\frac{u_0}{\Delta_k}\big(\Omega_2 - J_k\Omega_1\big)X_k(U) - 2\imath\frac{u_k}{\Delta_k}\,\Omega_2 X_0(U)\,,$$

$$[X_0(U), X_k(U)]_- = 2\imath\,\frac{J_j - J_i}{\Delta_k}\,[X_i(U), X_j(U)]_+ + \tag{2.41}$$
$$+2\imath\frac{u_i}{\Delta_k}\big(\Omega_2 - J_j\Omega_1\big)X_j(U) - 2\imath\frac{u_j}{\Delta_k}\big(\Omega_2 - J_i\Omega_1\big)X_i(U)\,,$$

$$[Y_i(U), Y_j(U)]_- = 2\imath\,\frac{J_k}{\Delta_k}\,[Y_0(U), Y_k(U)]_+ + \tag{2.42}$$
$$+2\imath\frac{u_0}{\Delta_k}\big(\Omega_2 - J_k\Omega_1\big)Y_k(U) - 2\imath\frac{u_k}{\Delta_k}\,\Omega_2 Y_0(U)\,,$$

$$[Y_0(U), Y_k(U)]_- = 2\imath\,\frac{J_j - J_i}{\Delta_k}\,[Y_i(U), Y_j(U)]_+ + \tag{2.43}$$
$$+2\imath\frac{u_i}{\Delta_k}\big(\Omega_2 - J_j\Omega_1\big)Y_j(U) - 2\imath\frac{u_j}{\Delta_k}\big(\Omega_2 - J_i\Omega_1\big)Y_i(U)\,,$$

$$X_i(U)Y_j(U) - X_j(U)Y_i(U) = \imath\big(X_0(U)Y_k(U) + X_k(U)Y_0(U)\big)\,, \tag{2.44}$$

$$X_0(U)Y_k(U) - X_k(U)Y_0(U) = \imath\frac{J_j - J_i}{J_k}\big(X_i(U)Y_j(U) + X_j(U)Y_i(U)\big)\,. \tag{2.45}$$

**Proof.** Using the results of Lemma 2.1, we multiply (2.23) on the left by $V$:

$$\begin{aligned} VG_1(U) &= \Delta_k\big(VS_iUS_j - VS_jUS_i\big) - 2\imath J_k\big(VS_0US_k + VS_kUS_0\big) = \\ &= 2\imath u_0 V\big(\Omega_2 - J_k\Omega_1\big)S_k - 2\imath u_k V\Omega_2 S_0 = 2\imath u_0\big(\Omega_2 - J_k\Omega_1\big)VS_k - 2\imath u_k\Omega_2 VS_0\,. \end{aligned} \tag{2.46}$$

In the last equality, we use the fact that $\Omega_1, \Omega_2$ are central: they commute with every element of $\mathcal{A}$. In terms of $X_\alpha(U)$, we obtain the relation

$$\begin{aligned} \Delta_k\big(X_i(V)X_j(U) - X_j(V)X_i(U)\big) - 2\imath J_k\big(X_0(V)X_k(U) + X_k(V)X_0(U)\big) = \\ = 2\imath u_0\big(\Omega_2 - J_k\Omega_1\big)X_k(V) - 2\imath u_k\Omega_2 X_0(V)\,, \end{aligned} \tag{2.47}$$

which coincides with (2.34). Similarly, multiplying (2.24) on the left by $V$ gives (2.35). To obtain (2.36) and (2.37), we multiply $G_1(V)$ and $G_2(V)$ (computed from (2.23) and (2.24)) on the right by $U$. Multiplying (2.1) and (2.2) on the left by $U$ and on the right by $V$ yields

$$\begin{aligned} US_iS_jV - US_jS_iV &= \imath\big(US_0S_kV + US_kS_0V\big)\,, \\ US_0S_kV - US_kS_0V &= \imath\frac{J_j - J_i}{J_k}\big(US_iS_jV + US_jS_iV\big)\,. \end{aligned} \tag{2.48}$$

In terms of $X_\alpha(U), Y_\beta(V)$, the two relations in (2.48) coincide with (2.38) and (2.39). This proves (2.34)–(2.39). ■

**Algebra of the normalized elements $X_\alpha(U)$.** Redefine the generators as follows:

$$Z_0 = X_0(U), \qquad Z_i = \frac{1}{2}\sqrt{\frac{\Delta_j\Delta_k}{J_jJ_k}}X_i(U)\,. \tag{2.49}$$

Then the commutation relations take the form:

$$\begin{gathered}
[Z_i, Z_j]_- = \imath\,[Z_0, Z_k]_+ + \imath u_0\mu_k Z_k - \imath u_k\lambda_k Z_0\,,\\
[Z_0, Z_k]_- = -\imath\,\frac{K_j - K_i}{K_k}\,[Z_i, Z_j]_+ + 2\imath\frac{u_i}{K_k}\,\rho_i\mu_j Z_j - 2\imath\frac{u_j}{K_k}\,\rho_j\mu_i Z_i\,,\\
\sqrt{K_i} = \Delta_i = J_j + J_k - J_i\,, \quad \rho_i = \sqrt{J_jJ_k\Delta_j\Delta_k}\,,\\
\mu_k = \frac{\rho_k\sqrt{K_k}}{\rho_i\rho_j}\,\Omega_2 - \Omega_1, \quad \lambda_k = \frac{\sqrt{K_iK_jK_k}}{2\rho_i\rho_j}\,\Omega_2
\end{gathered} \tag{2.50}$$

for every cyclic triple $i, j, k \in \{1, 2, 3\}$. The square roots in (2.49) and below are chosen consistently: $\sqrt{\Delta_j\Delta_k/(J_jJ_k)} = \rho_i/(J_jJ_k)$ and $\sqrt{K_iK_jK_k} = \Delta_i\Delta_j\Delta_k$, with $\rho_1\rho_2\rho_3 = J_1J_2J_3\Delta_1\Delta_2\Delta_3$. The above statement and all other results for the algebra (2.50) also hold if $X_\alpha(U)$ is replaced by $Y_\alpha(U)$ in (2.49).

Let us consider the algebra defined by the relations (2.50) separately. Define an algebra with four generators $\mathcal{Z}_0, \mathcal{Z}_1, \mathcal{Z}_2, \mathcal{Z}_3$ and relations

$$\begin{gathered}
[\mathcal{Z}_i, \mathcal{Z}_j]_- = \imath\,[\mathcal{Z}_0, \mathcal{Z}_k]_+ + \imath u_0\mu_k\mathcal{Z}_k - \imath u_k\lambda_k\mathcal{Z}_0\,,\\
[\mathcal{Z}_0, \mathcal{Z}_k]_- = -\imath\,\frac{K_j - K_i}{K_k}\,[\mathcal{Z}_i, \mathcal{Z}_j]_+ + 2\imath\frac{u_i}{K_k}\,\rho_i\mu_j\mathcal{Z}_j - 2\imath\frac{u_j}{K_k}\,\rho_j\mu_i\mathcal{Z}_i\,,\\
\sqrt{K_i} = \Delta_i = J_j + J_k - J_i\,, \quad \rho_i = \sqrt{J_jJ_k\Delta_j\Delta_k}\,,\\
\mu_k = \frac{\rho_k\sqrt{K_k}}{\rho_i\rho_j}\,\delta_2 - \delta_1, \quad \lambda_k = \frac{\sqrt{K_iK_jK_k}}{2\rho_i\rho_j}\,\delta_2\,, \text{ for every cyclic triple } i, j, k\, \in \{1, 2, 3\}\,.
\end{gathered} \tag{2.51}$$

These are the relations of $\mathcal{A}(\mathcal{A}_2)$ in (1.4). Here $\delta_1, \delta_2 \in \mathbb{C}$ are parameters corresponding to the values of the Casimir elements $\Omega_1$ and $\Omega_2$, respectively. The parameters $K_i = \Delta_i^2$ are defined in (2.51) and are distinct from the auxiliary parameters $\mathcal{K}_i$ in (2.7)–(2.9). Their elliptic parametrizations are related by doubling the argument, as shown in (3.10) and (3.27). After imposing $\Omega_1 = \delta_1$ and $\Omega_2 = \delta_2$, the substitution $\mathcal{Z}_\alpha = Z_\alpha$ turns (2.51) into an identity. Thus, Theorem 2.1 gives a realization of the algebra (2.51) by the classes of the quadratic elements $Z_\alpha$ in the quotient $\mathcal{A}/(\Omega_1 - \delta_1, \Omega_2 - \delta_2)$. The image lies in the even subalgebra of this quotient.

The algebra (2.51) is related by an affine change of generators to a specialization of the family introduced in [23] in the context of the quantization of the Zhukovsky-Volterra gyrostat. The explicit correspondence is given at the end of this section.

**Proposition 2.1** *The two elements*

$$\Omega_1^{\mathcal{Z}} = \mathcal{Z}_0^2 + a_0\mathcal{Z}_0 + \sum_{k=1}^{3}\left(\mathcal{Z}_k^2 - u_k\lambda_k\mathcal{Z}_k\right) \tag{2.52}$$

*with*

$$a_0 = \frac{u_0}{3}\left(\mu_i + \mu_j + \mu_k - \frac{K_i\lambda_i}{2\rho_i} - \frac{K_j\lambda_j}{2\rho_j} - \frac{K_k\lambda_k}{2\rho_k}\right) \tag{2.53}$$

*and*[4]

$$\Omega_2^{\mathcal{Z}} = \sum_{k=1}^{3} K_k\big(K_k - K_i - K_j\big)\mathcal{Z}_k^2 + \sum_{k=0}^{3} b_k \mathcal{Z}_k \tag{2.54}$$

*with*

$$\begin{aligned} b_0 &= u_0 \frac{K_i K_j K_k}{4 J_i J_j J_k}\,\delta_2\,, \\ b_k &= 2u_k\rho_k K_k\Big( -\delta_1 + \frac{1}{2}\Big(\frac{1}{J_i} + \frac{1}{J_j} - \frac{1}{J_k} + \frac{J_i}{2J_jJ_k} + \frac{J_j}{2J_iJ_k} + \frac{J_k}{2J_iJ_j}\Big)\delta_2\Big) \end{aligned} \tag{2.55}$$

*belong to the center of the algebra (2.51). Here $a_0, b_0$ are independent of the cyclic choice of $(i, j, k)$.*

In fact, since the algebra (2.51) is a particular case of the Mikhailov-Skrypnyk algebra, the Casimir elements (2.52), (2.54) for the algebra (2.51) can be deduced from the Casimir elements of the Mikhailov-Skrypnyk algebra given in [23].

As we have seen above, (2.17) and (2.49) define a map from the algebra $\mathcal{A}(\mathcal{A}_2)$ of (2.51) to the central quotient $\mathcal{A}/(\Omega_1 - \delta_1, \Omega_2 - \delta_2)$ of the Sklyanin algebra (2.1)–(2.2). This map determines the values of the central elements (2.52) and (2.54), as follows.

**Proposition 2.2** *Under the substitution $\mathcal{Z}_\alpha = Z_\alpha$, with $Z_\alpha$ defined by (2.49) in the central quotient $\mathcal{A}/(\Omega_1 - \delta_1, \Omega_2 - \delta_2)$, the Casimir elements $\Omega_1^{\mathcal{Z}}$ in (2.52) and $\Omega_2^{\mathcal{Z}}$ in (2.54) both vanish:*

$$\Omega_1^{\mathcal{Z}}\Big|_{\mathcal{Z}_\alpha = Z_\alpha} = \Omega_2^{\mathcal{Z}}\Big|_{\mathcal{Z}_\alpha = Z_\alpha} = 0\,. \tag{2.56}$$

**Proof.** The verification is straightforward. Substituting (2.49) into (2.52) and (2.54) and calculating in the free algebra gives a sum of terms proportional to the polynomials in (2.22), which vanish for the generators of the Sklyanin algebra. ■

## 2.3 Realizations of the quadratic algebras $\mathcal{A}$ and $\mathcal{B}C_1$ through $\mathcal{A}(\mathcal{A}_2)$

The algebra (2.51) has several important degenerations, which we will list below.

**Sklyanin algebra $\mathcal{A}$.** When $\delta_1 = \delta_2 = 0$, the algebra (2.51) reduces to the ordinary four-dimensional Sklyanin algebra

$$[Z_i, Z_j]_- = \imath\,[Z_0, Z_k]_+\,, \tag{2.57}$$

$$[Z_0, Z_k]_- = -\imath\,\frac{K_j - K_i}{K_k}\,[Z_i, Z_j]_+\,. \tag{2.58}$$

There is another way to obtain the Sklyanin algebra from the relations (2.51). One can shift the generators of the algebra by central elements $C_\alpha$:

$$Z_\alpha = Z'_\alpha + C_\alpha, \quad \alpha = 0, 1, 2, 3\,. \tag{2.59}$$

Then (2.50) becomes

$$\begin{aligned} [Z'_i, Z'_j]_- &= \imath\,[Z'_0, Z'_k]_+ + \imath u_0\mu_k Z'_k - \imath u_k\lambda_k Z'_0 + \\ &+ 2\imath\big(C_0 Z'_k + C_k Z'_0 + C_0 C_k\big) + \imath u_0\mu_k C_k - \imath u_k\lambda_k C_0\,, \end{aligned} \tag{2.60}$$

[4] In the sums over $k$ in (2.54), the indices $i, j$ are the two elements of $\{1, 2, 3\} \setminus \{k\}$.

$$
\begin{aligned}
[Z_0', Z_k']_- = -\imath\, \frac{K_j - K_i}{K_k}\, [Z_i', Z_j']_+ + 2\imath \frac{u_i}{K_k}\, \rho_i \mu_j Z_j' - 2\imath \frac{u_j}{K_k}\, \rho_j \mu_i Z_i' - \\
-2\imath\, \frac{K_j - K_i}{K_k}\, \big(C_i Z_j' + C_j Z_i' + C_i C_j\big) + 2\imath \frac{u_i}{K_k}\, \rho_i \mu_j C_j - 2\imath \frac{u_j}{K_k}\, \rho_j \mu_i C_i\,.
\end{aligned}
\tag{2.61}
$$

If

$$
C_0 = -\frac{1}{2}\, u_0 \mu_k, \quad C_k = \frac{1}{2}\, u_k \lambda_k, \quad C_0 C_k = 0\,. \tag{2.62}
$$

(2.60) turns into

$$
[Z_i', Z_j']_- = \imath\, [Z_0', Z_k']_+\,. \tag{2.63}
$$

In order to obtain the second set of relations that define the Sklyanin algebra, we put $u_k = 0$ for every $k \in \{1, 2, 3\}$. Then (2.61) takes the desired form

$$
[Z_0', Z_k']_- = -\imath\, \frac{K_j - K_i}{K_k}\, [Z_i', Z_j']_+\,. \tag{2.64}
$$

Thus, to obtain the Sklyanin algebra from (2.51), we only need to shift the generator $Z_0$. This shift is well defined if $\mu_1 = \mu_2 = \mu_3$. For example, this equality holds when $\delta_2 = 0$.

$BC_1$ **Sklyanin algebra.** This algebra can also be obtained from (2.51). In both reductions below, the parameters $K_i$ are kept fixed. The first reduction uses the shift (2.59)–(2.62) with $\delta_2 = 0$. Then $\lambda_k = 0$ and $\mu_k = -\delta_1$ for all $k \in \{1, 2, 3\}$. Therefore, the algebra homomorphism is given by the identification

$$
\mathbb{S}_0 = Z_0 - \frac{1}{2}\, u_0 \delta_1, \quad \mathbb{S}_k = Z_k, \quad \widetilde{\nu}_k = -u_k \delta_1, \quad k \in \{1, 2, 3\}\,. \tag{2.65}
$$

In the second case, we also apply the shift (2.59)–(2.62), but now we set $u_0 = 0$. The generators and structure constants of the algebra (2.11)–(2.13) can be represented as follows:

$$
\mathbb{S}_0 = Z_0, \quad \mathbb{S}_k = Z_k - \frac{1}{2}\, u_k \lambda_k, \quad \widetilde{\nu}_k = u_k \mu_k, \quad k \in \{1, 2, 3\}\,. \tag{2.66}
$$

In this case, the parameterization

$$
\begin{aligned}
\sqrt{K_k} = \Delta_k = J_j + J_i - J_k\,, \quad \rho_k = \sqrt{J_j J_i \Delta_j \Delta_i}\,, \\
\mu_k = \frac{\rho_k \sqrt{K_k}}{\rho_i \rho_j}\, \delta_2 - \delta_1\,, \quad \lambda_k = \frac{\sqrt{K_i K_j K_k}}{2 \rho_i \rho_j}\, \delta_2
\end{aligned}
\tag{2.67}
$$

is assumed for cyclic triples $(i, j, k)$ of $(1, 2, 3)$.

## 2.4 Exact relation of $\mathcal{A}(\mathcal{A}_2)$ to the Mikhailov-Skrypnyk algebra

Here we explain how the algebra $\mathcal{A}(\mathcal{A}_2)$ (2.51) is related to a more general case studied in [23].

**The general family.** Following [23], for every cyclic permutation $(i, j, k)$ of $(1, 2, 3)$, set

$$
\mathcal{D}_k = 1 + h^2 \mathbf{j}_k \big(\mathbf{c} + (\mathbf{j}_i + \mathbf{j}_j)\mathbf{d}\big) \neq 0\,. \tag{2.68}
$$

Then the Mikhailov-Skrypnyk algebra $\mathcal{M}S$ is defined by four generators $\mathcal{S}_0, \mathcal{S}_1, \mathcal{S}_2, \mathcal{S}_3$ satisfying the relations

$$\mathcal{M}S: \quad \begin{aligned} [\mathcal{S}_i, \mathcal{S}_j]_- &= -\frac{h}{\mathcal{D}_k}\Big((\mathbf{c}+\mathbf{d}\mathbf{j}_k)[\mathcal{S}_0, \mathcal{S}_k]_+ + 2(\mathbf{a}+\mathbf{b}\mathbf{j}_k)\mathcal{S}_k + 2\mathbf{d}\mathbf{k}_k\mathcal{S}_0 + 2\mathbf{b}\mathbf{k}_k\Big), \\ [\mathcal{S}_k, \mathcal{S}_0]_- &= -h\Big((\mathbf{j}_i - \mathbf{j}_j)[\mathcal{S}_i, \mathcal{S}_j]_+ + 2(\mathbf{k}_i\mathcal{S}_j - \mathbf{k}_j\mathcal{S}_i)\Big). \end{aligned} \tag{2.69}$$

Here $\mathbf{a}, \mathbf{b}, \mathbf{c}, \mathbf{d}, h, \mathbf{j}_i, \mathbf{k}_i$ are the parameters of [23]. In particular, $\mathbf{j}_i, \mathbf{k}_i$ should be distinguished from $J_i, K_i$ in (2.50). We do not know of an $R$-matrix quantization that gives this algebra in the general case.

**Identification with (2.51).** Fix the values $\Omega_1 = \delta_1$, $\Omega_2 = \delta_2$ of the original central elements and assume $J_i\Delta_i \neq 0$. Recall the definitions in (2.51) and choose the square roots consistently so that

$$\begin{aligned} \Delta_i = J_j + J_k - J_i, \qquad K_i = \Delta_i^2, \\ \rho_1\rho_2\rho_3 = J_1J_2J_3\Delta_1\Delta_2\Delta_3. \end{aligned} \tag{2.70}$$

The coefficients in (2.51) then obey

$$\mu_i = \frac{\delta_2}{J_i} - \delta_1, \qquad \lambda_i = \frac{\delta_2\rho_i}{2J_1J_2J_3}, \qquad \mu_i - \mu_j = \frac{\delta_2(K_i - K_j)}{4J_1J_2J_3}. \tag{2.71}$$

Define an invertible affine change of generators by

$$\mathcal{S}_0 = -\Delta_1\Delta_2\Delta_3\mathcal{Z}_0, \qquad \mathcal{S}_i = \Delta_i\Big(\mathcal{Z}_i - \frac{u_i\lambda_i}{2}\Big). \tag{2.72}$$

The corresponding parameters in (2.69) are

$$\begin{aligned} h = -\imath, \qquad \mathbf{c} = 1, \qquad \mathbf{d} = 0, \qquad \mathbf{j}_i = K_i + 1, \qquad \mathbf{k}_i = u_i\rho_i\Delta_i\mu_i, \\ \mathbf{a} = \frac{u_0\Delta_1\Delta_2\Delta_3}{2}\Big(\frac{\delta_2(K_i+1)}{4J_1J_2J_3} - \mu_i\Big), \qquad \mathbf{b} = -\frac{u_0\Delta_1\Delta_2\Delta_3\delta_2}{8J_1J_2J_3}. \end{aligned} \tag{2.73}$$

The expression for $\mathbf{a}$ is independent of $i$ by (2.71). Indeed, substitution of (2.72) into (2.50), using (2.71), gives

$$\begin{aligned} [\mathcal{S}_i, \mathcal{S}_j]_- &= -\frac{\imath}{K_k}[\mathcal{S}_0, \mathcal{S}_k]_+ + \frac{\imath u_0\Delta_1\Delta_2\Delta_3\mu_k}{K_k}\mathcal{S}_k + \frac{\imath u_0\Delta_1\Delta_2\Delta_3\delta_2\mathbf{k}_k}{4J_1J_2J_3K_k}, \\ [\mathcal{S}_k, \mathcal{S}_0]_- &= \imath(K_i - K_j)[\mathcal{S}_i, \mathcal{S}_j]_+ + 2\imath(\mathbf{k}_i\mathcal{S}_j - \mathbf{k}_j\mathcal{S}_i). \end{aligned} \tag{2.74}$$

On the other hand, (2.73) implies

$$\begin{aligned} \mathcal{D}_k = -K_k, \qquad 2(\mathbf{a}+\mathbf{b}\mathbf{j}_k) = -u_0\Delta_1\Delta_2\Delta_3\mu_k, \\ 2\mathbf{b} = -\frac{u_0\Delta_1\Delta_2\Delta_3\delta_2}{4J_1J_2J_3}. \end{aligned} \tag{2.75}$$

Thus (2.74) coincides with (2.69).

## 3 Representations by difference operators

In this section, we discuss representations of the algebras introduced above. For a representation $\pi$ by difference operators, the image of an abstract algebra element $A$ is denoted by $\hat{A} = \pi(A)$. This convention applies to the generators, their products and the Casimir elements. We write $\Omega_i$ for a central element, $\hat{\Omega}_i$ for its operator image and $\delta_i$ for its scalar value, so that $\hat{\Omega}_i = \delta_i \,\mathrm{Id}$ whenever the center acts by scalars.

### 3.1 Original representation of the Sklyanin algebra $\mathcal{A}$

As Sklyanin showed in his pioneering work [38], the algebra (2.1)–(2.2) admits a representation by difference operators acting on the space of meromorphic functions:

$$\begin{gathered}\hat{S}_\alpha(z,\hbar,\eta|\tau) := \hat{S}_\alpha = \imath^{\delta_{\alpha,2}} \frac{\theta_{\alpha+1}(\hbar|\tau)}{\theta_1(2z|\tau)} \Big( \theta_{\alpha+1}(2z-2\eta|\tau) e^{\hbar\partial_z} - \theta_{\alpha+1}(-2z-2\eta|\tau) e^{-\hbar\partial_z} \Big) = \\ = \imath^{\delta_{\alpha,2}} \frac{\theta_{\alpha+1}(\hbar)\theta_{\alpha+1}(-2\eta)}{\theta_1'(0)} \Big( \varphi_\alpha(2z,-2\eta+\omega_\alpha) e^{\hbar\partial_z} + \varphi_\alpha(-2z,-2\eta+\omega_\alpha) e^{-\hbar\partial_z} \Big)\end{gathered} \tag{3.1}$$

for $\alpha \in \{0,1,2,3\}$ and

$$J_1(\hbar) = \frac{\theta_2(2\hbar|\tau)\theta_2(0|\tau)}{\theta_2^2(\hbar|\tau)}, \quad J_2(\hbar) = \frac{\theta_3(2\hbar|\tau)\theta_3(0|\tau)}{\theta_3^2(\hbar|\tau)}, \quad J_3(\hbar) = \frac{\theta_4(2\hbar|\tau)\theta_4(0|\tau)}{\theta_4^2(\hbar|\tau)}, \tag{3.2}$$

where $\theta_\alpha(z|\tau)$ are the Jacobi theta functions (A.1) with modular parameter $\tau$, and the function $\varphi_\alpha(z,u)$ is defined in (A.26). In what follows, we omit the dependence on $\tau$. The shift operator $e^{\hbar\partial_z}$ acts by

$$e^{\hbar\partial_z} f(z) = f(z+\hbar)\,. \tag{3.3}$$

Introduce

$$\begin{gathered}\wp_i^\hbar = \wp(\hbar+\omega_i) - \wp(\hbar)\,, \\ \mathbf{K}_i(\hbar) = E_1(\hbar+\omega_i) - E_1(\hbar) - E_1(\omega_i)\,, \\ \omega_1 = \frac{1}{2}\,, \qquad \omega_2 = \frac{1+\tau}{2}\,, \qquad \omega_3 = \frac{\tau}{2}\,.\end{gathered} \tag{3.4}$$

The functions $E_1$ and $\wp$ are defined in (A.5). The notation $\mathbf{K}_i(\hbar)$ distinguishes these functions from the constants $K_i$ of the algebra (2.51). The addition formula for the Weierstrass zeta function gives

$$\mathbf{K}_i(\hbar) = \frac{\wp'(\hbar)}{2\big(\wp(\hbar)-e_i\big)}\,, \qquad e_i = \wp(\omega_i)\,. \tag{3.5}$$

Together with the Weierstrass identities

$$\begin{gathered}\wp'(\hbar)^2 = 4\prod_{\ell=1}^{3}\big(\wp(\hbar)-e_\ell\big)\,, \qquad e_1+e_2+e_3=0\,, \\ \wp(\hbar+\omega_i) = e_i + \frac{(e_i-e_j)(e_i-e_k)}{\wp(\hbar)-e_i}\,,\end{gathered} \tag{3.6}$$

this implies the exact quadratic relation

$$\wp_i^\hbar = \mathbf{K}_i(\hbar)\big(\mathbf{K}_i(\hbar) - \mathbf{K}_j(\hbar) - \mathbf{K}_k(\hbar)\big)\,. \tag{3.7}$$

Consequently,

$$\frac{\mathbf{K}_j(\hbar)-\mathbf{K}_i(\hbar)}{\mathbf{K}_k(\hbar)}=-\frac{\wp^\hbar_j-\wp^\hbar_i}{\wp^\hbar_k}\,, \tag{3.8}$$

as noted in the appendix of [20]. Here and below, $(i,j,k)$ is a cyclic permutation of $(1,2,3)$. In the normalization (3.2), the parameters $J_i$ satisfy

$$J_i(\hbar)=\frac{\theta_{i+1}(2\hbar)\theta_{i+1}(0)}{\theta_{i+1}(\hbar)^2}=-\frac{\theta_1(\hbar)^2}{\theta_1'(0)^2}\,\wp^\hbar_i\,. \tag{3.9}$$

The auxiliary parameters in (2.7) can therefore be chosen as

$$\mathcal{K}_i(\hbar)=\imath\,\frac{\theta_1(\hbar)}{\theta_1'(0)}\,\mathbf{K}_i(\hbar)\,. \tag{3.10}$$

Equations (3.7)–(3.10) give

$$\begin{gathered}J_i(\hbar)=\mathcal{K}_i(\hbar)\big(\mathcal{K}_i(\hbar)-\mathcal{K}_j(\hbar)-\mathcal{K}_k(\hbar)\big)\,,\\ \frac{\mathcal{K}_j(\hbar)-\mathcal{K}_i(\hbar)}{\mathcal{K}_k(\hbar)}=-\frac{J_j(\hbar)-J_i(\hbar)}{J_k(\hbar)}\,.\end{gathered} \tag{3.11}$$

Thus (2.7) holds. The structure constants $J_{\alpha\beta}=(J_\beta-J_\alpha)/J_\gamma$ can be represented in the following form

$$J_{12}(\hbar)=\frac{\theta_1^2(\hbar)\theta_4^2(\hbar)}{\theta_2^2(\hbar)\theta_3^2(\hbar)},\quad J_{23}(\hbar)=\frac{\theta_1^2(\hbar)\theta_2^2(\hbar)}{\theta_3^2(\hbar)\theta_4^2(\hbar)},\quad J_{31}(\hbar)=-\frac{\theta_1^2(\hbar)\theta_3^2(\hbar)}{\theta_2^2(\hbar)\theta_4^2(\hbar)}\,. \tag{3.12}$$

In this representation, the Casimir elements (2.3) take the following values:

$$\hat\Omega_1=\sum_{a=0}^3\hat S_a^2=\delta_1(\hbar,\eta)=4\theta_1^2(2\eta+\hbar)\,. \tag{3.13}$$

$$\hat\Omega_2=\sum_{a=1}^3 J_a\hat S_a^2=\delta_2(\hbar,\eta)=4\theta_1(2\eta)\theta_1(2\eta+2\hbar)\,. \tag{3.14}$$

The Sklyanin algebra also has the following representation [29]:

$$\begin{gathered}\hat S^+_\alpha(z,\hbar,\eta|\tau)=:\hat S^+_\alpha=\imath^{\delta_{\alpha,2}}\frac{\theta_{\alpha+1}(\hbar|\tau)}{\theta_1(2z|\tau)}\,\theta_{\alpha+1}(2z-2\eta|\tau)e^{\hbar\partial_z}=\\ =\imath^{\delta_{\alpha,2}}\frac{\theta_{\alpha+1}(\hbar)\theta_{\alpha+1}(-2\eta)}{\theta_1'(0)}\,\varphi_\alpha(2z,-2\eta+\omega_\alpha)e^{\hbar\partial_z}\,,\quad \alpha\in\{0,1,2,3\}\,.\end{gathered} \tag{3.15}$$

In this case, the structure constants are those of (3.2) and (3.12), and both Casimir elements vanish:

$$\delta_1=\delta_2=0\,. \tag{3.16}$$

The operators (3.15) define a representation of the Sklyanin algebra on the infinite-dimensional space of meromorphic functions.

### 3.2 Representations of $\mathcal{A}(\mathcal{A}_2)$

Applying the representation to the quadratic elements constructed in Section 2 gives

$$\hat{Z}_0^L = \hat{X}_0(\hat{U}), \qquad \hat{Z}_i^L = \frac{1}{2}\sqrt{\frac{\Delta_j\Delta_k}{J_jJ_k}}\hat{X}_i(\hat{U})\,, \tag{3.17}$$

or

$$\hat{Z}_0^R = \hat{Y}_0(\hat{U}), \qquad \hat{Z}_i^R = \frac{1}{2}\sqrt{\frac{\Delta_j\Delta_k}{J_jJ_k}}\hat{Y}_i(\hat{U})\,. \tag{3.18}$$

These operators satisfy the relations (2.50). Here $\hat{X}_\alpha(\hat{U}) = \hat{U}\hat{S}_\alpha$ and $\hat{Y}_\alpha(\hat{U}) = \hat{S}_\alpha\hat{U}$ are the images of $X_\alpha(U)$ and $Y_\alpha(U)$ from (2.17) and (2.20). If $\hat{U} = \sum\limits_{\alpha=0}^{3} u_\alpha\hat{S}_\alpha$ is invertible, then $\hat{Z}_\alpha^R$ and $\hat{Z}_\alpha^L$ are conjugate:

$$\hat{Z}_\alpha^R = \hat{U}^{-1}\hat{Z}_\alpha^L\hat{U}\,. \tag{3.19}$$

Otherwise, $\hat{Z}_\alpha^R$ and $\hat{Z}_\alpha^L$ are intertwined:

$$\hat{U}\hat{Z}_\alpha^R = \hat{Z}_\alpha^L\hat{U}\,. \tag{3.20}$$

For the representation (3.1)–(3.2) and $\alpha \in \{0,1,2,3\}$, we have

$$\begin{gathered}\hat{Z}_\alpha^L = \imath^{\delta_{\alpha,2}}\frac{\theta_1(\hbar)\theta_{\alpha+1}(2\hbar)\theta_{\alpha+1}(-2\eta)}{\theta_1'(0)\theta_1(2\hbar)}\Big(\varphi_\alpha(2z+2\hbar,-2\eta+\omega_\alpha)v(z,-2\eta)e^{2\hbar\partial_z}+\\ +\varphi_\alpha(-2z-2\hbar,-2\eta+\omega_\alpha)v(z,-2\eta)+\big(z\to -z\big)\Big)\,,\end{gathered} \tag{3.21}$$

$$\begin{gathered}\hat{Z}_\alpha^R = \imath^{\delta_{\alpha,2}}\frac{\theta_1(\hbar)\theta_{\alpha+1}(2\hbar)\theta_{\alpha+1}(-2\eta)}{\theta_1'(0)\theta_1(2\hbar)}\Big(\varphi_\alpha(2z,-2\eta+\omega_\alpha)v(z+\hbar,-2\eta)e^{2\hbar\partial_z}+\\ +\varphi_\alpha(-2z,-2\eta+\omega_\alpha)v(z-\hbar,-2\eta)+\big(z\to -z\big)\Big)\,.\end{gathered} \tag{3.22}$$

The constants $\nu_\alpha$ in the definition (A.32) of the function $v(z,u)$ are expressed in terms of the representation parameters as follows:

$$\nu_\alpha = \imath^{\delta_{\alpha,2}}\frac{\theta_{\alpha+1}(\hbar)\theta_{\alpha+1}(-2\eta)}{\theta_1'(0)}\,u_\alpha\,. \tag{3.23}$$

For the operators (3.21) and (3.22), the commutation relations (2.50) take the form

$$\begin{gathered}[\hat{Z}_i^{L,R},\hat{Z}_j^{L,R}]_- = \imath\,[\hat{Z}_0^{L,R},\hat{Z}_k^{L,R}]_+ + \imath u_0\mu_k\hat{Z}_k^{L,R} - \imath u_k\lambda_k\hat{Z}_0^{L,R}\,,\\ [\hat{Z}_k^{L,R},\hat{Z}_0^{L,R}]_- = \imath\,\frac{K_j-K_i}{K_k}\,[\hat{Z}_i^{L,R},\hat{Z}_j^{L,R}]_+ - 2\imath\frac{u_i}{K_k}\,\rho_i\mu_j\hat{Z}_j^{L,R} + 2\imath\frac{u_j}{K_k}\,\rho_j\mu_i\hat{Z}_i^{L,R}\,.\end{gathered} \tag{3.24}$$

The parameters $K_i$ are obtained from the original triple $J_i$ by

$$K_i(\hbar) = \Delta_i(\hbar)^2\,, \qquad \Delta_i(\hbar) = J_j(\hbar) + J_k(\hbar) - J_i(\hbar)\,. \tag{3.25}$$

Their theta-function expressions and the remaining coefficients are

$$\sqrt{K_k} = \Delta_k = J_i + J_j - J_k = \frac{\theta_{i+1}(2\hbar)\theta_{j+1}(2\hbar)\theta_{i+1}(0)\theta_{j+1}(0)\theta_{k+1}^2(0)}{\theta_2^2(\hbar)\theta_3^2(\hbar)\theta_4^2(\hbar)} = J_i(\hbar)J_j(\hbar)\frac{\theta_{k+1}^2(0)}{\theta_{k+1}^2(\hbar)}\,,$$

$$\begin{aligned}
\rho_i &= \sqrt{J_j J_k \Delta_j \Delta_k} = J_i(\hbar)J_j(\hbar)J_k(\hbar)\frac{\theta_{j+1}(0)\theta_{k+1}(0)}{\theta_{j+1}(\hbar)\theta_{k+1}(\hbar)}\,,\\
\mu_k &= \frac{\rho_k\sqrt{K_k}}{\rho_i\rho_j}\,\delta_2 - \delta_1 = -\frac{4\theta_1^2(\hbar)\theta_{k+1}(2\eta+2\hbar)\theta_{k+1}(2\eta)}{\theta_{k+1}(2\hbar)\theta_{k+1}(0)}\,,\\
\lambda_k &= \frac{\sqrt{K_iK_jK_k}}{2\rho_i\rho_j}\,\delta_2 = \frac{2\theta_1(2\eta)\theta_1(2\eta+2\hbar)\theta_{j+1}(0)\theta_{i+1}(0)}{\theta_{j+1}(\hbar)\theta_{i+1}(\hbar)}\,.
\end{aligned} \tag{3.26}$$

In calculating (3.21), (3.22) and (3.26), we use the relations (A.11)–(A.23) and (A.10). The duplication formulae applied to the first line of (3.26), together with (3.5), also give

$$\begin{aligned}
K_i(\hbar) &= 16\,\frac{\theta_1(\hbar)^4}{\theta_1'(0)^4}\,\big(\wp(2\hbar)-e_j\big)\big(\wp(2\hbar)-e_k\big) =\\
&= 8\,\frac{\theta_1(\hbar)^4}{\theta_1'(0)^4}\,\wp'(2\hbar)\,\mathbf{K}_i(2\hbar)\,.
\end{aligned} \tag{3.27}$$

The common prefactor cancels in the structure constants. Using (3.8) and (3.9) at $2\hbar$ yields

$$\begin{aligned}
&\frac{K_j(\hbar)-K_i(\hbar)}{K_k(\hbar)} = \frac{\mathbf{K}_j(2\hbar)-\mathbf{K}_i(2\hbar)}{\mathbf{K}_k(2\hbar)} =\\
&= -\frac{\wp_j^{2\hbar}-\wp_i^{2\hbar}}{\wp_k^{2\hbar}} = -\frac{J_j(2\hbar)-J_i(2\hbar)}{J_k(2\hbar)} = -J_{ij}(2\hbar)\,.
\end{aligned} \tag{3.28}$$

Thus $\mathcal{K}_i(\hbar)$ and $K_i(\hbar)$ are proportional, with factors independent of $i$, to $\mathbf{K}_i(\hbar)$ and $\mathbf{K}_i(2\hbar)$, respectively. The passage to the quadratic generators doubles the argument of the elliptic parametrization. Formula (3.27) retains the exact normalization $K_i = \Delta_i^2$.

By Proposition 2.2, the Casimir elements (2.52) and (2.54) vanish for the operators (3.21) and (3.22):

$$\delta_1^{\mathcal{Z}} = \delta_2^{\mathcal{Z}} = 0\,. \tag{3.29}$$

### 3.3 One-parameter deformation of the representations of $\mathcal{A}(\mathcal{A}_2)$

It is natural to expect a representation with nonzero Casimir values. Such a representation can be obtained by a slight deformation of (3.22).

**Proposition 3.1** *The following operators*

$$\begin{aligned}
\hat{Z}_\alpha^\xi = i^{\delta_{\alpha,2}}\frac{\theta_1(\hbar)\theta_{\alpha+1}(2\hbar)\theta_{\alpha+1}(-2\eta)}{\theta_1'(0)\theta_1(2\hbar)}\Big(&\varphi_\alpha(2z,-2\eta+\omega_\alpha)v(z+\hbar,-2\xi)e^{2\hbar\partial_z}+\\
&+\varphi_\alpha(-2z,-2\eta+\omega_\alpha)v(z-\hbar,-2\eta)+\big(z\to -z\big)\Big)\,, \qquad \alpha\in\{0,1,2,3\}
\end{aligned} \tag{3.30}$$

*satisfy the defining relations of the algebra (2.51) with the structure constants and parameters defined in (3.26) and (3.28). In this representation, the Casimir elements (2.52) and (2.54) take the following values:*

$$\delta_1^{\mathcal{Z}} = \delta_1(2\hbar,\eta)\frac{\theta_1^2(\hbar)}{\theta_1^2(2\hbar)}\sum_{\alpha=0}^{3}\nu_\alpha^2\Big(\wp(2\xi+\omega_\alpha)-\wp(2\eta+\omega_\alpha)\Big)\,, \tag{3.31}$$

$$\delta_2^{\mathcal{Z}} = -64\,\delta_2(2\hbar,\eta)\frac{\theta_1^{10}(\hbar)\theta_1^2(4\hbar)}{\theta_1^{12}(2\hbar)}\sum_{\alpha=0}^{3}\nu_\alpha^2\Big(\wp(2\xi+\omega_\alpha)-\wp(2\eta+\omega_\alpha)\Big)\,, \tag{3.32}$$

*where*

$$\nu_\alpha = \imath^{\delta_{\alpha,2}}\frac{\theta_{\alpha+1}(\hbar)\theta_{\alpha+1}(-2\eta)}{\theta_1'(0)}\,u_\alpha\,. \tag{3.33}$$

*These are the parameters of the function $v(z,\eta)$. When $\xi=\eta$, we recover the representation $\hat{Z}_\alpha^R$ in (3.22).*

**Proof.** The commutation relations for $\hat{Z}_\alpha^\xi$ follow from (3.24) for $\hat{Z}_\alpha^R$ and the Sklyanin relations (2.1), (2.2) for $\hat{S}_\alpha(z,2\hbar,\eta|\tau)$. Indeed, (3.30) differs from (3.22) only by replacing $\eta$ with $\xi$ in the factors $v(\pm z+\hbar,-2\eta)$ multiplying $e^{\pm 2\hbar\partial_z}$. Thus $\hat{Z}_\alpha^\eta=\hat{Z}_\alpha^R$, while the terms without shifts are independent of $\xi$. The coefficient of $e^{\pm 2\hbar\partial_z}$ in $\hat{Z}_\alpha^\xi$ equals the corresponding coefficient in $\hat{S}_\alpha(z,2\hbar,\eta|\tau)$ multiplied by

$$\frac{\theta_1(\hbar)}{\theta_1(2\hbar)}\,v(\pm z+\hbar,-2\xi)\,. \tag{3.34}$$

The relations (3.24) have degree at most two, and their scalar coefficients (3.26), (3.28) do not depend on $\xi$. All contributions to the coefficients of $e^{\pm 2\hbar\partial_z}$ have the common factor $v(\pm z+\hbar,-2\xi)$. Similarly, the coefficients of $e^{\pm 4\hbar\partial_z}$ have the common factor

$$v(\pm z+\hbar,-2\xi)v(\pm z+3\hbar,-2\xi)\,. \tag{3.35}$$

Once these factors are removed, the remaining expressions are independent of $\xi$. They vanish because the relations hold for $\hat{Z}_\alpha^R$. The argument applies first to generic parameters and then extends meromorphically. Hence all coefficients of nonzero shifts vanish.

In the terms without shifts, only products of opposite shifts depend on $\xi$. Their common factors are

$$v(\pm z+\hbar,-2\xi)v(\mp z-\hbar,-2\xi)\,. \tag{3.36}$$

By (A.36) or, equivalently, by (A.34) and (A.35),

$$v(\pm z+\hbar,-2\xi)v(\mp z-\hbar,-2\xi)=\sum_{\alpha=0}^{3}\Big(\nu_\alpha^2\wp(-2\xi+\omega_\alpha)-\check{\nu}_\alpha^2\wp(\pm z+\hbar+\omega_\alpha)\Big)\,. \tag{3.37}$$

Consequently, for any $\alpha,\beta$, the term without a shift in $\hat{Z}_\alpha^\xi\hat{Z}_\beta^\xi-\hat{Z}_\alpha^R\hat{Z}_\beta^R$ is the term without a shift in $\hat{S}_\alpha(z,2\hbar,\eta|\tau)\hat{S}_\beta(z,2\hbar,\eta|\tau)$ multiplied by

$$\frac{\theta_1^2(\hbar)}{\theta_1^2(2\hbar)}\sum_{\alpha=0}^{3}\nu_\alpha^2\Big(\wp(-2\xi+\omega_\alpha)-\wp(-2\eta+\omega_\alpha)\Big)\,. \tag{3.38}$$

The terms without shifts in the linear parts of (3.24) do not change. Thus the remaining differences between the relations for $\hat{Z}_\alpha^\xi$ and $\hat{Z}_\alpha^R$ are proportional to the terms without shifts in

$$[\hat{S}_i(z,2\hbar,\eta|\tau),\hat{S}_j(z,2\hbar,\eta|\tau)]_- - \imath[\hat{S}_0(z,2\hbar,\eta|\tau),\hat{S}_k(z,2\hbar,\eta|\tau)]_+\,, \tag{3.39}$$

and

$$[\hat{S}_k(z,2\hbar,\eta|\tau),\hat{S}_0(z,2\hbar,\eta|\tau)]_- - \imath\frac{K_j-K_i}{K_k}[\hat{S}_i(z,2\hbar,\eta|\tau),\hat{S}_j(z,2\hbar,\eta|\tau)]_+ \,. \tag{3.40}$$

Both expressions vanish by (2.1), (2.2) and (3.28). This proves the required commutation relations.

For the Casimir elements (2.52), (2.54), the same factorization shows that all coefficients of nonzero shifts vanish. Their values for $\hat{Z}^R_\alpha$ are zero by (3.29). It therefore remains to evaluate their homogeneous quadratic parts on $\hat{S}_\alpha(z,2\hbar,\eta|\tau)$ and multiply by (3.38). For the first Casimir, (3.13) gives

$$\delta_1^{\mathcal{Z}} = \delta_1(2\hbar,\eta)\frac{\theta_1^2(\hbar)}{\theta_1^2(2\hbar)}\sum_{\alpha=0}^3 \nu_\alpha^2\Big(\wp(2\xi+\omega_\alpha)-\wp(2\eta+\omega_\alpha)\Big)\,. \tag{3.41}$$

Here the evenness and double periodicity of $\wp$ have been used. For the second Casimir, the duplication formulae (A.11)–(A.23), together with (3.26), give

$$K_i(K_i-K_j-K_k) = -64\,\frac{\theta_1^8(\hbar)\theta_1^2(4\hbar)}{\theta_1^{10}(2\hbar)}\,J_i(2\hbar)\,. \tag{3.42}$$

Using (3.14) and multiplying by (3.38) yields

$$\delta_2^{\mathcal{Z}} = -64\,\delta_2(2\hbar,\eta)\frac{\theta_1^{10}(\hbar)\theta_1^2(4\hbar)}{\theta_1^{12}(2\hbar)}\sum_{\alpha=0}^3 \nu_\alpha^2\Big(\wp(2\xi+\omega_\alpha)-\wp(2\eta+\omega_\alpha)\Big)\,. \tag{3.43}$$

This proves (3.31), (3.32). ■

Next we consider the representations that our construction gives in particular cases of the algebra (2.50).

### 3.4 Back to the Sklyanin algebra $\mathcal{A}$

As noted above, if $\delta_1=\delta_2=0$, then the algebra (2.50) reduces to the ordinary four-dimensional Sklyanin algebra. For the representation (3.1)–(3.2), setting the Casimir values (3.13) and (3.14) to zero is equivalent to imposing both $\hbar\in\mathbb{Z}+\tau\mathbb{Z}$ and $2\eta\in\mathbb{Z}+\tau\mathbb{Z}$. This is a degenerate lattice-shift case, which will not be considered here. However, the nontrivial representation (3.15) has zero Casimir values. For these generators, the elements (3.17) and (3.18) reduce to

$$\hat{Z}^L_\alpha = \imath^{\delta_{\alpha,2}}\frac{\theta_1(\hbar)\theta_{\alpha+1}(2\hbar)\theta_{\alpha+1}(-2\eta)}{\theta_1'(0)\theta_1(2\hbar)}\,\varphi_\alpha(2z+2\hbar,-2\eta+\omega_\alpha)v(z,-2\eta)e^{2\hbar\partial_z}\,, \tag{3.44}$$

$$\hat{Z}^R_\alpha = \imath^{\delta_{\alpha,2}}\frac{\theta_1(\hbar)\theta_{\alpha+1}(2\hbar)\theta_{\alpha+1}(-2\eta)}{\theta_1'(0)\theta_1(2\hbar)}\,\varphi_\alpha(2z,-2\eta+\omega_\alpha)v(z+\hbar,-2\eta)e^{2\hbar\partial_z}\,. \tag{3.45}$$

At generic points, $\hat{Z}^R_\alpha$ and $\hat{Z}^L_\alpha$ are intertwined by the operator $\hat{U}=\pi(U)$:

$$\hat{U}\hat{Z}^R_\alpha = \hat{Z}^L_\alpha\hat{U}\,. \tag{3.46}$$

If the operator $\hat{U}$ is invertible, then $\hat{Z}^R_\alpha$ and $\hat{Z}^L_\alpha$ are gauge equivalent:

$$\hat{Z}^R_\alpha = \hat{U}^{-1}\hat{Z}^L_\alpha\hat{U},\quad \hat{U}=v(z,-2\eta)e^{\hbar\partial_z},\quad \hat{U}^{-1}=\frac{1}{v(z-\hbar,-2\eta)}\,e^{-\hbar\partial_z}\,. \tag{3.47}$$

Note that

$$\hat{Z}^R_\alpha = \frac{\theta_1(\hbar)}{\theta_1(2\hbar)}\,v(z+\hbar,-2\eta)\hat{S}^+_\alpha(z,2\hbar,\eta|\tau)\,, \tag{3.48}$$

so we obtain a representation of type (3.15), but now all generators have the additional dynamical factor $\theta_1(\hbar)v(z+\hbar,-2\eta)/\theta_1(2\hbar)$. One might wonder how important this factor is. In fact, for any representation of the form $\hat{S}_\alpha = f_\alpha(z)e^{\hbar\partial_z}$, the commutation relations (2.1)–(2.2) reduce to

$$
\begin{aligned}
& f_i(z)f_j(z+\hbar) - f_j(z)f_i(z+\hbar) = \imath\big(f_0(z)f_k(z+\hbar) + f_k(z)f_0(z+\hbar)\big)\,, \\
& f_0(z)f_k(z+\hbar) - f_k(z)f_0(z+\hbar) = \imath\frac{J_j - J_i}{J_k}\,\big(f_i(z)f_j(z+\hbar) + f_j(z)f_i(z+\hbar)\big)\,.
\end{aligned}
\tag{3.49}
$$

This system of functional relations is invariant under the transformation

$$
f_\alpha(z) \to g(z)f_\alpha(z)\,, \tag{3.50}
$$

where $g(z)$ is an arbitrary meromorphic function. Thus, the factor $\theta_1(\hbar)v(z+\hbar,-2\eta)/\theta_1(2\hbar)$ in (3.48) can be removed using the rescaling symmetry (3.50). We conclude that, starting from (3.15), our construction gives the same representation as (3.15) with $\hbar$ doubled.

Now let us turn to the second method for obtaining the Sklyanin algebra from (2.50). In this case, we shift the generator $\hat{Z}_0$

$$
\hat{Z}'_0 = \hat{Z}_0 - C_0\,, \quad C_0 = -\frac{1}{2}\,u_0\mu_k\,, \tag{3.51}
$$

and set $u_k = 0$ for $k \in \{1,2,3\}$. For $C_0$ to be well defined, we set $\mu_1 = \mu_2 = \mu_3$. By (2.71), $\mu_1 = \mu_2 = \mu_3$ if and only if $\delta_2 = 0$ or $J_1 = J_2 = J_3$. In the latter case, (3.12) implies $\theta_1(\hbar) = 0$, which is excluded in the nondegenerate elliptic parametrization considered here. We therefore impose

$$
\delta_2 = 4\theta_1(2\eta)\theta_1(2\eta + 2\hbar) = 0\,. \tag{3.52}
$$

The last equality holds when $\eta = n\omega_1 + m\omega_3$ or $\eta = -\hbar + k\omega_1 + l\omega_3$ for $m,\ n,\ k,\ l\ \in \mathbb{Z}$. In the following, for the sake of brevity, we will consider only the cases $n = m = 0$ and $k = l = 0$.

- $\eta = 0$

For the operators (3.30), there are no nontrivial degenerations at $\eta = 0$ (for fixed generic $\xi$, the shifted generators obtained from (3.30) tend to zero as $\eta \to 0$). For (3.21) and (3.22), this limit gives

$$
\hat{Z}'^{L}_0 = u_0\theta_1^2(\hbar)\Big(e^{2\hbar\partial_z} + e^{-2\hbar\partial_z}\Big) = u_0\frac{\theta_1^2(\hbar)}{\theta_1(2\hbar)}\,\hat{S}_0(z,2\hbar,0|\tau)\,, \tag{3.53}
$$

$$
\begin{aligned}
\hat{Z}'^{L}_k = \imath^{\delta_{k,2}}u_0\frac{\theta_1^2(\hbar)\theta_{k+1}(2\hbar)\theta_{k+1}(0)}{\theta'_1(0)\theta_1(2\hbar)}\,\Big(&\varphi_k(2z+2\hbar)e^{2\hbar\partial_z} + \varphi_k(-2z+2\hbar)e^{-2\hbar\partial_z} + \\
&+\varphi_k(-2z-2\hbar) + \varphi_k(2z-2\hbar)\Big)\,,
\end{aligned}
\tag{3.54}
$$

$$
\hat{Z}'^{R}_0 = u_0\theta_1^2(\hbar)\Big(e^{2\hbar\partial_z} + e^{-2\hbar\partial_z}\Big) = u_0\frac{\theta_1^2(\hbar)}{\theta_1(2\hbar)}\,\hat{S}_0(z,2\hbar,0|\tau)\,, \tag{3.55}
$$

$$
\begin{aligned}
\hat{Z}'^{R}_k &= \imath^{\delta_{k,2}}u_0\frac{\theta_1^2(\hbar)\theta_{k+1}(2\hbar)\theta_{k+1}(0)}{\theta'_1(0)\theta_1(2\hbar)}\,\Big(\varphi_k(2z)e^{2\hbar\partial_z} + \varphi_k(-2z)e^{-2\hbar\partial_z} + \varphi_k(-2z) + \varphi_k(2z)\Big) = \\
&= \imath^{\delta_{k,2}}u_0\frac{\theta_1^2(\hbar)\theta_{k+1}(2\hbar)\theta_{k+1}(0)}{\theta'_1(0)\theta_1(2\hbar)}\,\Big(\varphi_k(2z)e^{2\hbar\partial_z} + \varphi_k(-2z)e^{-2\hbar\partial_z}\Big) = u_0\frac{\theta_1^2(\hbar)}{\theta_1(2\hbar)}\,\hat{S}_k(z,2\hbar,0|\tau)\,,
\end{aligned}
\tag{3.56}
$$

where the functions $\varphi_k(z)$ are defined in (A.28). In this case, the operator

$$
\hat{U} = u_0\theta_1(\hbar)(e^{\hbar\partial_z} + e^{-\hbar\partial_z})\,, \tag{3.57}
$$

is not invertible on the whole space of meromorphic functions, and the operators (3.53)–(3.54) and (3.55)–(3.56) are intertwined as in (3.20) or (3.46). The Casimir values in this case can be expressed in terms of (3.13) and (3.14) (throughout this subsection, $\delta_2'^{\mathcal{Z}}$ denotes the scalar value of $\sum_{k=1}^3 J_k(2\hbar)(\hat{Z}_k')^2$) :

$$\delta_1'^{\mathcal{Z}} = u_0^2 \frac{\theta_1^4(\hbar)}{\theta_1^2(2\hbar)}\, \delta_1(2\hbar, 0) = 4u_0^2\theta_1^4(\hbar)\,, \tag{3.58}$$

$$\delta_2'^{\mathcal{Z}} = u_0^2 \frac{\theta_1^4(\hbar)}{\theta_1^2(2\hbar)}\, \delta_2(2\hbar, 0) = 0\,. \tag{3.59}$$

- $\eta = -\hbar$

In this case, the operators (3.30) reduce to

$$\hat{Z}_0'^{\xi} = -u_0 \frac{\theta_1^2(\hbar)\theta_1(2\hbar)\theta_1(2\hbar - 2\xi)}{\theta_1'(0)\theta_1(2\xi)} \Big(\phi(2z, 2\hbar - 2\xi)e^{2\hbar\partial_z} + \phi(-2z, 2\hbar - 2\xi)e^{-2\hbar\partial_z}\Big) = \tag{3.60}$$

$$= -u_0 \frac{\theta_1^2(\hbar)}{\theta_1(2\xi)}\, \hat{S}_0(z, 2\hbar, \xi - \hbar)\,,$$

$$\hat{Z}_k'^{\xi} = \imath^{\delta_{k,2}} u_0 \frac{\theta_1^2(\hbar)\theta_{k+1}^2(2\hbar)}{\theta_1'(0)^2} \Big(\varphi_k(2z, 2\hbar + \omega_k)\phi(2z + 2\hbar, -2\xi)e^{2\hbar\partial_z} + \tag{3.61}$$

$$+\varphi_k(-2z, 2\hbar + \omega_k)\phi(2z - 2\hbar, 2\hbar) + \big(z \to -z\big)\Big)\,.$$

The Casimir elements for these generators take the values

$$\delta_1'^{\mathcal{Z}} = 4u_0^2\theta_1^4(\hbar)\,, \tag{3.62}$$

$$\delta_2'^{\mathcal{Z}} = -4u_0^2 \frac{\theta_1^4(\hbar)\theta_1^2(2\hbar)}{\theta_1'(0)^2} \big(\wp(2\xi) - \wp(2\hbar)\big)\,. \tag{3.63}$$

The operators (3.21) and (3.22) become

$$\hat{Z}_0'^{L} = u_0 \frac{\theta_1^2(\hbar)\theta_1(4\hbar)}{\theta_1'(0)} \Big(\phi(2z, 4\hbar)e^{2\hbar\partial_z} + \phi(-2z, 4\hbar)e^{-2\hbar\partial_z}\Big) = u_0 \frac{\theta_1^2(\hbar)}{\theta_1(2\hbar)}\, \hat{S}_0(z, 2\hbar, -2\hbar|\tau)\,, \tag{3.64}$$

$$\hat{Z}_k'^{L} = \imath^{\delta_{k,2}} u_0 \frac{\theta_1^2(\hbar)\theta_{k+1}^2(2\hbar)}{\theta_1'(0)^2} \Big(\varphi_k(2z + 2\hbar, 2\hbar + \omega_k)\phi(2z, 2\hbar)e^{2\hbar\partial_z} + \tag{3.65}$$

$$+\varphi_k(-2z - 2\hbar, 2\hbar + \omega_k)\phi(2z, 2\hbar) + \big(z \to -z\big)\Big) =$$

$$= \imath^{\delta_{k,2}} u_0 \frac{\theta_1^2(\hbar)\theta_{k+1}(2\hbar)\theta_{k+1}(4\hbar)}{\theta_1'(0)\theta_1(2\hbar)} \Big(\varphi_k(2z, 4\hbar + \omega_k)e^{2\hbar\partial_z} + \varphi_k(-2z, 4\hbar + \omega_k)e^{-2\hbar\partial_z}\Big) =$$

$$= u_0 \frac{\theta_1^2(\hbar)}{\theta_1(2\hbar)}\, \hat{S}_k(z, 2\hbar, -2\hbar|\tau)\,,$$

$$\hat{Z}_0'^{R} = u_0 \frac{\theta_1^2(\hbar)\theta_1(4\hbar)}{\theta_1'(0)} \Big(\phi(2z, 4\hbar)e^{2\hbar\partial_z} + \phi(-2z, 4\hbar)e^{-2\hbar\partial_z}\Big) = u_0 \frac{\theta_1^2(\hbar)}{\theta_1(2\hbar)}\, \hat{S}_0(z, 2\hbar, -2\hbar|\tau)\,, \tag{3.66}$$

$$\hat{Z}_k'^{R} = \imath^{\delta_{k,2}} u_0 \frac{\theta_1^2(\hbar)\theta_{k+1}^2(2\hbar)}{\theta_1'(0)^2} \Big(\varphi_k(2z, 2\hbar + \omega_k)\phi(2z + 2\hbar, 2\hbar)e^{2\hbar\partial_z} + \tag{3.67}$$

$$+\varphi_k(-2z,2\hbar+\omega_k)\phi(2z-2\hbar,2\hbar)+\big(z\to -z\big)\Big)\,,$$

and

$$\hat U\hat Z'^{R}_\alpha = u_0\frac{\theta_1^2(\hbar)}{\theta_1(2\hbar)}\,\hat S_\alpha(z,2\hbar,-2\hbar|\tau)\hat U\,,\quad \alpha=0,1,2,3\,, \tag{3.68}$$

$$\hat U = u_0\frac{\theta_1(\hbar)}{\theta_1(2z)}\left(\theta_1(2z+2\hbar)e^{\hbar\partial_z}-\theta_1(-2z+2\hbar)e^{-\hbar\partial_z}\right). \tag{3.69}$$

In all the remaining nontrivial cases listed above, the Casimir values are

$$\delta_1'^{\mathcal Z} = 4u_0^2\theta_1^4(\hbar)\,, \tag{3.70}$$

$$\delta_2'^{\mathcal Z} = 0\,. \tag{3.71}$$

### 3.5 $BC_1$ Sklyanin algebra

Here we construct representations of this algebra in terms of elliptic difference operators. In the previous sections, we showed that there are two ways, (2.65) and (2.66), to obtain the $BC_1$ Sklyanin algebra from the algebra (2.50).

The first is to deform the representations (3.53)–(3.67) of the ordinary Sklyanin algebra with scalar Casimir value $\delta_2=0$.

- $\eta=0$

In this limit ($\eta\to 0$ for fixed $\xi$), the operators (3.30) reduce to

$$\hat{\mathbb S}_0^\xi = \theta_1(\hbar)\sum_{j=1}^{3}\imath^{\delta_{j,2}}u_j\theta_{j+1}(\hbar)\Big(\frac{\theta_{j+1}(0)\theta_{j+1}(2z+2\hbar-2\xi)}{\theta_1(2z+2\hbar)\theta_{j+1}(2\xi)}\,e^{2\hbar\partial_z}- \tag{3.72}$$

$$-\frac{\theta_{j+1}(2z+2\hbar)}{\theta_1(2z+2\hbar)}\Big)+\big(z\to -z\big)\,,$$

$$\hat{\mathbb S}_k^\xi = \imath^{\delta_{k,2}}\frac{\theta_1(\hbar)\theta_{k+1}(2\hbar)\theta_{k+1}(0)}{\theta_1'(0)\theta_1(2\hbar)}\,\varphi_k(2z)\sum_{j=1}^{3}\imath^{\delta_{j,2}}u_j\theta_{j+1}(\hbar)\Big(\frac{\theta_{j+1}(0)\theta_{j+1}(2z+2\hbar-2\xi)}{\theta_1(2z+2\hbar)\theta_{j+1}(2\xi)}\,e^{2\hbar\partial_z}- \tag{3.73}$$

$$-\frac{\theta_{j+1}(2z+2\hbar)}{\theta_1(2z+2\hbar)}\Big)+\big(z\to -z\big)\,,$$

with the following values of the Casimir elements (2.14):

$$\delta_1^{\mathcal{B}C_1} = 4\frac{\theta_1^2(\hbar)}{\theta_1'(0)^2}\sum_{\alpha=1}^{3}(-1)^{\delta_{\alpha,2}}u_\alpha^2\theta_{\alpha+1}^2(\hbar)\theta_{\alpha+1}^2(0)\Big(\wp(2\xi+\omega_\alpha)-\wp(\omega_\alpha)\Big)\,, \tag{3.74}$$

$$\delta_2^{\mathcal{B}C_1} = 0\,. \tag{3.75}$$

For (3.21) and (3.22), this degeneration (here $\xi=\eta$ and $\eta\to 0$) leads to

$$\hat{\mathbb S}_0^L = \theta_1(\hbar)\sum_{\alpha=0}^{3}\imath^{\delta_{\alpha,2}}u_\alpha\frac{\theta_{\alpha+1}(\hbar)\theta_{\alpha+1}(2z)}{\theta_1(2z)}\,e^{2\hbar\partial_z}+\big(z\to -z\big)\,, \tag{3.76}$$

$$\hat{\mathbb{S}}_k^L = \imath^{\delta_{k,2}} \frac{\theta_1(\hbar)\theta_{k+1}(2\hbar)\theta_{k+1}(0)}{\theta_1'(0)\theta_1(2\hbar)}\, \varphi_k(2z+2\hbar) \sum_{\alpha=0}^{3} \imath^{\delta_{\alpha,2}} u_\alpha \frac{\theta_{\alpha+1}(\hbar)\theta_{\alpha+1}(2z)}{\theta_1(2z)} \Big(e^{2\hbar\partial_z} - 1\Big) + \tag{3.77}$$

$$+ \big(z \to -z\big)\,,$$

$$\hat{\mathbb{S}}_0^R = u_0 \theta_1^2(\hbar) e^{2\hbar\partial_z} + \theta_1(\hbar) \sum_{j=1}^{3} \imath^{\delta_{j,2}} u_j \frac{\theta_{j+1}(\hbar)\theta_{j+1}(2z+2\hbar)}{\theta_1(2z+2\hbar)} \Big(e^{2\hbar\partial_z} - 1\Big) + \big(z \to -z\big)\,, \tag{3.78}$$

$$\hat{\mathbb{S}}_k^R = \imath^{\delta_{k,2}} \frac{\theta_1(\hbar)\theta_{k+1}(2\hbar)\theta_{k+1}(0)}{\theta_1'(0)\theta_1(2\hbar)}\, \varphi_k(2z) \Big(u_0 \theta_1(\hbar) e^{2\hbar\partial_z} + \tag{3.79}$$

$$+ \sum_{j=1}^{3} \imath^{\delta_{j,2}} u_j \frac{\theta_{j+1}(\hbar)\theta_{j+1}(2z+2\hbar)}{\theta_1(2z+2\hbar)} \left(e^{2\hbar\partial_z} - 1\right)\Big) + \big(z \to -z\big)\,,$$

and the Casimir values coincide with those in (3.58)–(3.59):

$$\delta_1^{\mathcal{B}C_1} = 4u_0^2 \theta_1^4(\hbar)\,, \tag{3.80}$$

$$\delta_2^{\mathcal{B}C_1} = 0\,. \tag{3.81}$$

- $\eta = -\hbar$

In this limit, (3.30) yields one of the most remarkable representations of the quantum $BC_1$ Sklyanin algebra (2.11)–(2.14):

$$\hat{\mathbb{S}}_0^\xi = \frac{\theta_1(\hbar)\theta_1(2\hbar)}{\theta_1'(0)} \Big(\phi(2z, 2\hbar) v(z+\hbar, -2\xi) e^{2\hbar\partial_z} + \phi(-2z, 2\hbar) v(-z+\hbar, -2\xi) e^{-2\hbar\partial_z}\Big)\,, \tag{3.82}$$

$$\hat{\mathbb{S}}_k^\xi = \imath^{\delta_{k,2}} \frac{\theta_1(\hbar)\theta_{k+1}^2(2\hbar)}{\theta_1'(0)\theta_1(2\hbar)} \Big(\varphi_k(2z, 2\hbar+\omega_k) v(z+\hbar, -2\xi) e^{2\hbar\partial_z} + \tag{3.83}$$

$$+ \varphi_k(-2z, 2\hbar+\omega_k) v(z-\hbar, 2\hbar) + \big(z \to -z\big)\Big)\,.$$

The Casimir values are

$$\delta_1^{\mathcal{B}C_1} = 4u_0^2 \theta_1^4(\hbar)\,, \tag{3.84}$$

$$\delta_2^{\mathcal{B}C_1} = 256\, \frac{\theta_1^{10}(\hbar)\theta_1^2(4\hbar)}{\theta_1'(0)^2 \theta_1^{10}(2\hbar)} \sum_{\alpha=0}^{3} \imath^{2\delta_{\alpha,2}} u_\alpha^2 \theta_{\alpha+1}^2(\hbar) \theta_{\alpha+1}^2(2\hbar) \Big(\wp(2\xi+\omega_\alpha) - \wp(2\hbar+\omega_\alpha)\Big)\,. \tag{3.85}$$

In the following sections, we will see its connection with the quantum $BC_1$ van Diejen system and, in the classical limit, with the Zhukovsky-Volterra gyrostat. The degenerations of (3.21) and (3.22) have the form

$$\hat{\mathbb{S}}_0^L = \frac{\theta_1(\hbar)\theta_1(2\hbar)}{\theta_1'(0)} \Big(\phi(2z+2\hbar, 2\hbar) v(z, 2\hbar) e^{2\hbar\partial_z} + \tag{3.86}$$

$$+ \sum_{j=1}^{3} \imath^{\delta_{j,2}} u_j \frac{\theta_{j+1}(\hbar)\theta_{j+1}(0)}{\theta_1(2\hbar)}\, \varphi_j(2z+2\hbar) + \big(z \to -z\big)\Big)\,,$$

$$\hat{\mathbb{S}}_k^L = \imath^{\delta_{k,2}} \frac{\theta_1(\hbar)\theta_{k+1}(2\hbar)\theta_{k+1}(2\hbar)}{\theta_1'(0)\theta_1(2\hbar)} \Big(\varphi_k(2z+2\hbar, 2\hbar+\omega_k) v(z, 2\hbar) e^{2\hbar\partial_z} - \tag{3.87}$$

$$- \frac{\theta_{k+1}(2z)}{\theta_1(2z)} \sum_{j=1}^{3} \imath^{\delta_{j,2}} u_j \frac{\theta_{j+1}(\hbar)\theta_{j+1}(0)}{\theta_{k+1}(2\hbar)}\, \varphi_j(2z+2\hbar) + \big(z \to -z\big)\Big)\,,$$

$$\hat{\mathbb{S}}_0^R = \frac{\theta_1(\hbar)\theta_1(2\hbar)}{\theta_1'(0)}\Big(\phi(2z,2\hbar)v(z+\hbar,2\hbar)e^{2\hbar\partial_z} + \phi(-2z,2\hbar)v(-z+\hbar,2\hbar)e^{-2\hbar\partial_z}\Big)\,, \tag{3.88}$$

$$\hat{\mathbb{S}}_k^R = \imath^{\delta_{k,2}}\frac{\theta_1(\hbar)\theta_{k+1}^2(2\hbar)}{\theta_1'(0)\theta_1(2\hbar)}\Big(\varphi_k(2z,2\hbar+\omega_k)v(z+\hbar,2\hbar)e^{2\hbar\partial_z}+ \tag{3.89}$$

$$+\varphi_k(-2z,2\hbar+\omega_k)v(z-\hbar,2\hbar)+\big(z\to -z\big)\Big)\,.$$

$$\delta_1^{\mathcal{B}C_1} = 4u_0^2\theta_1^4(\hbar)\,, \tag{3.90}$$

$$\delta_2^{\mathcal{B}C_1} = 0\,. \tag{3.91}$$

The parameters $\nu_a$ of the function $v(z,\eta)$ in (3.82)–(3.89) are given by

$$\nu_\alpha = \imath^{\delta_{\alpha,2}}\frac{\theta_{\alpha+1}(\hbar)\theta_{\alpha+1}(2\hbar)}{\theta_1'(0)}\,u_\alpha\,. \tag{3.92}$$

In all cases, the structure constants $\widetilde{\nu}_k$ in (2.11)–(2.13) are given by

$$\widetilde{\nu}_k = -u_k\delta_1 = -4\theta_1^2(\hbar)u_k\,. \tag{3.93}$$

The second method gives a deformation of the trivial representation of the Sklyanin algebra. By (2.66), we have

$$\hat{\mathbb{S}}_0^\xi = -\frac{\theta_1(\hbar)\theta_1(2\eta)}{\theta_1'(0)}\Big(\phi(2z,-2\eta)v(z+\hbar,-2\xi|0,\nu_1,\nu_2,\nu_3)e^{2\hbar\partial_z}+ \tag{3.94}$$

$$+\phi(-2z,-2\eta)v(z-\hbar,-2\eta|0,\nu_1,\nu_2,\nu_3)+\big(z\to -z\big)\Big)\,,$$

$$\hat{\mathbb{S}}_k^\xi = \imath^{\delta_{k,2}}\frac{\theta_1(\hbar)\theta_{k+1}(2\hbar)\theta_{k+1}(2\eta)}{\theta_1'(0)\theta_1(2\hbar)}\Big(\varphi_k(2z,-2\eta+\omega_k)v(z+\hbar,-2\xi|0,\nu_1,\nu_2,\nu_3)e^{2\hbar\partial_z}+ \tag{3.95}$$

$$+\varphi_k(-2z,-2\eta+\omega_k)v(z-\hbar,-2\eta|0,\nu_1,\nu_2,\nu_3)+\big(z\to -z\big)\Big)-$$

$$-u_k\frac{\theta_1(2\eta)\theta_1(2\eta+2\hbar)\theta_{j+1}(0)\theta_{i+1}(0)}{\theta_{j+1}(\hbar)\theta_{i+1}(\hbar)}\,,$$

$$\delta_1^{\mathcal{B}C_1} = \delta_1(2\hbar,\eta)\frac{\theta_1^2(\hbar)}{\theta_1^2(2\hbar)}\sum_{k=1}^3\nu_k^2\Big(\wp(2\xi+\omega_k)-\wp(2\eta+\omega_k)\Big)+\frac{\delta_2^2(\hbar,\eta)\theta_1^2(\hbar)}{4\theta_1^2(2\hbar)}\sum_{k=1}^3 u_k^2\frac{\theta_{k+1}^2(\hbar)}{\theta_{k+1}^2(0)}\,, \tag{3.96}$$

$$\delta_2^{\mathcal{B}C_1} = -64\,\delta_2(2\hbar,\eta)\frac{\theta_1^{10}(\hbar)\theta_1^2(4\hbar)}{\theta_1^{12}(2\hbar)}\sum_{k=1}^3\nu_k^2\Big(\wp(2\xi+\omega_k)-\wp(2\eta+\omega_k)\Big)- \tag{3.97}$$

$$-64\,\delta_2(\hbar,\eta)\theta_1(2\eta-2\hbar)\theta_1(2\eta+4\hbar)\frac{\theta_1^{10}(\hbar)\theta_1^2(4\hbar)}{\theta_1^{12}(2\hbar)}\sum_{k=1}^3 u_k^2\frac{\theta_{k+1}^2(\hbar)}{\theta_{k+1}^2(2\hbar)}\,,$$

$$\hat{\mathbb{S}}_0^L = -\frac{\theta_1(\hbar)\theta_1(2\eta)}{\theta_1'(0)}\Big(\phi(2z+2\hbar,-2\eta)v(z,-2\eta|0,\nu_1,\nu_2,\nu_3)e^{2\hbar\partial_z}+ \tag{3.98}$$

$$+\phi(-2z-2\hbar,-2\eta)v(z,-2\eta|0,\nu_1,\nu_2,\nu_3)+\big(z\to -z\big)\Big)\,,$$

$$\hat{\mathbb{S}}_k^L = \imath^{\delta_{k,2}}\frac{\theta_1(\hbar)\theta_{k+1}(2\hbar)\theta_{k+1}(2\eta)}{\theta_1'(0)\theta_1(2\hbar)}\Big(\varphi_k(2z+2\hbar,-2\eta+\omega_k)v(z,-2\eta|0,\nu_1,\nu_2,\nu_3)e^{2\hbar\partial_z}+ \tag{3.99}$$

$$+\varphi_k(-2z-2\hbar,-2\eta+\omega_k)v(z,-2\eta|0,\nu_1,\nu_2,\nu_3)+\big(z\to -z\big)\Big)-$$

$$-u_k\frac{\theta_1(2\eta)\theta_1(2\eta+2\hbar)\theta_{j+1}(0)\theta_{i+1}(0)}{\theta_{j+1}(\hbar)\theta_{i+1}(\hbar)}\,,$$

$$\hat{\mathbb{S}}^R_0=-\frac{\theta_1(\hbar)\theta_1(2\eta)}{\theta_1'(0)}\Big(\phi(2z,-2\eta)v(z+\hbar,-2\eta|0,\nu_1,\nu_2,\nu_3)e^{2\hbar\partial_z}+ \tag{3.100}$$

$$+\phi(-2z,-2\eta)v(z-\hbar,-2\eta|0,\nu_1,\nu_2,\nu_3)+\big(z\to -z\big)\Big)\,,$$

$$\hat{\mathbb{S}}^R_k=\imath^{\delta_{k,2}}\frac{\theta_1(\hbar)\theta_{k+1}(2\hbar)\theta_{k+1}(2\eta)}{\theta_1'(0)\theta_1(2\hbar)}\Big(\varphi_k(2z,-2\eta+\omega_k)v(z+\hbar,-2\eta|0,\nu_1,\nu_2,\nu_3)e^{2\hbar\partial_z}+ \tag{3.101}$$

$$+\varphi_k(-2z,-2\eta+\omega_k)v(z-\hbar,-2\eta|0,\nu_1,\nu_2,\nu_3)+\big(z\to -z\big)\Big)-$$

$$-u_k\frac{\theta_1(2\eta)\theta_1(2\eta+2\hbar)\theta_{j+1}(0)\theta_{i+1}(0)}{\theta_{j+1}(\hbar)\theta_{i+1}(\hbar)}\,,$$

$$\delta_1^{\mathcal{B}C_1}=\frac{\delta_2^2(\hbar,\eta)\theta_1^2(\hbar)}{4\theta_1^2(2\hbar)}\sum_{k=1}^3 u_k^2\frac{\theta_{k+1}^2(\hbar)}{\theta_{k+1}^2(0)}\,, \tag{3.102}$$

$$\delta_2^{\mathcal{B}C_1}=-64\,\delta_2(\hbar,\eta)\theta_1(2\eta-2\hbar)\theta_1(2\eta+4\hbar)\frac{\theta_1^{10}(\hbar)\theta_1^2(4\hbar)}{\theta_1^{12}(2\hbar)}\sum_{k=1}^3 u_k^2\frac{\theta_{k+1}^2(\hbar)}{\theta_{k+1}^2(2\hbar)}\,. \tag{3.103}$$

For all these representations, the structure constants $\widetilde{\nu}_k$ are given by

$$\widetilde{\nu}_k=-4u_k\frac{\theta_1^2(\hbar)\theta_{k+1}(2\eta+2\hbar)\theta_{k+1}(2\eta)}{\theta_{k+1}(2\hbar)\theta_{k+1}(0)}\,. \tag{3.104}$$

# 4 Elliptic modular doubles and the modular triple

In this section, we apply the results of the previous section to the elliptic modular double algebra. We then construct commuting operators, which we use in the next section to study the Ruijsenaars-van Diejen model.

## 4.1 Spiridonov's two elliptic modular doubles

**Definition.** This algebra was introduced by V. Spiridonov in [42] as an elliptic deformation of Faddeev's modular double [11]. The elliptic modular double is defined in terms of two copies $\mathcal{A}$ and $\widetilde{\mathcal{A}}$ of the Sklyanin algebra with defining relations

$$[S_i,S_j]_-=\imath[S_0,S_k]_+,\qquad [S_0,S_k]_-=\imath\frac{J_j-J_i}{J_k}\,[S_i,S_j]_+\,, \tag{4.1}$$

$$[\widetilde{S}_i,\widetilde{S}_j]_-=\imath[\widetilde{S}_0,\widetilde{S}_k]_+,\qquad [\widetilde{S}_0,\widetilde{S}_k]_-=\imath\frac{\widetilde{J}_j-\widetilde{J}_i}{\widetilde{J}_k}\,[\widetilde{S}_i,\widetilde{S}_j]_+\,, \tag{4.2}$$

for every cyclic permutation $(i,j,k)$ of $(1,2,3)$. At the elliptic level, there are two versions of the modular double with different cross-commutation relations. In the first case, the cross-commutation relations between the generators $S_\alpha$ and $\widetilde{S}_\alpha$ of the two Sklyanin algebras have the form

$$\begin{aligned} S_\alpha \widetilde{S}_\beta &= \widetilde{S}_\beta S_\alpha, \quad \alpha,\beta \in \{0,3\} \text{ or } \alpha,\beta \in \{1,2\}\,, \\ S_\alpha \widetilde{S}_\beta &= -\widetilde{S}_\beta S_\alpha, \quad \alpha \in \{0,3\}, \beta \in \{1,2\} \text{ or } \alpha \in \{1,2\}, \beta \in \{0,3\}\,. \end{aligned} \tag{4.3}$$

In the second case, we have

$$\begin{aligned} S_\alpha \widetilde{S}_\beta &= \widetilde{S}_\beta S_\alpha, \quad \alpha,\beta \in \{0,1\} \text{ or } \alpha,\beta \in \{2,3\}\,, \\ S_\alpha \widetilde{S}_\beta &= -\widetilde{S}_\beta S_\alpha, \quad \alpha \in \{0,1\}, \beta \in \{2,3\} \text{ or } \alpha \in \{2,3\}, \beta \in \{0,1\}\,. \end{aligned} \tag{4.4}$$

Equivalently, both cases can be written as

$$S_\alpha \widetilde{S}_\beta = \varepsilon_{1,2}^{\alpha\beta}\, \widetilde{S}_\beta S_\alpha\,, \tag{4.5}$$

where $\varepsilon_{1,2}^{\alpha\beta}$ are the matrix elements of

$$\varepsilon_1 = \begin{pmatrix} 1 & -1 & -1 & 1 \\ -1 & 1 & 1 & -1 \\ -1 & 1 & 1 & -1 \\ 1 & -1 & -1 & 1 \end{pmatrix}, \tag{4.6}$$

in the first case, and

$$\varepsilon_2 = \begin{pmatrix} 1 & 1 & -1 & -1 \\ 1 & 1 & -1 & -1 \\ -1 & -1 & 1 & 1 \\ -1 & -1 & 1 & 1 \end{pmatrix}, \tag{4.7}$$

in the second. Both matrices have the form

$$\varepsilon^{\alpha\beta} = a^\alpha b^\beta\,, \qquad a^\alpha\,, b^\beta \in \{-1,1\}\,. \tag{4.8}$$

**Representation by difference operators.** In the original paper [42], Spiridonov proposed a representation of the elliptic modular doubles by difference operators. Derkachov and Spiridonov [8] constructed finite-dimensional representations of the first elliptic modular double, which corresponds to the interchange $2\hbar \leftrightarrow \tau$.

Consider the standard Sklyanin operators $\hat{S}_\alpha$ in (3.1) and the following dual operators. In the first case, the cross-relations are

$$S_\alpha \overset{1}{\widetilde{S}}_\beta = \varepsilon_1^{\alpha\beta}\, \overset{1}{\widetilde{S}}_\beta S_\alpha\,, \tag{4.9}$$

with the sign matrix (4.6). The dual operators are

$$\begin{gathered} \overset{1}{\hat{\widetilde{S}}}_\alpha := \hat{S}_\alpha(z, \frac{\tau}{2}, \eta|2\hbar) = \\ \imath^{\delta_{\alpha,2}} \frac{\theta_{\alpha+1}(\frac{\tau}{2}|2\hbar)}{\theta_1(2z|2\hbar)} \Big(\theta_{\alpha+1}(2z-2\eta|2\hbar) e^{\frac{\tau}{2}\partial_z} - \theta_{\alpha+1}(-2z-2\eta|2\hbar) e^{-\frac{\tau}{2}\partial_z}\Big) = \\ = \imath^{\delta_{\alpha,2}} \frac{\theta_{\alpha+1}(\frac{\tau}{2}|2\hbar)\theta_{\alpha+1}(-2\eta|2\hbar)}{\theta_1'(0|2\hbar)} \Big(\varphi_\alpha(2z, -2\eta + \widetilde{\omega}_\alpha^1|2\hbar) e^{\frac{\tau}{2}\partial_z} + \varphi_\alpha(-2z, -2\eta + \widetilde{\omega}_\alpha^1|2\hbar) e^{-\frac{\tau}{2}\partial_z}\Big)\,. \end{gathered} \tag{4.10}$$

In the second case, the cross-relations are

$$S_\alpha \widetilde{S}^2_\beta = \varepsilon_2^{\alpha\beta}\, \widetilde{S}^2_\beta S_\alpha\,, \tag{4.11}$$

with the sign matrix (4.7). The dual operators are

$$\hat{\widetilde{S}}^2_\alpha := \hat{S}_\alpha\Big(\frac{z}{2\hbar}, \frac{1}{4\hbar}, \frac{\eta}{2\hbar}\Big|\frac{\tau}{2\hbar}\Big) = \tag{4.12}$$

$$= \imath^{\delta_{\alpha,2}} \frac{\theta_{\alpha+1}\big(\frac{1}{4\hbar}\,\big|\,\frac{\tau}{2\hbar}\big)}{\theta_1\big(\frac{z}{\hbar}\,\big|\,\frac{\tau}{2\hbar}\big)} \left(\theta_{\alpha+1}\Big(\frac{z-\eta}{\hbar}\,\Big|\,\frac{\tau}{2\hbar}\Big)\, e^{\frac{1}{2}\partial_z} - \theta_{\alpha+1}\Big(\frac{-z-\eta}{\hbar}\,\Big|\,\frac{\tau}{2\hbar}\Big)\, e^{-\frac{1}{2}\partial_z}\right) =$$

$$= \imath^{\delta_{\alpha,2}} \frac{\theta_{\alpha+1}\left(\frac{1}{4\hbar}|\frac{\tau}{2\hbar}\right)\theta_{\alpha+1}\left(-\frac{\eta}{\hbar}|\frac{\tau}{2\hbar}\right)}{\theta_1'\left(0|\frac{\tau}{2\hbar}\right)} \left(\varphi_\alpha\left(\frac{z}{\hbar}, -\frac{\eta}{\hbar} + \widetilde{\omega}^2_\alpha\Big|\frac{\tau}{2\hbar}\right) e^{\frac{1}{2}\partial_z} + \varphi_\alpha\left(-\frac{z}{\hbar}, -\frac{\eta}{\hbar} + \widetilde{\omega}^2_\alpha\Big|\frac{\tau}{2\hbar}\right) e^{-\frac{1}{2}\partial_z}\right).$$

The dual half-periods have the form

$$\begin{aligned} &\widetilde{\omega}^1_1 = \frac{1}{2}\,, \qquad \widetilde{\omega}^1_2 = \hbar + \frac{1}{2}\,, \qquad \widetilde{\omega}^1_3 = \hbar\,, \\ &\widetilde{\omega}^2_1 = \frac{1}{2}\,, \qquad \widetilde{\omega}^2_2 = \frac{1}{2} + \frac{\tau}{4\hbar}\,, \qquad \widetilde{\omega}^2_3 = \frac{\tau}{4\hbar}\,. \end{aligned} \tag{4.13}$$

Notice that (4.10) and (4.12) are well-defined for $\mathrm{Im}(\hbar) > 0$ and $\mathrm{Im}\,(\tau/\hbar) > 0$. The transformation of parameters, by which (4.10) and (4.12) were obtained, can be associated with actions from $\mathrm{SL}(3,\mathbb{Z})$, for this one needs to reparameterize $z$, $\hbar$, $\tau$ as follows

$$z = \frac{w}{w_3}\,, \qquad \tau = \frac{w_1}{w_3}\,, \qquad \hbar = \frac{w_2}{2w_3}\,, \qquad \mathbf{w} = (w_1, w_2, w_3)^{\mathbf{T}} \in \mathbb{C}^3\,. \tag{4.14}$$

Operators (3.1) are invariant under change $z \to -z$, in terms of $\mathbf{w}$ this change is equivalent to $\mathbf{w} \to -\mathbf{w}$. Then the interchange $2\hbar \leftrightarrow \tau$ is induced by the change $(w_1, w_2, w_3) \to -(w_2, w_1, w_3)$, which can be represented as the action on $\mathbf{w}$ by

$$T = \begin{pmatrix} 0 & -1 & 0 \\ -1 & 0 & 0 \\ 0 & 0 & -1 \end{pmatrix} \in \mathrm{SL}(3,\mathbb{Z})\,. \tag{4.15}$$

Analogously, the change $(z, \hbar, \tau) \to (\frac{z}{2\hbar}, \frac{1}{4\hbar}, \frac{\tau}{2\hbar})$ is the change $(w_1, w_2, w_3) \to -(w_1, w_3, w_2)$ represented by

$$P = \begin{pmatrix} -1 & 0 & 0 \\ 0 & 0 & -1 \\ 0 & -1 & 0 \end{pmatrix} \in \mathrm{SL}(3,\mathbb{Z})\,. \tag{4.16}$$

## 4.2 Commutative families of quadratic elements

Consider the following elements of the second graded component $\mathcal{A}_2$ of the Sklyanin algebra $\mathcal{A}$:

$$X_\alpha(U) = US_\alpha, \quad Y_\alpha(U) = S_\alpha U, \quad U = \sum_{\alpha=0}^{3} u_\alpha S_\alpha, \qquad \alpha = 0,1,2,3\,. \tag{4.17}$$

The analogous elements for the second copy $\widetilde{\mathcal{A}}$ are

$$\widetilde{X}_\alpha(\widetilde{U}) = \widetilde{U}\widetilde{S}_\alpha, \quad \widetilde{Y}_\alpha(\widetilde{U}) = \widetilde{S}_\alpha\widetilde{U}, \quad \widetilde{U} = \sum_{\alpha=0}^{3} \widetilde{u}_\alpha \widetilde{S}_\alpha, \qquad \alpha = 0,1,2,3\,. \tag{4.18}$$

By Theorem 2.1, the elements $X_\alpha(U)$ and $Y_\alpha(U)$ generate quadratic algebras. In particular cases, these algebras degenerate into the standard Sklyanin algebra and the $BC_1$ Sklyanin algebra.

**Theorem 4.1** *For the elements defined in (4.17) and (4.18) and any matrix of cross-commutation signs of the form (4.8), the following relations hold:*

$$
\begin{aligned}
&[X_\alpha(U), \widetilde{X}_\beta(\widetilde{U})]_- = 0, \quad [Y_\alpha(U), \widetilde{Y}_\beta(\widetilde{U})]_- = 0\,,\\
&[X_\alpha(U), \widetilde{Y}_\beta(\widetilde{U})]_- = 0, \quad [\widetilde{X}_\alpha(\widetilde{U}), Y_\beta(U)]_- = 0\,.
\end{aligned}
\tag{4.19}
$$

**Proof.** Let us consider the first commutator

$$
\begin{aligned}
[X_\alpha(U), \widetilde{X}_\beta(\widetilde{U})]_- &= \sum_{\mu,\nu=0}^{3} u_\mu \widetilde{u}_\nu \big(S_\mu S_\alpha \widetilde{S}_\nu \widetilde{S}_\beta - \widetilde{S}_\nu \widetilde{S}_\beta S_\mu S_\alpha\big) \stackrel{(4.5)}{=}\\
&\stackrel{(4.5)}{=} \sum_{\mu,\nu=0}^{3} u_\mu \widetilde{u}_\nu S_\mu S_\alpha \widetilde{S}_\nu \widetilde{S}_\beta \big(1 - \varepsilon^{\mu\beta}\,\varepsilon^{\mu\nu}\,\varepsilon^{\alpha\beta}\,\varepsilon^{\alpha\nu}\big)\,.
\end{aligned}
\tag{4.20}
$$

It is easy to see that $\varepsilon^{\mu\beta}\,\varepsilon^{\mu\nu}\,\varepsilon^{\alpha\beta}\,\varepsilon^{\alpha\nu} = a^{\mu 2} a^{\alpha 2} b^{\nu 2} b^{\beta 2} = 1$ for all $\alpha, \beta, \mu, \nu \in \{0, 1, 2, 3\}$. Therefore,

$$
[X_\alpha(U), \widetilde{X}_\beta(\widetilde{U})]_- = 0\,. \tag{4.21}
$$

The remaining identities are verified in the same way. ■

### 4.3 Elliptic modular triple

In the original formulation, the generators (3.1), (4.10) and (4.12) share a common spin parameter $\eta$. However, the realizations of the algebras (4.9) and (4.11), with sign matrices (4.6) and (4.7), respectively, can be generalized to distinct spin parameters.

**Proposition 4.1** *For $\hat{S}_\alpha = \hat{S}_\alpha(z, \hbar, \eta|\tau)$ and*

$$
\hat{\tilde{S}}^1_\alpha = \hat{S}_\alpha(z, \frac{\tau}{2}, \eta_1|2\hbar)\,, \qquad \eta_1 = \eta + \frac{k}{4}\,, k \in \mathbb{Z}\,, \tag{4.22}
$$

$$
\hat{\tilde{S}}^2_\alpha = \hat{S}_\alpha(\frac{z}{2\hbar}, \frac{1}{4\hbar}, \frac{\eta_2}{2\hbar}|\frac{\tau}{2\hbar})\,, \tag{4.23}
$$

*the relations (4.11) with the sign matrix (4.7) hold, while (4.9) with the sign matrix (4.6) becomes*

$$
\hat{S}_\alpha \hat{\tilde{S}}^1_\beta = (-1)^k \varepsilon_1^{\alpha\beta}\, \hat{\tilde{S}}^1_\beta \hat{S}_\alpha\,. \tag{4.24}
$$

*Moreover, when $\eta_2 = \frac{N\hbar}{2}$, $N \in \mathbb{Z}$, there are closed cross-commutation relations between $\hat{\tilde{S}}^1_\alpha$ and $\hat{\tilde{S}}^2_\beta$:*

$$
\hat{\tilde{S}}^1_\alpha \hat{\tilde{S}}^2_\beta = (-1)^N \varepsilon^{\alpha\beta}\, \hat{\tilde{S}}^2_\beta \hat{\tilde{S}}^1_\alpha\,. \tag{4.25}
$$

$$
\varepsilon = \begin{pmatrix} 1 & -1 & -1 & 1 \\ 1 & -1 & -1 & 1 \\ -1 & 1 & 1 & -1 \\ -1 & 1 & 1 & -1 \end{pmatrix}\,. \tag{4.26}
$$

**Proof.** The proof uses the quasi-periodicity relations (A.7) and (A.8) for $\phi(z,u)$. Consider first $\hat{\tilde{S}}^1_0$ and $\hat{\tilde{S}}^2_0$. Set

$$d_1 = \frac{\theta_1(\tau/2|2\hbar)\theta_1(-2\eta_1|2\hbar)}{\theta_1'(0|2\hbar)}\,, \qquad d_2 = \frac{\theta_1(1/(4\hbar)|\tau/(2\hbar))\theta_1(-\eta_2/\hbar|\tau/(2\hbar))}{\theta_1'(0|\tau/(2\hbar))}\,. \tag{4.27}$$

Then the representations (3.1), (4.22) and (4.23) give

$$\hat{\tilde{S}}^1_0 = d_1 \sum_{\sigma=\pm1} \phi(2\sigma z, -2\eta_1|2\hbar)e^{\sigma\tau\partial_z/2}\,, \qquad \hat{\tilde{S}}^2_0 = d_2 \sum_{\rho=\pm1} \phi(\rho z/\hbar, -\eta_2/\hbar|\tau/(2\hbar))e^{\rho\partial_z/2}\,. \tag{4.28}$$

Multiplying the operators in this order yields

$$\begin{aligned}\hat{\tilde{S}}^1_0\hat{\tilde{S}}^2_0 &= d_1 d_2 \sum_{\sigma,\rho=\pm1} \phi(2\sigma z, -2\eta_1|2\hbar)\phi\Big(\frac{\rho z}{\hbar} + \frac{\sigma\rho\tau}{2\hbar}, -\frac{\eta_2}{\hbar}\Big|\frac{\tau}{2\hbar}\Big)e^{(\sigma\tau+\rho)\partial_z/2} = \\ &= d_1 d_2 \sum_{\sigma,\rho=\pm1} e^{2\pi\imath\sigma\rho\eta_2/\hbar}\phi(2\sigma z, -2\eta_1|2\hbar)\phi\Big(\frac{\rho z}{\hbar}, -\frac{\eta_2}{\hbar}\Big|\frac{\tau}{2\hbar}\Big)e^{(\sigma\tau+\rho)\partial_z/2}\,.\end{aligned} \tag{4.29}$$

In the opposite order, the first argument of $\phi(2\sigma z, -2\eta_1|2\hbar)$ is shifted by $\sigma\rho$. This is an integer period, so

$$\hat{\tilde{S}}^2_0\hat{\tilde{S}}^1_0 = d_1 d_2 \sum_{\sigma,\rho=\pm1} \phi(2\sigma z, -2\eta_1|2\hbar)\phi\Big(\frac{\rho z}{\hbar}, -\frac{\eta_2}{\hbar}\Big|\frac{\tau}{2\hbar}\Big)e^{(\sigma\tau+\rho)\partial_z/2}\,. \tag{4.30}$$

The two phase factors in (4.29) coincide when

$$e^{2\pi\imath\eta_2/\hbar} = e^{-2\pi\imath\eta_2/\hbar} \quad \Longleftrightarrow \quad \eta_2 = \frac{N\hbar}{2}\,, \qquad N\in\mathbb{Z}\,. \tag{4.31}$$

Their common value is $(-1)^N$, as required by (4.25) for $\alpha=\beta=0$. The remaining relations can be checked in the same manner. Before setting $\eta_2 = N\hbar/2$, we use the definition (A.26) to write the coefficients of $\hat{\tilde{S}}^2_\beta$ in the theta-function form (4.12), which is regular in $\eta_2$ for every $\beta$.

■

Notice that the matrix $\varepsilon$ in (4.26) has the form (4.8), as do $\varepsilon^1$ and $\varepsilon^2$.

The generator $\hat{S}_0$ of the Sklyanin algebra can be viewed as the two-body Hamiltonian of the elliptic Ruijsenaars-Schneider model in the center-of-mass frame [35]. Therefore, the operators $\hat{\tilde{S}}^1_0$ and $\hat{\tilde{S}}^2_0$ (with even $k$ for $\hat{\tilde{S}}^1_0$, arbitrary $\eta_2$ for $\hat{\tilde{S}}^2_0$ and even $N$ for mutual commutativity when $\eta_2 = N\hbar/2$) are commuting Hamiltonians of the same type. A similar statement for the $BC_1$ van Diejen model will be formulated in the next section.

## 4.4 Commuting difference operators

We call the algebra with cross-relations (4.11), (4.24) and (4.25) the *elliptic modular triple algebra.* Applying the results of the previous section, we can construct quadratic combinations of the dual generators (4.22) and (4.23) of the elliptic modular triple algebra. For example, the operators dual to

(3.22) have the following form:

$$\hat{\tilde{Z}}^{1,R}_\alpha = \hat{Z}^R_\alpha\left(z, \frac{\tau}{2}, \eta_1|2\hbar;\widetilde{\nu}^1\right) =$$

$$= \imath^{\delta_{\alpha,2}}\frac{\theta_1(\frac{\tau}{2}|2\hbar)\theta_{\alpha+1}(\tau|2\hbar)\theta_{\alpha+1}(-2\eta_1|2\hbar)}{\theta_1'(0|2\hbar)\theta_1(\tau|2\hbar)}\Big(\varphi_\alpha\left(2z,-2\eta_1+\widetilde{\omega}^1_\alpha|2\hbar\right)\widetilde{v}^1\left(z+\frac{\tau}{2},-2\eta_1\right)e^{\tau\partial_z}+ \tag{4.32}$$

$$+\varphi_\alpha\left(-2z,-2\eta_1+\widetilde{\omega}^1_\alpha|2\hbar\right)\widetilde{v}^1\left(z-\frac{\tau}{2},-2\eta_1\right)+\left(z\to -z\right)\Big),$$

$$\hat{\tilde{Z}}^{2,R}_\alpha = \hat{Z}^R_\alpha\left(\frac{z}{2\hbar}, \frac{1}{4\hbar}, \frac{\eta_2}{2\hbar}|\frac{\tau}{2\hbar};\widetilde{\nu}^2\right) = \tag{4.33}$$

$$= \imath^{\delta_{\alpha,2}}\frac{\theta_1\left(\frac{1}{4\hbar}|\frac{\tau}{2\hbar}\right)\theta_{\alpha+1}\left(\frac{1}{2\hbar}|\frac{\tau}{2\hbar}\right)\theta_{\alpha+1}\left(-\frac{\eta_2}{\hbar}|\frac{\tau}{2\hbar}\right)}{\theta_1'\left(0|\frac{\tau}{2\hbar}\right)\theta_1\left(\frac{1}{2\hbar}|\frac{\tau}{2\hbar}\right)}\Big(\varphi_\alpha\left(\frac{z}{\hbar},-\frac{\eta_2}{\hbar}+\widetilde{\omega}^2_\alpha|\frac{\tau}{2\hbar}\right)\widetilde{v}^2\left(\frac{z}{2\hbar}+\frac{1}{4\hbar},-\frac{\eta_2}{\hbar}\right)e^{\partial_z}+$$

$$+\varphi_\alpha\left(-\frac{z}{\hbar},-\frac{\eta_2}{\hbar}+\widetilde{\omega}^2_\alpha|\frac{\tau}{2\hbar}\right)\widetilde{v}^2\left(\frac{z}{2\hbar}-\frac{1}{4\hbar},-\frac{\eta_2}{\hbar}\right)+\left(z\to -z\right)\Big),$$

where $\widetilde{v}^1(z,\eta) = v\left(z,\eta|2\hbar;\widetilde{\nu}^1\right)$, $\widetilde{v}^2(z,\eta) = v\left(z,\eta|\frac{\tau}{2\hbar};\widetilde{\nu}^2\right)$ and

$$\widetilde{\nu}^1_\alpha = \imath^{\delta_{\alpha,2}}\frac{\theta_{\alpha+1}\left(\frac{\tau}{2}|2\hbar\right)\theta_{\alpha+1}(-2\eta_1|2\hbar)}{\theta_1'(0|2\hbar)}\,\widetilde{u}_\alpha^{\,1}, \qquad \widetilde{\nu}^2_\alpha = \imath^{\delta_{\alpha,2}}\frac{\theta_{\alpha+1}\left(\frac{1}{4\hbar}|\frac{\tau}{2\hbar}\right)\theta_{\alpha+1}\left(-\frac{\eta_2}{\hbar}|\frac{\tau}{2\hbar}\right)}{\theta_1'\left(0|\frac{\tau}{2\hbar}\right)}\,\widetilde{u}^2_\alpha\,. \tag{4.34}$$

The cross-relations of the elliptic modular triple algebra satisfy the hypotheses of Theorem 4.1. Therefore,

$$\left[\hat{Z}^R_\alpha, \hat{\tilde{Z}}^{1,R}_\beta\right]_- = 0\,, \qquad \text{for } \eta_1 = \eta + \frac{k}{4}\,, k\in\mathbb{Z}\,,$$

$$\left[\hat{Z}^R_\alpha, \hat{\tilde{Z}}^{2,R}_\beta\right]_- = 0\,, \qquad \text{for any } \eta,\eta_2, \tag{4.35}$$

$$\left[\hat{\tilde{Z}}^{1,R}_\alpha, \hat{\tilde{Z}}^{2,R}_\beta\right]_- = 0\,, \qquad \text{for any } \eta_2 = \frac{N\hbar}{2}\,, N\in\mathbb{Z}\,.$$

All dual generators are obtained by transforming the arguments $z$, $\hbar$ and $\tau$ of the original generators $\hat{S}_\alpha(z,\hbar,\eta|\tau)$ or $\hat{Z}_a(z,\hbar,\eta|\tau;\nu)$. The same transformations can be applied to construct operators dual to (3.30):

$$\hat{Z}^\xi_\alpha = \imath^{\delta_{\alpha,2}}\frac{\theta_1(\hbar)\theta_{\alpha+1}(2\hbar)\theta_{\alpha+1}(-2\eta)}{\theta_1'(0)\theta_1(2\hbar)}\Big(\varphi_\alpha(2z,-2\eta+\omega_\alpha)v(z+\hbar,-2\xi)e^{2\hbar\partial_z}+$$
$$+\varphi_\alpha(-2z,-2\eta+\omega_\alpha)v(z-\hbar,-2\eta)+\left(z\to -z\right)\Big), \qquad \alpha\in\{0,1,2,3\}\,. \tag{4.36}$$

**Proposition 4.2** *The dual operators*

$$\hat{\tilde{Z}}^{1,\xi_1}_\alpha = \hat{Z}^{\xi_1}_\alpha\left(z, \frac{\tau}{2}, \eta_1|2\hbar;\widetilde{\nu}^1\right) =$$

$$= \imath^{\delta_{\alpha,2}}\frac{\theta_1(\frac{\tau}{2}|2\hbar)\theta_{\alpha+1}(\tau|2\hbar)\theta_{\alpha+1}(-2\eta_1|2\hbar)}{\theta_1'(0|2\hbar)\theta_1(\tau|2\hbar)}\Big(\varphi_\alpha\left(2z,-2\eta_1+\widetilde{\omega}^1_\alpha|2\hbar\right)\widetilde{v}^1\left(z+\frac{\tau}{2},-2\xi_1\right)e^{\tau\partial_z}+ \tag{4.37}$$

$$+\varphi_\alpha\left(-2z,-2\eta_1+\widetilde{\omega}^1_\alpha|2\hbar\right)\widetilde{v}^1\left(z-\frac{\tau}{2},-2\eta_1\right)+\left(z\to -z\right)\Big),$$

$$\hat{\tilde{Z}}_\alpha^{2,\xi_2/2\hbar} = \hat{Z}_\alpha^{\xi_2/2\hbar}\left(\frac{z}{2\hbar},\frac{1}{4\hbar},\frac{\eta_2}{2\hbar}|\frac{\tau}{2\hbar};\tilde{\nu}^2\right) = \tag{4.38}$$

$$= \imath^{\delta_{\alpha,2}}\frac{\theta_1\left(\frac{1}{4\hbar}|\frac{\tau}{2\hbar}\right)\theta_{\alpha+1}\left(\frac{1}{2\hbar}|\frac{\tau}{2\hbar}\right)\theta_{\alpha+1}\left(-\frac{\eta_2}{\hbar}|\frac{\tau}{2\hbar}\right)}{\theta_1'\left(0|\frac{\tau}{2\hbar}\right)\theta_1\left(\frac{1}{2\hbar}|\frac{\tau}{2\hbar}\right)}\Bigg(\varphi_\alpha\left(\frac{z}{\hbar},-\frac{\eta_2}{\hbar}+\tilde{\omega}_\alpha^2|\frac{\tau}{2\hbar}\right)\tilde{v}^2\left(\frac{z}{2\hbar}+\frac{1}{4\hbar},-\frac{\xi_2}{\hbar}\right)e^{\partial_z}+$$

$$+\varphi_\alpha\left(-\frac{z}{\hbar},-\frac{\eta_2}{\hbar}+\tilde{\omega}_\alpha^2|\frac{\tau}{2\hbar}\right)\tilde{v}^2\left(\frac{z}{2\hbar}-\frac{1}{4\hbar},-\frac{\eta_2}{\hbar}\right)+\left(z\to -z\right)\Bigg),$$

*together with $\hat{Z}_\alpha^\xi$ in (4.36), form a commutative family in the following sense:*

$$\left[\hat{Z}_\alpha^\xi,\hat{\tilde{Z}}_\beta^{1,\xi_1}\right]_- = 0\,,\qquad \textit{for } \eta_1+\xi_1=\eta+\xi+\frac{k}{4}\,,k\in\mathbb{Z}\,,$$

$$\left[\hat{Z}_\alpha^\xi,\hat{\tilde{Z}}_\beta^{2,\xi_2/2\hbar}\right]_- = 0\,,\qquad \textit{for any } \eta,\xi,\eta_2,\xi_2\,, \tag{4.39}$$

$$\left[\hat{\tilde{Z}}_\alpha^{1,\xi_1},\hat{\tilde{Z}}_\beta^{2,\xi_2/2\hbar}\right]_- = 0\,,\qquad \textit{for any } \eta_1,\xi_1 \textit{ and } \eta_2+\xi_2=\frac{N\hbar}{2}\,,N\in\mathbb{Z}\,.$$

**Proof.** The proof uses the quasi-periodicity relations (A.7), (A.8) from Appendix A, together with the definitions (A.25), (A.32) and, in this sense, is completely similar to the proof of Proposition 4.1.

■

# 5 Applications to the $BC_1$ Ruijsenaars-van Diejen model

In this section, we consider the quantum elliptic $BC_1$ Ruijsenaars-van Diejen model [9]. Its Hamiltonian can be written in several forms (see [17]). We use the form proposed in [4]:

$$\hat{H}^{\rm 8vD} := \hat{H}^{\rm 8vD}(z|\tau;\hbar,\nu,\bar{\nu},\eta,\xi) =$$

$$= \bar{v}(z,-2\eta)v(z+\hbar,-2\xi)e^{2\hbar\partial_z}+\bar{v}(-z,-2\eta)v(-z+\hbar,-2\xi)e^{-2\hbar\partial_z}+ \tag{5.1}$$

$$+\bar{v}(-z,-2\eta)v(z-\hbar,-2\eta)+\bar{v}(z,-2\eta)v(-z-\hbar,-2\eta)\,.$$

The function $v(z,\eta)$ is defined in (A.32), and $\bar{v}(z,\eta)=v(z,\eta|\tau;\bar{\nu})$ is the same function with the parameters $\bar{\nu}_a$ in place of $\nu_a$. The Hamiltonian (5.1) is the one-body case ($n=1$) of the elliptic $n$-body $BC_n$ Ruijsenaars-van Diejen model with nine coupling constants. Its quantum integrability was proved in part in [9] and in full by Y. Komori and K. Hikami in [17] using Dunkl-Cherednik-type operators. In the one-body case, the expression (5.1) contains ten parameters: $\nu_0,...,\nu_3$, $\bar{\nu}_0,...,\bar{\nu}_3$, $\eta$ and $\xi$ (the parameters $\hbar$ and $\tau$ are fixed). The Hamiltonian is linear in each of the two sets of four parameters, $\nu_a$ and $\bar{\nu}_a$. Thus there are $10-2=8$ independent parameters, which explains the notation $\hat{H}^{\rm 8vD}$.

In [32], Rains and Ruijsenaars related van Diejen operators to quadratic combinations of Sklyanin generators (1.10). We also propose a connection between the van Diejen operator (5.1) and realizations of the algebras mentioned above.

## 5.1 The one-body Ruijsenaars-van Diejen Hamiltonian in terms of quadratic algebras

**The one-body Ruijsenaars-van Diejen Hamiltonian in terms of $\mathcal{A}(\mathcal{A}_2)$.**

**Proposition 5.1** *A linear combination of the generators (3.30) yields the $BC_1$ van Diejen operator (5.1):*

$$\sum_{\alpha=0}^{3}\bar{u}_\alpha \hat{Z}^{\xi}_\alpha=\sum_{\alpha=0}^{3}\imath^{\delta_{\alpha,2}}\bar{u}_\alpha\frac{\theta_1(\hbar)\theta_{\alpha+1}(2\hbar)\theta_{\alpha+1}(-2\eta)}{\theta_1'(0)\theta_1(2\hbar)}\Big(\varphi_\alpha(2z,-2\eta+\omega_\alpha)v(z+\hbar,-2\xi)e^{2\hbar\partial_z}+$$
$$+\varphi_\alpha(-2z,-2\eta+\omega_\alpha)v(z-\hbar,-2\eta)+\big(z\to -z\big)\Big)=$$
$$=\bar{v}(z,-2\eta)v(z+\hbar,-2\xi)e^{2\hbar\partial_z}+\bar{v}(-z,-2\eta)v(z-\hbar,-2\eta)+\big(z\to -z\big)=\hat{H}^{8vD}\,, \tag{5.2}$$

*with*

$$\nu_\alpha=\imath^{\delta_{\alpha,2}}\frac{\theta_{\alpha+1}(\hbar)\theta_{\alpha+1}(-2\eta)}{\theta_1'(0)}\,u_\alpha\,,\qquad \bar{\nu}_\alpha=\imath^{\delta_{\alpha,2}}\frac{\theta_1(\hbar)\theta_{\alpha+1}(2\hbar)\theta_{\alpha+1}(-2\eta)}{\theta_1'(0)\theta_1(2\hbar)}\,\bar{u}_\alpha\,. \tag{5.3}$$

The van Diejen Hamiltonian with seven independent parameters

$$\hat{H}^{7\mathrm{vD}}:=\hat{H}^{7\mathrm{vD}}(z|\tau;\hbar,\nu,\bar{\nu},\eta)=\hat{H}^{8\mathrm{vD}}(z|\tau;\hbar,\nu,\bar{\nu},\eta,\eta)\,, \tag{5.4}$$

can be expressed in terms of the generators (3.22) as a product of two linear combinations of the Sklyanin generators (3.1).

**Proposition 5.2** *A linear combination of the generators (3.22) yields*

$$\hat{H}^{7vD}=\sum_{\alpha=0}^{3}\bar{u}_\alpha\hat{Z}^R_\alpha\,, \tag{5.5}$$

*with*

$$\nu_\alpha=\imath^{\delta_{\alpha,2}}\frac{\theta_{\alpha+1}(\hbar)\theta_{\alpha+1}(-2\eta)}{\theta_1'(0)}\,u_\alpha\,,\qquad \bar{\nu}_\alpha=\imath^{\delta_{\alpha,2}}\frac{\theta_1(\hbar)\theta_{\alpha+1}(2\hbar)\theta_{\alpha+1}(-2\eta)}{\theta_1'(0)\theta_1(2\hbar)}\,\bar{u}_\alpha\,. \tag{5.6}$$

*Since* $\hat{Z}^R_\alpha=\frac{\theta_1(\hbar)\theta_{\alpha+1}(2\hbar)}{\theta_1(2\hbar)\theta_{\alpha+1}(\hbar)}\,\hat{S}_\alpha\sum_{\beta=0}^{3}u_\beta\hat{S}_\beta$, *(5.5) can be rewritten as follows:*

$$\hat{H}^{7vD}=\sum_{\alpha,\beta=0}^{3}\frac{\theta_1(\hbar)\theta_{\alpha+1}(2\hbar)}{\theta_1(2\hbar)\theta_{\alpha+1}(\hbar)}\,\bar{u}_\alpha u_\beta\hat{S}_\alpha\hat{S}_\beta\,. \tag{5.7}$$

**The one-body Ruijsenaars-van Diejen Hamiltonian in terms of $\mathcal{B}C_1$.** In the above propositions, we represent the van Diejen Hamiltonian as a linear combination of generators of the algebra (2.50) or as a product of two linear combinations of generators of the four-dimensional Sklyanin algebra. Furthermore, the representations (3.82)–(3.83) of the $BC_1$ Sklyanin algebra can be related to $\hat{H}^{7vD}$.

**Proposition 5.3** *The van Diejen Hamiltonian (5.4) with seven independent parameters can be written as the following quadratic expression in the generators (3.82):*

$$\hat{H}^{7vD}(z|\tau;2\hbar,\nu,\bar{\nu},\xi-\hbar)=\hat{\mathbb{S}}^{\xi}_0(z|\bar{u})\hat{\mathbb{S}}^{\xi}_0(z|u)\,, \tag{5.8}$$

*with*

$$\nu_\alpha=\imath^{\delta_{\alpha,2}}\frac{\theta_1(\hbar)\theta_{\alpha+1}(\hbar)\theta_{\alpha+1}(2\hbar)\theta_{\alpha+1}(2\hbar-2\xi)}{\theta_1'(0)\theta_{\alpha+1}(-2\xi)}\,u_\alpha\,,$$
$$\bar{\nu}_\alpha=\imath^{\delta_{\alpha,2}}\frac{\theta_1(\hbar)\theta_{\alpha+1}(\hbar)\theta_{\alpha+1}(2\hbar)\theta_{\alpha+1}(2\hbar-2\xi)}{\theta_1'(0)\theta_{\alpha+1}(-2\xi)}\,\bar{u}_\alpha\,. \tag{5.9}$$

**Proof.** Note that $\hat{\mathbb{S}}_0^\xi$ in (3.82) is a linear combination of the standard Sklyanin generators (3.1):

$$\hat{\mathbb{S}}_0^\xi(z|u) = \frac{\theta_1(\hbar)\theta_1(2\hbar)}{\theta_1'(0)}\Big(\phi(2z,2\hbar)v(z+\hbar,-2\xi)e^{2\hbar\partial_z} + \phi(-2z,2\hbar)v(-z+\hbar,-2\xi)e^{-2\hbar\partial_z}\Big) =$$
$$= \sum_{\alpha=0}^{3} u_\alpha \frac{\theta_1(\hbar)\theta_{\alpha+1}(\hbar)}{\theta_{\alpha+1}(-2\xi)}\,\hat{S}_\alpha(z,2\hbar,\xi-\hbar|\tau) = \sum_{\alpha=0}^{3} u'_\alpha \hat{S}_\alpha(z,2\hbar,\xi-\hbar|\tau)\,, \tag{5.10}$$

where

$$u'_\alpha = u_\alpha \frac{\theta_1(\hbar)\theta_{\alpha+1}(\hbar)}{\theta_{\alpha+1}(-2\xi)}\,. \tag{5.11}$$

Then, with the parameters (5.9),

$$\hat{\mathbb{S}}_0^\xi(z|\bar{u})\hat{\mathbb{S}}_0^\xi(z|u) = \sum_{\alpha,\beta=0}^{3} \bar{u}'_\alpha u'_\beta \hat{S}_\alpha(z,2\hbar,\xi-\hbar|\tau)\hat{S}_\beta(z,2\hbar,\xi-\hbar|\tau) \stackrel{(5.4)}{=} \hat{H}^{7\text{vD}}(z|\tau;2\hbar,\nu,\bar{\nu},\xi-\hbar)\,. \tag{5.12}$$

This completes the proof. ∎

## 5.2 Commuting Hamiltonians from the modular triple

The relations (4.39) and (5.2) allow us to construct additional Hamiltonians of van Diejen type

$$\hat{H}_1^{8\text{vD}} := \hat{H}^{8\text{vD}}\left(z|2\hbar;\frac{\tau}{2},\widetilde{\nu}^1,\widetilde{\bar{\nu}}^1,\eta_1,\xi_1\right)\,, \tag{5.13}$$

$$\hat{H}_2^{8\text{vD}} := \hat{H}^{8\text{vD}}\left(\frac{z}{2\hbar}|\frac{\tau}{2\hbar};\frac{1}{4\hbar},\widetilde{\nu}^2,\widetilde{\bar{\nu}}^{\,2},\frac{\eta_2}{2\hbar},\frac{\xi_2}{2\hbar}\right)\,, \tag{5.14}$$

which commute with (5.1) under the following conditions on the parameters:

$$\begin{aligned}
&\left[\hat{H}^{8\text{vD}},\hat{H}_1^{8\text{vD}}\right]_- = 0\,, && \text{for } \eta_1+\xi_1 = \eta+\xi+\frac{k}{4}\,,k\in\mathbb{Z}\,,\\
&\left[\hat{H}^{8\text{vD}},\hat{H}_2^{8\text{vD}}\right]_- = 0\,, && \text{for any } \eta,\xi,\eta_2,\xi_2\,,\\
&\left[\hat{H}_1^{8\text{vD}},\hat{H}_2^{8\text{vD}}\right]_- = 0\,, && \text{for any } \eta_1,\xi_1 \text{ and } \eta_2+\xi_2 = \frac{N\hbar}{2}\,,N\in\mathbb{Z}\,.
\end{aligned} \tag{5.15}$$

By definition (1.7), each Hamiltonian is bilinear in its two sets of coupling constants, and commutativity holds term by term. Indeed, writing the two Hamiltonians as $\hat{H}^{8\text{vD}} = \sum_{a,b=0}^{3}\bar{\nu}_a\nu_b\,D_{ab}$ and $\hat{H}_1^{8\text{vD}} = \sum_{c,d=0}^{3}\widetilde{\bar{\nu}}_c^{\,1}\widetilde{\nu}_d^1\,\widetilde{D}_{cd}^1$, we obtain from the first line of (5.15) that $[D_{ab},\widetilde{D}_{cd}^1] = 0$ for all $a,b,c,d = 0,1,2,3$.

Rains and Ruijsenaars [32] introduced (5.13) as an additional Hamiltonian commuting with (5.1). They did not consider the second version (4.11) of the elliptic modular double algebra or the operator (5.14). The first commutation relation in (5.15) is the modular-pair relation of Rains and Ruijsenaars. To the best of our knowledge, the simultaneous realization of the three van Diejen operators with the parameter conditions in (5.15) has not been described previously.

### 5.3 Classical limit

In our recent work [24], we established the gauge equivalence between the four-parameter van Diejen model and the Zhukovsky-Volterra gyrostat. In our notation, the four-parameter van Diejen model has the following Hamiltonian:

$$H^{\rm 4vD} = v(\eta, z|\tau;\nu)e^{p/2c} + v(\eta, -z|\tau;\nu)e^{-p/2c}\,, \tag{5.16}$$

where $p, z \in \mathbb{C}$ are the momentum and coordinate of the particle, respectively, with canonical Poisson bracket

$$\{p, z\} = 1\,. \tag{5.17}$$

The Hamiltonian (5.16) can be considered as classical limit of $\Delta$-Rains operator

$$\Delta^{\rm R}_{\nu} = v(z, -2\eta)e^{\hbar\partial_z} + v(-z, -2\eta)e^{-\hbar\partial_z}\,. \tag{5.18}$$

To recover (5.16), on should set $\eta_{\rm class} = -2\eta_{\rm quant}$ and $\nu_{\rm quant} = I\nu_{\rm class}$, and replace $e^{\pm\hbar\partial_z}$ by $e^{\pm p/(2c)}$. The identity (A.33) then gives (5.16).

The Zhukovsky-Volterra gyrostat can be described by the following quadratic Poisson structure[5]

$$c\{s_\alpha, s_\beta\} = -\imath\varepsilon_{\alpha\beta\gamma}s_0 s_\gamma\,, \tag{5.19}$$

$$c\{s_0, s_\alpha\} = -\imath\varepsilon_{\alpha\beta\gamma}s_\beta s_\gamma(\wp_\beta - \wp_\gamma) + \imath\varepsilon_{\alpha\beta\gamma}(s_\beta\lambda_\gamma - \lambda_\beta s_\gamma)\,, \tag{5.20}$$

where $s_1, s_2, s_3$ are the components of the angular momentum of a rigid body in three-dimensional space, and $\lambda_1, \lambda_2, \lambda_3$ are the components of a constant vector. The variable $s_0$ plays the role of the Hamiltonian for this system. The Poisson structure (5.19)–(5.20) has two Casimir functions

$$\begin{gathered}\mathbf{C}_1 = s_1^2 + s_2^2 + s_3^2\,,\\ \mathbf{C}_2 = s_0^2 + \sum_{k=1}^{3}\Big(s_k^2\wp_k + 2s_k\lambda_k\Big)\,, \quad \wp_k = \wp(\omega_k)\,.\end{gathered} \tag{5.21}$$

The Poisson map from (5.17) to (5.19), (5.20), with $\breve{\nu} = I\nu$ defined in (A.33), is given by the following change of variables

$$s_0(p, z) = \frac{1}{2}\Big(v(z, \eta|\tau;\breve{\nu})e^{p/2c} + v(-z, \eta|\tau;\breve{\nu})e^{-p/2c}\Big)\,, \tag{5.22}$$

and

$$\begin{gathered}s_\alpha(p, z|\eta, \breve{\nu}_0, \breve{\nu}_1, \breve{\nu}_2, \breve{\nu}_3) = s_\alpha(p, z) = \frac{c_\alpha}{2}\Big(v(z, \eta|\tau;\breve{\nu})e^{p/2c} - v(-z, \eta|\tau;\breve{\nu})e^{-p/2c}\Big)\varphi_\alpha(2z) + \\ + c_\alpha\breve{\nu}_0\varphi_\beta(2z)\varphi_\gamma(2z) + c_\alpha\varphi_\alpha(2z)\sum_{k=1}^{3}\breve{\nu}_k\varphi_k(2z)\end{gathered} \tag{5.23}$$

for distinct $\alpha, \beta, \gamma \in \{1, 2, 3\}$. The identification of parameters is given by

$$\lambda_\alpha = \frac{\breve{\nu}_\alpha}{c_\alpha}\,, \quad \alpha = 1, 2, 3\,, \tag{5.24}$$

with the constants $c_\alpha$ defined in (A.31). Under this change of variables, the Casimir functions (5.21) take the values

$$\mathbf{C}_1 = \breve{\nu}_0^2\,, \quad \mathbf{C}_2 = \sum_{a=0}^{3}\breve{\nu}_a^2\wp(\eta + \omega_a) - \sum_{a=1}^{3}\breve{\nu}_a^2\wp(\omega_a)\,. \tag{5.25}$$

[5] Summation over repeated indices is not implied in (5.19)–(5.20).

The Poisson structure (5.19)–(5.20) is a classical counterpart of the quantum $BC_1$ Sklyanin algebra (2.11)–(2.13). Our goal in this section is to find the classical limit of the representations of the quantum $BC_1$ Sklyanin algebra and compare it with (5.22) and (5.23). Let us consider the representation (3.82)–(3.83)

$$\hat{\mathbb{S}}_0^\xi = \frac{\theta_1(\hbar)\theta_1(2\hbar)}{\theta_1'(0)} \Big(\phi(2z, 2\hbar)v(z+\hbar, -2\xi)e^{2\hbar\partial_z} + \phi(-2z, 2\hbar)v(-z+\hbar, -2\xi)e^{-2\hbar\partial_z}\Big)\,,$$

$$\hat{\mathbb{S}}_k^\xi = \imath^{\delta_{k,2}} \frac{\theta_1(\hbar)\theta_{k+1}^2(2\hbar)}{\theta_1'(0)\theta_1(2\hbar)} \Big(\varphi_k(2z, 2\hbar+\omega_k)v(z+\hbar, -2\xi)e^{2\hbar\partial_z}+ \tag{5.26}$$

$$+\varphi_k(-2z, 2\hbar+\omega_k)v(z-\hbar, 2\hbar) + \big(z\to -z\big)\Big)\,.$$

The generators of the $BC_1$ Sklyanin algebra satisfy

$$\imath[\hat{\mathbb{S}}_0, \hat{\mathbb{S}}_k]_+ = [\hat{\mathbb{S}}_i, \hat{\mathbb{S}}_j]_-\,, \tag{5.27}$$

$$[\hat{\mathbb{S}}_k, \hat{\mathbb{S}}_0]_- = \imath\frac{K_j - K_i}{K_k}\,[\hat{\mathbb{S}}_i, \hat{\mathbb{S}}_j]_+ - 2\imath\frac{1}{K_k}\,(\widetilde{\nu}_i\rho_i\hat{\mathbb{S}}_j - \widetilde{\nu}_j\rho_j\hat{\mathbb{S}}_i)\,. \tag{5.28}$$

The structure constants are

$$K_k = J_i^2(\hbar)J_j^2(\hbar)\frac{\theta_{k+1}^4(0)}{\theta_{k+1}^4(\hbar)}\,, \qquad \frac{K_j - K_i}{K_k} = -J_{ij}(2\hbar)\,,$$

$$\rho_i = J_i(\hbar)J_j(\hbar)J_k(\hbar)\frac{\theta_{j+1}(0)\theta_{k+1}(0)}{\theta_{j+1}(\hbar)\theta_{k+1}(\hbar)}\,, \qquad \widetilde{\nu}_k = -4\theta_1^2(\hbar)u_k\,. \tag{5.29}$$

The Casimir elements are

$$\hat{\Omega}_1^{\mathcal{B}C_1} = \hat{\mathbb{S}}_0^2 + \sum_{k=1}^3 \hat{\mathbb{S}}_k^2 = \delta_1^{\mathcal{B}C_1} = 4u_0^2\theta_1^4(\hbar)\,,$$

$$\hat{\Omega}_2^{\mathcal{B}C_1} = \sum_{k=1}^3 \Big(K_k\big(K_k - K_i - K_j\big)\hat{\mathbb{S}}_k^2 + 2\widetilde{\nu}_k\rho_k K_k\hat{\mathbb{S}}_k\Big) = \delta_2^{\mathcal{B}C_1} = \tag{5.30}$$

$$= 256\,\frac{\theta_1^{10}(\hbar)\theta_1^2(4\hbar)}{\theta_1'(0)^2\theta_1^{10}(2\hbar)}\,\sum_{\alpha=0}^3 \imath^{2\delta_{\alpha,2}}u_\alpha^2\theta_{\alpha+1}^2(\hbar)\theta_{\alpha+1}^2(2\hbar)\Big(\wp(2\xi+\omega_\alpha) - \wp(2\hbar+\omega_\alpha)\Big)\,.$$

**Proposition 5.4** *In the classical limit $\hbar\to 0$, $\hbar\partial_z \to p/4c$, $[\hat{A}, \hat{B}] \to 4\hbar c\{A, B\}$, the quantum algebra (5.27)–(5.29) becomes the Poisson algebra (5.19)–(5.20). If the parameters are chosen as*

$$u_0 = \frac{\theta_1'(0)\check{\nu}_0}{\theta_1(\hbar)\theta_1(2\hbar)}\,, \qquad u_k = -\frac{\check{\nu}_k}{\theta_1'(0)c_k}\,, \qquad \eta = -2\xi\,, \tag{5.31}$$

*the following limits hold:*

$$-\frac{\hat{\mathbb{S}}_0^\xi}{2\theta_1(\hbar)} \to s_0(p + 2\imath\pi c, z)\,, \qquad \frac{\hat{\mathbb{S}}_k^\xi}{\theta_1'(0)} \to s_k(p + 2\imath\pi c, z)\,, \tag{5.32}$$

$$\frac{\hat{\Omega}_1^{\mathcal{B}C_1}}{\theta_1'(0)^2} \to \mathbf{C}_1\,, \qquad \frac{\hat{\Omega}_1^{\mathcal{B}C_1} + \hat{\Omega}_2^{\mathcal{B}C_1}}{\theta_1^2(2\hbar)} + \frac{\theta_1'''(0)}{3\theta_1'(0)}\,\frac{\hat{\Omega}_1^{\mathcal{B}C_1}}{\theta_1'(0)^2} \to \mathbf{C}_2\,. \tag{5.33}$$

**Proof.** The parameters in this proof are those of (5.31). To distinguish the function occurring in the quantum operators from its limit, write $v_\hbar(z,u)$ for $v(z,u)$ with the coefficients (3.92). Substitution of (5.31) gives

$$v_\hbar(z,u) = \breve{\nu}_0\phi(2z,u) - \sum_{k=1}^{3} \frac{\imath^{\delta_{k,2}}\breve{\nu}_k\theta_{k+1}(\hbar)\theta_{k+1}(2\hbar)}{\theta_1'(0)^2 c_k}\,\varphi_k(2z,u+\omega_k)\,. \tag{5.34}$$

By (A.31), $-\imath^{\delta_{k,2}}\theta_{k+1}(0)^2/(\theta_1'(0)^2c_k) = 1$. Hence

$$v_\hbar(z,u) \to v(z,u|\tau;\breve{\nu})\,, \qquad \hbar \to 0\,. \tag{5.35}$$

The limiting coefficients are $\breve{\nu}$, as in (5.22), (5.23).

Now $\theta_1(2\hbar)\phi(\pm 2z, 2\hbar) \to \theta_1'(0)$. Replacing $e^{\pm 2\hbar\partial_z}$ by $e^{\pm p/(2c)}$ in the first line of (5.26) therefore gives

$$\begin{aligned} -\frac{\hat{\mathbb{S}}_0^\xi}{2\theta_1(\hbar)} \to -\frac{1}{2}\Big(v(z,\eta|\tau;\breve{\nu})e^{p/(2c)} + v(-z,\eta|\tau;\breve{\nu})e^{-p/(2c)}\Big) = \\ = s_0(p+2\pi\imath c, z)\,. \end{aligned} \tag{5.36}$$

For $k=1,2,3$, the prefactor in the second line of (5.26), after division by $\theta_1'(0)$, tends to $-c_k/2$. Since $\varphi_k(-2z) = -\varphi_k(2z)$, its terms with shifts tend to

$$-\frac{c_k}{2}\Big(v(z,\eta|\tau;\breve{\nu})e^{p/(2c)} - v(-z,\eta|\tau;\breve{\nu})e^{-p/(2c)}\Big)\varphi_k(2z)\,. \tag{5.37}$$

In the terms without shifts, the contributions proportional to $\breve{\nu}_0$ must be combined before taking the limit. The theta-function expressions (A.3) and (A.26) imply the exact identity

$$\begin{aligned} \varphi_k(-2z, 2\hbar+\omega_k)\phi(2z-2\hbar, 2\hbar) + \varphi_k(2z, 2\hbar+\omega_k)\phi(-2z-2\hbar, 2\hbar) = \\ = \frac{\theta_1'(0)\theta_{k+1}(0)}{\theta_{k+1}(2\hbar)\theta_1(2\hbar)}\Big(\varphi_k(2z+2\hbar) - \varphi_k(2z-2\hbar)\Big)\,. \end{aligned} \tag{5.38}$$

The limit of the right-hand side is $2\varphi_k'(2z)$. The identity obtained from (A.28),

$$\varphi_k'(x) = -\varphi_i(x)\varphi_j(x)\,, \qquad \{i,j,k\} = \{1,2,3\}\,, \tag{5.39}$$

therefore gives the contribution $c_k\breve{\nu}_0\varphi_i(2z)\varphi_j(2z)$. The remaining terms without shifts have a regular limit, equal to $c_k\varphi_k(2z)\sum_{a=1}^{3}\breve{\nu}_a\varphi_a(2z)$. Together with (5.37), this proves

$$\begin{aligned} \frac{\hat{\mathbb{S}}_k^\xi}{\theta_1'(0)} \to -\frac{c_k}{2}\Big(v(z,\eta|\tau;\breve{\nu})e^{p/(2c)} - v(-z,\eta|\tau;\breve{\nu})e^{-p/(2c)}\Big)\varphi_k(2z) + \\ + c_k\breve{\nu}_0\varphi_i(2z)\varphi_j(2z) + c_k\varphi_k(2z)\sum_{a=1}^{3}\breve{\nu}_a\varphi_a(2z) = s_k(p+2\pi\imath c, z)\,. \end{aligned} \tag{5.40}$$

This is (5.32).

To pass from the commutation relations to the Poisson brackets, set

$$\gamma_0 = \frac{\theta_1'''(0)}{3\theta_1'(0)}\,. \tag{5.41}$$

Equations (3.2), (5.29), (5.31) and (A.5) give the expansions

$$J_k(\hbar) = 1 + \hbar^2(\gamma_0 - \wp_k) + O(\hbar^4)\,, \tag{5.42}$$

$$\rho_i = 1 + O(\hbar^2)\,, \qquad K_k = 1 + 2\hbar^2(\gamma_0 + 2\wp_k) + O(\hbar^4)\,, \tag{5.43}$$

$$\frac{K_j - K_i}{K_k} = -J_{ij}(2\hbar) = 4\hbar^2(\wp_j - \wp_i) + O(\hbar^4)\,, \tag{5.44}$$

$$\widetilde{\nu}_k = -4\theta_1^2(\hbar)u_k = \frac{4\theta_1^2(\hbar)\breve{\nu}_k}{\theta_1'(0)c_k} = 4\hbar^2\theta_1'(0)\lambda_k + O(\hbar^4)\,. \tag{5.45}$$

In particular, the $O(\hbar^2)$ term in $K_k$ is retained in (5.43). Substitution into (5.28) gives, for every cyclic triple $(i, j, k)$,

$$[\hat{\mathbb{S}}_k, \hat{\mathbb{S}}_0]_- = 4\hbar^2\Big(\imath(\wp_j - \wp_i)[\hat{\mathbb{S}}_i, \hat{\mathbb{S}}_j]_+ - -2\imath\theta_1'(0)(\lambda_i\hat{\mathbb{S}}_j - \lambda_j\hat{\mathbb{S}}_i)\Big) + O(\hbar^4)\,. \tag{5.46}$$

If we suppress the common arguments $(p + 2\pi\imath c, z)$ of $s_\alpha$, the leading term of the left-hand side is $-8\hbar^2 c\,\theta_1'(0)^2\{s_k, s_0\}$, while that of the right-hand side is

$$8\hbar^2\theta_1'(0)^2\Big(\imath(\wp_j - \wp_i)s_is_j - \imath(\lambda_is_j - \lambda_js_i)\Big)\,. \tag{5.47}$$

Consequently,

$$c\{s_0, s_k\} = \imath(\wp_j - \wp_i)s_is_j - \imath(\lambda_is_j - \lambda_js_i)\,, \tag{5.48}$$

which is (5.20). Similarly, (5.27) gives

$$c\{s_i, s_j\} = -\imath s_0s_k\,, \tag{5.49}$$

and hence (5.19).

Finally, the first relation in (5.33) follows directly from the first line of (5.30) and (5.32). For the second relation, use the exact identity (3.42) and its expansion

$$K_k(K_k - K_i - K_j) = -64\,\frac{\theta_1^8(\hbar)\theta_1^2(4\hbar)}{\theta_1^{10}(2\hbar)}\,J_k(2\hbar) = -1 - 4\hbar^2(\gamma_0 - \wp_k) + O(\hbar^4)\,. \tag{5.50}$$

The terms of order one in $\hat{\Omega}_1^{\mathcal{B}C_1} + \hat{\Omega}_2^{\mathcal{B}C_1}$ cancel. Equations (5.45), (5.50) then give

$$\frac{\hat{\Omega}_1^{\mathcal{B}C_1} + \hat{\Omega}_2^{\mathcal{B}C_1}}{\theta_1^2(2\hbar)} \to s_0^2 + \sum_{k=1}^{3}\Big((\wp_k - \gamma_0)s_k^2 + 2\lambda_ks_k\Big)\,. \tag{5.51}$$

Adding $\gamma_0\hat{\Omega}_1^{\mathcal{B}C_1}/\theta_1'(0)^2$ yields $\mathbf{C}_2$ from (5.21), as stated in (5.33). ■

# 6 Three-dimensional Sklyanin algebra

The three-dimensional Sklyanin algebra is defined as the quotient

$$Skl_3 = Skl_3(a, b, c) = \mathbb{C}\langle x, y, z\rangle/(r_x, r_y, r_z)\,, \tag{6.1}$$

where

$$r_x = ayz + bzy + cx^2\,, \tag{6.2}$$

$$r_y = azx + bxz + cy^2\,, \tag{6.3}$$

$$r_z = axy + byx + cz^2\,. \tag{6.4}$$

Set

$$\alpha = a(b^3 - c^3), \quad \beta = b(c^3 - a^3), \quad \gamma = c(a^3 - b^3)\,, \tag{6.5}$$

and define

$$G = \alpha(xyz + yzx + zxy) + \beta(yxz + xzy + zyx) + \gamma(x^3 + y^3 + z^3)\,. \tag{6.6}$$

It was shown in [7] that, for generic parameters, $G$ is a nonzero central element of degree three, canonical up to a nonzero scalar factor. For this algebra, there are analogues of Lemma 2.1 and Theorem 2.1 of Section 2. For a degree-one element $w = w_x x + w_y y + w_z z$, introduce three cubic expressions

$$\mathcal{E}_1(w) := \beta ywz + \alpha zwy + \gamma xwx\,, \tag{6.7}$$

$$\mathcal{E}_2(w) := \beta zwx + \alpha xwz + \gamma ywy\,, \tag{6.8}$$

$$\mathcal{E}_3(w) := \beta xwy + \alpha ywx + \gamma zwz\,. \tag{6.9}$$

Let us remark that from viewpoint of Hesse cubic, the relation between algebras $Skl_3(a,b,c)$ and $Skl_3(\beta,\alpha,\gamma)$ was studied in [7].

We can now state the following lemma.

**Lemma 6.1** *In the algebra (6.1), the values of the expressions (6.7)–(6.9) on the generators are given in the following table:*

| $w$ | $\mathcal{E}_1(w)$ | $\mathcal{E}_2(w)$ | $\mathcal{E}_3(w)$ |
|---|---|---|---|
| $x$ | $\frac{1}{3}G$ | 0 | 0 |
| $y$ | 0 | $\frac{1}{3}G$ | 0 |
| $z$ | 0 | 0 | $\frac{1}{3}G$ |

*Consequently, in the central quotient $Skl_3/(G)$, all $\mathcal{E}_i(w)$ vanish.*

**Proof.** One can check by direct computation that

$$\begin{aligned} 3\mathcal{E}_1(x) - G = (c^3 - b^3)xr_x + 2(c^3 - a^3)yr_y + (b^3 - a^3)zr_z + \\ +(2a^3 - b^3 - c^3)r_x x + (a^3 + b^3 - 2c^3)r_y y\,. \end{aligned} \tag{6.10}$$

$$\mathcal{E}_1(y) = -b^2cxr_z - a^2byr_x - ac^2zr_y + ab^2r_xy + bc^2r_yz + a^2cr_zx\,, \tag{6.11}$$

$$\mathcal{E}_1(z) = a^2cxr_y + bc^2yr_z + ab^2zr_x - a^2br_xz - b^2cr_yx - ac^2r_zy\,. \tag{6.12}$$

Thus, the right-hand sides of (6.10)–(6.12) vanish in the algebra $Skl_3$. To find the values of $\mathcal{E}_2(w)$ and $\mathcal{E}_3(w)$, note that the cyclic permutation of $(x,y,z)$ permutes the relations (6.2)–(6.4) and likewise the expressions (6.7)–(6.9). The central element (6.6) is invariant under this permutation. Therefore,

$$\mathcal{E}_1(x) = \mathcal{E}_2(y) = \mathcal{E}_3(z) = \frac{1}{3}G, \quad \mathcal{E}_1(y) = \mathcal{E}_2(z) = \mathcal{E}_3(x) = 0, \quad \mathcal{E}_1(z) = \mathcal{E}_2(x) = \mathcal{E}_3(y) = 0 \tag{6.13}$$

modulo the relations (6.2)–(6.4). ■

As in the case of the four-dimensional Sklyanin algebra in Section 2, we wish to describe the relations in the second graded component of $Skl_3$. To this end, define the matrix of quadratic monomials

$$Q = \begin{pmatrix} x \\ y \\ z \end{pmatrix} \begin{pmatrix} x & y & z \end{pmatrix} = \begin{pmatrix} x^2 & xy & xz \\ yx & y^2 & yz \\ zx & zy & z^2 \end{pmatrix}\,. \tag{6.14}$$

For $w = w_x x + w_y y + w_z z$, introduce

$$\begin{aligned} \big(w_x, w_y, w_z\big)\, Q &= \big(wx, wy, wz\big) := \big(L_x(w), L_y(w), L_z(w)\big) \;, \\ Q\,\big(w_x, w_y, w_z\big)^{\mathbf{T}} &= \big(xw, yw, zw\big)^{\mathbf{T}} := \big(R_x(w), R_y(w), R_z(w)\big)^{\mathbf{T}} \;. \end{aligned} \tag{6.15}$$

**Theorem 6.1** *In the central quotient $Skl_3/G$, the quadratic elements $L_x(u)$, $L_y(u)$, $L_z(u)$ and $R_x(u)$, $R_y(u)$, $R_z(u)$ satisfy homogeneous quadratic relations*

$$\begin{aligned} \beta L_y(u)L_z(v) + \alpha L_z(u)L_y(v) + \gamma L_x(u)L_x(v) &= 0\,, \\ \beta L_z(u)L_x(v) + \alpha L_x(u)L_z(v) + \gamma L_y(u)L_y(v) &= 0\,, \\ \beta L_x(u)L_y(v) + \alpha L_y(u)L_x(v) + \gamma L_z(u)L_z(v) &= 0\,. \end{aligned} \tag{6.16}$$

$$\begin{aligned} \beta R_y(u)R_z(v) + \alpha R_z(u)R_y(v) + \gamma R_x(u)R_x(v) &= 0\,, \\ \beta R_z(u)R_x(v) + \alpha R_x(u)R_z(v) + \gamma R_y(u)R_y(v) &= 0\,, \\ \beta R_x(u)R_y(v) + \alpha R_y(u)R_x(v) + \gamma R_z(u)R_z(v) &= 0\,. \end{aligned} \tag{6.17}$$

*In particular, setting $u = v$, one also has*

$$\begin{aligned} \beta L_y(u)L_z(u) + \alpha L_z(u)L_y(u) + \gamma L_x(u)L_x(u) &= 0\,, \\ \beta L_z(u)L_x(u) + \alpha L_x(u)L_z(u) + \gamma L_y(u)L_y(u) &= 0\,, \\ \beta L_x(u)L_y(u) + \alpha L_y(u)L_x(u) + \gamma L_z(u)L_z(u) &= 0\,. \end{aligned} \tag{6.18}$$

$$\begin{aligned} \beta R_y(u)R_z(u) + \alpha R_z(u)R_y(u) + \gamma R_x(u)R_x(u) &= 0\,, \\ \beta R_z(u)R_x(u) + \alpha R_x(u)R_z(u) + \gamma R_y(u)R_y(u) &= 0\,, \\ \beta R_x(u)R_y(u) + \alpha R_y(u)R_x(u) + \gamma R_z(u)R_z(u) &= 0\,, \end{aligned} \tag{6.19}$$

*where the constants $\alpha$, $\beta$ and $\gamma$ are defined in (6.5).*

**Proof.** The proof is similar to that of Theorem 2.1. By the preceding lemma and linearity, $\mathcal{E}_1(w) = w_x G/3$, $\mathcal{E}_2(w) = w_y G/3$, and $\mathcal{E}_3(w) = w_z G/3$. The three left-hand sides of (6.16) are $u\mathcal{E}_i(v)$, and those of (6.17) are $\mathcal{E}_i(u)v$. All six expressions vanish in $Skl_3/(G)$. Setting $v = u$ gives (6.18) and (6.19). ∎

## Main notation

$[A, B]_\pm = AB \pm BA$ – anticommutator (+) and commutator (−).

$\mathcal{A}$, $S_\alpha$ – the four-dimensional Sklyanin algebra and its generators ($\alpha = 0, 1, 2, 3$). In the defining relations (1.1), $(i, j, k)$ runs over the cyclic permutations of $(1, 2, 3)$.

$\mathcal{A}_n$ – the homogeneous component of degree $n$ in (1.2). In particular, $\mathcal{A}_1$ is spanned by the $S_\alpha$ and $\mathcal{A}_2$ by the products $S_\alpha S_\beta$.

$\mathcal{BC}_1$, $\mathbb{S}_\alpha$ – the three-parameter deformation of $\mathcal{A}$ and its generators (1.3).

$U = \sum_{\alpha=0}^{3} u_\alpha S_\alpha$, $Q_{\alpha\beta} = S_\alpha S_\beta$ – a general element of $\mathcal{A}_1$ and the entries of the matrix of quadratic monomials (2.15)–(2.16).

$X_\alpha(U) = US_\alpha$, $Y_\alpha(U) = S_\alpha U$ – the two families of quadratic elements in (1.20).

$\mathcal{F}_{ij}$, $\mathcal{F}_{0k}$ – the six quadratic polynomials defining $\mathcal{A}$, with $(i, j, k)$ cyclic in $(1, 2, 3)$ (2.22).

$\mathcal{A}(\mathcal{A}_2)$, $\mathcal{Z}_\alpha$ – the auxiliary quadratic algebra and its abstract generators (2.51).

$Z_\alpha$ – the normalized quadratic elements (2.49) that realize $\mathcal{Z}_\alpha$ in the central quotient $\mathcal{A}/(\Omega_1 - \delta_1, \Omega_2 - \delta_2)$.

$J_i$, $J_{ij} = (J_j - J_i)/J_k$ – parameters and structure constants of $\mathcal{A}$ (2.21). In the alternative parametrization, $J_i = \mathcal{K}_i(\mathcal{K}_i - \mathcal{K}_j - \mathcal{K}_k)$ (2.7).

$\Delta_i = J_j + J_k - J_i$, $K_i = \Delta_i^2$ – parameters of $\mathcal{A}(\mathcal{A}_2)$ and its $\mathcal{BC}_1$ specialization (3.25).

$\mathbf{K}_i(\hbar)$ – the elliptic functions (3.4) for which $\mathcal{K}_i(\hbar) \propto \mathbf{K}_i(\hbar)$ and $K_i(\hbar) \propto \mathbf{K}_i(2\hbar)$ with factors independent of $i$ (3.10), (3.27).

$\rho_i$, $\mu_i$, $\lambda_i$ – coefficients in the quadratic relations (2.51). Their elliptic parametrization is given in (3.26).

$\Omega_r$, $\Omega_r^{\mathcal{Z}}$, $\Omega_r^{\mathcal{BC}_1}$ ($r = 1, 2$) – Casimir elements of $\mathcal{A}$ (2.3), $\mathcal{A}(\mathcal{A}_2)$ (2.52), (2.54), and $\mathcal{BC}_1$ (2.14), respectively.

$\hat{A} = \pi(A)$ – the difference operator representing the abstract element $A$ (Section 3).

$\delta_r$, $\delta_r^{\mathcal{Z}}$, $\delta_r^{\mathcal{BC}_1}$ ($r = 1, 2$) – scalar values of $\Omega_r$, $\Omega_r^{\mathcal{Z}}$, and $\Omega_r^{\mathcal{BC}_1}$ in a representation, e.g. $\hat{\Omega}_r = \delta_r\,\mathrm{Id}$.

$\mathcal{MS}$, $\mathcal{S}_\alpha$ – the Mikhailov-Skrypnyk algebra and its generators (2.69). The affine change from $\mathcal{Z}_\alpha$ to $\mathcal{S}_\alpha$ is given in (2.72).

$\mathcal{U}_q(sl_2)$, $\mathcal{AW}$ – the quantum algebra (1.16) and the Askey-Wilson-type algebra (1.18) in the trigonometric construction.

$z$, $\hbar$, $e^{\hbar\partial_z}$ – the variable, shift parameter, and translation operator in Sklyanin's representation, $e^{\hbar\partial_z} f(z) = f(z + \hbar)$ (3.3).

$\hat{S}_\alpha$, $\hat{\mathbb{S}}_\alpha^\xi$ – elliptic difference operators representing the generators of $\mathcal{A}$ (3.1) and $\mathcal{BC}_1$ (3.82)–(3.83), respectively.

$\hat{Z}_\alpha^L$, $\hat{Z}_\alpha^R$ – the normalized products $\hat{U}\hat{S}_\alpha$ and $\hat{S}_\alpha\hat{U}$, respectively (3.17)–(3.18).

$\hat{Z}_\alpha^\xi$ – the deformed generators (3.30), which reduce to $\hat{Z}_\alpha^R$ at $\xi = \eta$.

$\overset{\sim}{S}{}^{1}_\alpha$, $\overset{\sim}{S}{}^{2}_\alpha$ – generators of the two dual copies of $\mathcal{A}$ in the elliptic modular triple (1.33). Their difference-operator realizations are given in (1.29).

$\varepsilon_1$, $\varepsilon_2$, $\varepsilon$ – the sign matrices (1.34) in the cross-relations of the modular triple. The integers $k, N$ enter through the factors $(-1)^k, (-1)^N$ in (1.33).

$\hat{\tilde{Z}}{}_\alpha^{1,\xi_1}$, $\hat{\tilde{Z}}{}_\alpha^{2,\xi_2/(2\hbar)}$ – the two modular partners of $\hat{Z}_\alpha^\xi$ in (4.37) and (4.38).

$\hat{H}^{\mathrm{RS}} = \hat{S}_0$ – the one-body quantum Ruijsenaars-Schneider Hamiltonian (1.6).

$\hat{H}^{8\mathrm{vD}}$ – the one-body Ruijsenaars-van Diejen Hamiltonian with eight independent parameters (5.1).

$\hat{H}^{7\mathrm{vD}}$, $\Delta_\nu^{\mathrm{R}}$ – the specialization $\xi = \eta$ and the Rains difference operator in its factorization $\hat{H}^{7\mathrm{vD}} = \Delta_{\bar{\nu}}^{\mathrm{R}}\Delta_\nu^{\mathrm{R}}$ (1.9).

$\hat{H}_1^{\rm 8vD}$, $\hat{H}_2^{\rm 8vD}$ – the modular partners (5.13)–(5.14) with shifts $\tau$ and 1 in $z$, respectively. Their commutation conditions are given in (5.15).

$p$, $z$ – the canonical momentum and coordinate, $\{p,z\}=1$, in the classical four-parameter Hamiltonian $H^{\rm 4vD}$ (5.16)–(5.17).

$s_0$, $s_i$, $\mathbf{C}_1$, $\mathbf{C}_2$ – the Hamiltonian, angular momentum components, and two Casimir functions of the Zhukovsky-Volterra gyrostat (5.19)–(5.21).

$Skl_3(a,b,c)$, $G$ – the three-dimensional Sklyanin algebra on $x,y,z$ (6.1) and its central cubic (6.6) for generic $a,b,c$.

$L_t(w)=wt$, $R_t(w)=tw$ $(t=x,y,z)$ – quadratic elements (6.15) for $w$ of degree one. In $Skl_3/(G)$ the transformed parameters are $(\beta,\alpha,\gamma)$, with $\alpha,\beta,\gamma$ given by (6.5).

$\tau$, $\theta_a(z|\tau)$ $(a=1,2,3,4)$ – the elliptic modulus, $\mathrm{Im}(\tau)>0$, and the Jacobi theta functions (A.1).

$\omega_0,\omega_1,\omega_2,\omega_3$ – the half-periods $0,\frac{1}{2},\frac{1+\tau}{2},\frac{\tau}{2}$, in this order (A.6).

$\phi(z,u)$, $\varphi_\alpha(z,u+\omega_\alpha)$ – the Kronecker function (A.3) and its half-period shifts with the exponential factors in (A.25). We write $\varphi_i(z)=\varphi_i(z,\omega_i)$ for $i=1,2,3$ (A.28).

$E_1$, $E_2$, $\zeta$, $\wp$ – the Eisenstein and Weierstrass functions in (A.5). We write $\wp_i=e_i=\wp(\omega_i)$ and $\wp_i^\hbar=\wp(\hbar+\omega_i)-\wp(\hbar)$ (A.29), (3.4)–(3.5).

$v(z,u)$, $\bar{v}(z,u)$ – linear combinations of $\varphi_a(2z,u+\omega_a)$ with coefficients $\nu_a$ and $\bar{\nu}_a$, respectively (1.8).

$R_{12}^\hbar(z)$, $P_{12}$ – Baxter's eight-vertex $R$-matrix with spectral parameter $z$ and the permutation operator (B.1). Here $\sigma_0=I_2$ and $\sigma_1,\sigma_2,\sigma_3$ are the Pauli matrices.

$\mathcal{L}^{\mathcal{A}}(z)$, $\mathcal{L}^{\mathcal{B}C_1}(z)$ – the Lax operators (B.3) and (B.8) in the exchange relation and the reflection equation, respectively.

Notation for coupling constants

$\eta$, $\xi$ – the spin parameter of Sklyanin's representation (3.1) and the deformation parameter in (3.30). Both enter the Hamiltonian (5.1).

$u_a$, $\bar{u}_a$ $(a=0,1,2,3)$ – coefficients of $U=\sum_a u_a S_a$ and of $\hat{H}^{\rm 8vD}=\sum_a \bar{u}_a\hat{Z}_a^\xi$, respectively (2.16), (1.35).

$\nu_a$, $\bar{\nu}_a$ $(a=0,1,2,3)$ – the two independent sets of couplings in $v$ and $\bar{v}$. Their relation to $u_a$ and $\bar{u}_a$ is given in (5.3).

$\breve{\nu}=I\nu$ – the transformed couplings obtained by interchanging the arguments of $v$, with $I^2=1$ (A.33).

$\widetilde{\nu}_i$ $(i=1,2,3)$ – the three central parameters of $\mathcal{B}C_1$, entering its defining relations through $\widetilde{\nu}_i\rho_i$ (2.11)–(2.13).

$\eta_r$, $\xi_r$, $\widetilde{\nu}^r$, $\widetilde{\bar{\nu}}^r$ $(r=1,2)$ – the spin and deformation parameters and the two coupling quadruples of the $r$-th dual Hamiltonian (5.13)–(5.14).

$c_i$, $\lambda_i=\breve{\nu}_i/c_i$ $(i=1,2,3)$ – the theta constants (A.31) and the components of the constant gyrostatic vector in the classical limit (5.24).

# 7 Appendix A: elliptic functions

**Main definitions.** We use the standard definitions of the Jacobi theta functions [25]:

$$\theta_1(u|\tau) = -\imath\sum_{k\in\mathbb{Z}}(-1)^k q^{(k+\frac{1}{2})^2}e^{\pi\imath(2k+1)u}, \qquad \theta_2(u|\tau) = \sum_{k\in\mathbb{Z}} q^{(k+\frac{1}{2})^2}e^{\pi\imath(2k+1)u}\,,$$
$$\theta_3(u|\tau) = \sum_{k\in\mathbb{Z}} q^{k^2}e^{2\pi\imath ku}, \qquad \theta_4(u|\tau) = \sum_{k\in\mathbb{Z}}(-1)^k q^{k^2}e^{2\pi\imath ku}\,, \tag{A.1}$$

where $q = e^{\pi\imath\tau}$ and $\mathrm{Im}(\tau) > 0$. The theta functions are related by shifts of the argument by half-periods:

$$\theta_2(u|\tau) = \theta_1(u+\frac{1}{2}|\tau),\;\; \theta_3(u|\tau) = q^{\frac{1}{4}}e^{\pi\imath u}\theta_1(u+\frac{\tau+1}{2}|\tau),\;\; \theta_4(u|\tau) = -\imath q^{\frac{1}{4}}e^{\pi\imath u}\theta_1(u+\frac{\tau}{2}|\tau)\,. \tag{A.2}$$

The elliptic Kronecker function is defined by

$$\phi(z,u) = \phi(z,u|\tau) = \frac{\theta_1'(0)\theta_1(z+u)}{\theta_1(z)\theta_1(u)} = \phi(u,z)\,, \quad \operatorname*{Res}_{z=0}\phi(z,u) = 1\,, \tag{A.3}$$

$$\phi(-z,-u) = -\phi(z,u)\,. \tag{A.4}$$

The first and second Eisenstein functions are defined as follows:

$$E_1(z) = \frac{\theta_1'(z)}{\theta_1(z)} = \zeta(z) + \frac{z}{3}\frac{\theta_1'''(0)}{\theta_1'(0)}\,, \quad E_2(z) = -\partial_z E_1(z) = \wp(z) - \frac{\theta_1'''(0)}{3\theta_1'(0)}\,, \tag{A.5}$$

where $\wp(z)$ and $\zeta(z)$ are the Weierstrass functions.

**Periodic properties.** Define the half-periods as

$$\omega_0 = 0\,, \quad \omega_1 = \frac{1}{2}\,, \quad \omega_2 = \frac{1+\tau}{2}\,, \quad \omega_3 = \frac{\tau}{2}\,. \tag{A.6}$$

The functions defined above satisfy the following quasi-periodicity relations:

$$\phi(z+1,u) = \phi(z,u)\,, \qquad \phi(z+\tau,u) = e^{-2\pi\imath u}\phi(z,u)\,, \tag{A.7}$$

or, equivalently,

$$\phi(z\pm 2\omega_a,u) = e^{\mp 4\pi\imath\partial_\tau\omega_a u}\phi(z,u)\,, \tag{A.8}$$

which come from

$$\theta_1(z+2\omega_a) = -e^{-4\pi\imath(z+\omega_a)\partial_\tau\omega_a}\theta_1(z)\,, \qquad \theta_1(z+\omega_a) = -e^{-4\pi\imath z\partial_\tau\omega_a}\theta_1(z-\omega_a)\,, \quad a=1,2,3\,. \tag{A.9}$$

We also use the following Weierstrass identities (see, for example, [16]):

$$\theta_1(u+x)\theta_1(u-x)\theta_r(v+y)\theta_r(v-y) - \theta_1(v+x)\theta_1(v-x)\theta_r(u+y)\theta_r(u-y) =$$
$$= \theta_1(u+v)\theta_1(u-v)\theta_r(x+y)\theta_r(x-y), \quad r=1,2,3,4\,. \tag{A.10}$$

**Duplication formulae.**

$$\theta_1(2u)\theta_2(0)\theta_3(0)\theta_4(0) = 2\theta_1(u)\theta_2(u)\theta_3(u)\theta_4(u)\,, \tag{A.11}$$

$$\theta_2(2u)\theta_2(0)\theta_3^2(0) = \theta_2^2(u)\theta_3^2(u) - \theta_1^2(u)\theta_4^2(u)\,, \tag{A.12}$$

$$\theta_2(2u)\theta_2(0)\theta_4^2(0) = \theta_2^2(u)\theta_4^2(u) - \theta_1^2(u)\theta_3^2(u)\,, \tag{A.13}$$

$$\theta_2(2u)\theta_2^3(0) = \theta_2^4(u) - \theta_1^4(u)\,, \tag{A.14}$$

$$\theta_2(2u)\theta_2^3(0) = \theta_3^4(u) - \theta_4^4(u)\,, \tag{A.15}$$

$$\theta_3(2u)\theta_3(0)\theta_2^2(0) = \theta_2^2(u)\theta_3^2(u) + \theta_1^2(u)\theta_4^2(u)\,, \tag{A.16}$$

$$\theta_3(2u)\theta_3(0)\theta_4^2(0) = \theta_3^2(u)\theta_4^2(u) - \theta_1^2(u)\theta_2^2(u)\,, \tag{A.17}$$

$$\theta_3(2u)\theta_3^3(0) = \theta_1^4(u) + \theta_3^4(u)\,, \tag{A.18}$$

$$\theta_3(2u)\theta_3^3(0) = \theta_2^4(u) + \theta_4^4(u)\,, \tag{A.19}$$

$$\theta_4(2u)\theta_4(0)\theta_2^2(0) = \theta_2^2(u)\theta_4^2(u) + \theta_1^2(u)\theta_3^2(u)\,, \tag{A.20}$$

$$\theta_4(2u)\theta_4(0)\theta_3^2(0) = \theta_1^2(u)\theta_2^2(u) + \theta_3^2(u)\theta_4^2(u)\,, \tag{A.21}$$

$$\theta_4(2u)\theta_4^3(0) = \theta_4^4(u) - \theta_1^4(u)\,, \tag{A.22}$$

$$\theta_4(2u)\theta_4^3(0) = \theta_3^4(u) - \theta_2^4(u)\,, \tag{A.23}$$

$$\phi(z,u)\phi(z,-u) = \wp(z) - \wp(u) = E_2(z) - E_2(u)\,, \tag{A.24}$$

**The functions $\varphi_k(z, x+\omega_k)$.** Using the function $\phi$ and the notation in (A.3), define the following four functions:

$$\varphi_k(z, x+\omega_k) = e^{2\pi\imath z\partial_\tau\omega_k}\phi(z, x+\omega_k)\,, \quad k = 0,1,2,3\,, \tag{A.25}$$

or, equivalently,

$$\varphi_k(z, x+\omega_k) = \frac{\theta_1'(0)\theta_{k+1}(z+x)}{\theta_1(z)\theta_{k+1}(x)}\,, \quad k = 0,1,2,3\,. \tag{A.26}$$

By (A.4) and (A.25), we have

$$\varphi_k(-z, x+\omega_k) = -\varphi_k(z, -x+\omega_k)\,, \quad k = 0,1,2,3\,. \tag{A.27}$$

**The functions $\varphi_k(z)$ and theta constants.** Here we discuss the properties of the functions (A.26) at $x=0$ for $k=1,2,3$. We write $\varphi_k(z) = \varphi_k(z,\omega_k) = e^{2\pi\imath z\partial_\tau\omega_k}\phi(z,\omega_k)$ or, equivalently,

$$\varphi_1(z) = \frac{\theta_1'(0)\theta_2(z)}{\theta_1(z)\theta_2(0)}\,, \qquad \varphi_2(z) = \frac{\theta_1'(0)\theta_3(z)}{\theta_1(z)\theta_3(0)}\,, \qquad \varphi_3(z) = \frac{\theta_1'(0)\theta_4(z)}{\theta_1(z)\theta_4(0)}\,. \tag{A.28}$$

With the half-periods numbered as in (A.6), introduce the notation $\wp_k = \wp(\omega_k)$, $k=1,2,3$:

$$\wp_1 = \wp\Big(\frac{1}{2}\Big)\,, \qquad \wp_2 = \wp\Big(\frac{1+\tau}{2}\Big)\,, \qquad \wp_3 = \wp\Big(\frac{\tau}{2}\Big)\,. \tag{A.29}$$

From (A.24) one gets

$$\varphi_k^2(z) = \wp(z) - \wp_k\,, \quad k = 1,2,3\,. \tag{A.30}$$

Introduce the following theta constants $c_k = c_k(\tau)$, $k=1,2,3$:

$$c_1(\tau) = -\Big(\frac{\theta_2(0)}{\theta_1'(0)}\Big)^2, \quad c_2(\tau) = -\imath\Big(\frac{\theta_3(0)}{\theta_1'(0)}\Big)^2, \quad c_3(\tau) = -\Big(\frac{\theta_4(0)}{\theta_1'(0)}\Big)^2. \tag{A.31}$$

**Properties of the function $v(z,u)$.** Define

$$v(z,u) \equiv v(z,u|\tau;\nu_0,\nu_1,\nu_2,\nu_3) =$$

$$= \sum_{a=0}^{3} \nu_a \varphi_a(2z, u+\omega_a) = \sum_{a=0}^{3} \nu_a \exp(4\pi\imath z\partial_\tau \omega_a)\phi(2z, u+\omega_a)\,. \tag{A.32}$$

Interchanging the arguments is equivalent to the following change of constants:

$$\breve{\nu} = I\nu, \qquad I = \frac{1}{2}\begin{pmatrix} 1 & 1 & 1 & 1 \\ 1 & 1 & -1 & -1 \\ 1 & -1 & 1 & -1 \\ 1 & -1 & -1 & 1 \end{pmatrix}, \qquad I^2 = 1\,,$$

$$v(u,z|\tau;\nu) = v(z,u|\tau;\breve{\nu})\,. \tag{A.33}$$

Below we list some identities for the function (A.32). Further properties can be found in [17, 4, 24] and in the appendix of [22]. First,

$$v(z,u)v(z,-u) = \sum_{a=0}^{3}\Big(\breve{\nu}_a^2 \wp(z+\omega_a) - \nu_a^2 \wp(u+\omega_a)\Big)\,. \tag{A.34}$$

By (A.27) we have

$$v(-z,-u) = -v(z,u)\,, \tag{A.35}$$

and, therefore,

$$v(z,u)v(-z,u) = \sum_{a=0}^{3}\Big(\nu_a^2 \wp(u+\omega_a) - \breve{\nu}_a^2 \wp(z+\omega_a)\Big) = -v(z,u)v(z,-u)\,. \tag{A.36}$$

Further elliptic identities for these functions can be found, for example, in [24].

# 8 Appendix B: $R$-matrix derivation of the Sklyanin algebras

Here we briefly recall (in our notation) how the original Sklyanin algebra $\mathcal{A}$ (2.1)-(2.2) follows from the quantum exchange relations [37], and how the $\mathcal{B}C_1$ Sklyanin algebra (2.11)-(2.13) follows from the quantum reflection equation [39, 19].

**The Baxter $R$-matrix.** Let $\sigma_0 = 1_{2\times 2}$ and $\sigma_i$ be the Pauli matrices, with $\sigma_i^2 = 1_{2\times 2}$ and $\sigma_i\sigma_j = \imath\sigma_k = -\sigma_j\sigma_i$. The elliptic $R$-matrix is [3]

$$R_{12}^{\hbar}(z) = \sum_{a=0}^{3} \varphi_a(z, \hbar+\omega_a)\,\sigma_a \otimes \sigma_a\,,$$

$$\operatorname*{Res}_{z=0} R_{12}^{\hbar}(z) = 2P_{12}\,, \qquad P_{12} = \frac{1}{2}\sum_{a=0}^{3} \sigma_a \otimes \sigma_a\,. \tag{B.1}$$

It satisfies the quantum Yang–Baxter equation. The main scalar identity used in the component calculations below is the Fay identity

$$\phi(x,u)\phi(y,v) = \phi(x-y,u)\phi(y,u+v) + \phi(y-x,v)\phi(x,u+v)\,. \tag{B.2}$$

**The original Sklyanin algebra.** Introduce the Lax operator

$$\mathcal{L}^{\mathcal{A}}(z) = \sum_{a=0}^{3} S_a\, \varphi_a(z, \hbar + \omega_a)\, \sigma_a\,, \qquad \mathcal{L}_1 = \mathcal{L} \otimes I_2\,, \quad \mathcal{L}_2 = I_2 \otimes \mathcal{L}\,, \tag{B.3}$$

and consider the exchange relation

$$R^{\hbar}_{12}(z-w)\mathcal{L}^{\mathcal{A}}_1(z)\mathcal{L}^{\mathcal{A}}_2(w) = \mathcal{L}^{\mathcal{A}}_2(w)\mathcal{L}^{\mathcal{A}}_1(z)R^{\hbar}_{12}(z-w)\,. \tag{B.4}$$

Then (B.4) is equivalent to

$$\begin{gathered} [S_i, S_j]_- = \imath [S_0, S_k]_+\,, \\ [S_0, S_k]_- = \imath \frac{\wp^{\hbar}_j - \wp^{\hbar}_i}{\wp^{\hbar}_k}\, [S_i, S_j]_+\,. \end{gathered} \tag{B.5}$$

Finally, (3.9) gives

$$\frac{\wp^{\hbar}_j - \wp^{\hbar}_i}{\wp^{\hbar}_k} = \frac{J_j(\hbar) - J_i(\hbar)}{J_k(\hbar)}\,, \qquad J_i(\hbar) = -\frac{\theta_1(\hbar)^2}{\theta_1'(0)^2}\, \wp^{\hbar}_i\,. \tag{B.6}$$

Thus the component equations are precisely (2.1)–(2.2) with the theta-function parametrization (3.2).

Higher dimensional generalization of the Sklyanin algebra can be found in [28, 29, 30] and [5].

**The $BC_1$ Sklyanin algebra.** In accordance with (3.27), set

$$\begin{gathered} K_i(\hbar) = c(\hbar)\, \mathbf{K}_i(2\hbar)\,, \qquad c(\hbar) = 8\, \frac{\theta_1(\hbar)^4}{\theta_1'(0)^4}\, \wp'(2\hbar)\,, \\ \rho_i = J_1(\hbar)J_2(\hbar)J_3(\hbar) \frac{\theta_{j+1}(0)\theta_{k+1}(0)}{\theta_{j+1}(\hbar)\theta_{k+1}(\hbar)}\,. \end{gathered} \tag{B.7}$$

Here $K_i = \Delta_i^2$ are as in (3.25), and the $\widetilde{\nu}_i$ are free parameters. The normalization of $\rho_i$ agrees with (5.29). Define the Lax operator:

$$\mathcal{L}^{\mathcal{B}C_1}(z) = \mathbb{S}_0\, \varphi_0(z, 2\hbar)\, \sigma_0 + \sum_{i=1}^{3} \left( \mathbb{S}_i\, \varphi_i(z, 2\hbar + \omega_i) + \frac{\widetilde{\nu}_i \rho_i}{c(\hbar)}\, \frac{\varphi_i(z - \omega_i, 2\hbar + \omega_i)}{\varphi_i(-\omega_i, 2\hbar + \omega_i)} \right) \sigma_i\,. \tag{B.8}$$

Then the reflection equation

$$R^{\hbar}_{12}(z-w)\mathcal{L}^{\mathcal{B}C_1}_1(z)R^{\hbar}_{12}(z+w)\mathcal{L}^{\mathcal{B}C_1}_2(w) = \mathcal{L}^{\mathcal{B}C_1}_2(w)R^{\hbar}_{12}(z+w)\mathcal{L}^{\mathcal{B}C_1}_1(z)R^{\hbar}_{12}(z-w) \tag{B.9}$$

is equivalent to

$$\begin{gathered} [\mathbb{S}_i, \mathbb{S}_j]_- = \imath [\mathbb{S}_0, \mathbb{S}_k]_+\,, \\ [\mathbb{S}_0, \mathbb{S}_k]_- = -\imath \frac{K_j - K_i}{K_k}\, [\mathbb{S}_i, \mathbb{S}_j]_+ + \frac{2\imath}{K_k}\, (\widetilde{\nu}_i \rho_i \mathbb{S}_j - \widetilde{\nu}_j \rho_j \mathbb{S}_i)\,, \end{gathered} \tag{B.10}$$

which gives (2.11)–(2.13). When $\widetilde{\nu}_i = 0$, the resulting algebra is the original Sklyanin algebra with the argument $\hbar$ doubled, since

$$-\frac{K_j(\hbar) - K_i(\hbar)}{K_k(\hbar)} = \frac{J_j(2\hbar) - J_i(2\hbar)}{J_k(2\hbar)}\,. \tag{B.11}$$

### Acknowledgments

We are grateful to S. Derkachev, A. Mikhailov, V. Rubtsov, V. Spiridonov, A. Zabrodin and especially N. Zamolotov for useful discussions and comments.

This work was supported by the Russian Science Foundation under grant no. 25-11-00081, https://rscf.ru/en/project/25-11-00081/, and was carried out at the Steklov Mathematical Institute of the Russian Academy of Sciences.

During the preparation of this work, the authors used ChatGPT for language, grammar and editorial refinement, typographical and presentational corrections, LaTeX assistance, and supplementary checking for possible errors and gaps. After using this tool, the authors reviewed and edited the content as needed and take full responsibility for the content of the published article.